\documentclass{article}

\usepackage{arxiv}
\usepackage[onehalfspacing]{setspace}

\usepackage[utf8]{inputenc} % allow utf-8 input
\usepackage[T1]{fontenc}    % use 8-bit T1 fonts
\usepackage{hyperref}       % hyperlinks
\usepackage{url}            % simple URL typesetting
\usepackage{booktabs}       % professional-quality tables
\usepackage{amsmath,amsfonts,amssymb}       % blackboard math symbols
\usepackage{nicefrac}       % compact symbols for 1/2, etc.
\usepackage{microtype}      % microtypography
\usepackage{cleveref}       % smart cross-referencing
\usepackage{lipsum}         % Can be removed after putting your text content
\usepackage[square,numbers]{natbib}
\usepackage{doi}
\usepackage{xcolor}
\usepackage{bm}
\usepackage{derivative}
\usepackage{caption}
\usepackage{cancel}
\usepackage{comment}
\usepackage{algorithm}
\usepackage{float}
\usepackage{tikz}
\usetikzlibrary{positioning, shapes, decorations.markings, decorations.pathreplacing, decorations.pathmorphing, arrows, calc, matrix, fit, backgrounds}
\usepackage{subcaption}
\usepackage{mathtools}
\usepackage{makecell}
\usepackage[graphicx]{realboxes}
\usepackage{pdflscape}
\usepackage[bottom]{footmisc}
\usepackage{enumitem}

\DeclarePairedDelimiter{\ceil}{\lceil}{\rceil}
\DeclarePairedDelimiter{\floor}{\lfloor}{\rfloor}

\newcommand{\rmax}{r_{\mathrm{max}}}

\newcommand{\frulr}{f_{R U L_{k}}\left(r\right)}
\newcommand{\EC}{\mathbb{E}\left[C\right]}
\newcommand{\ECrep}{\mathbb{E}\left[C_{\mathrm{rep}}\right]}
\newcommand{\ECord}{\mathbb{E}\left[C_{\mathrm{ord}}\right]}
\newcommand{\ETlc}{\mathbb{E}\left[T_{\mathrm{lc}}\right]}
\newcommand{\Ectratio}{\frac{\EC}{\ETlc}}
\newcommand{\Ectrep}{\frac{\ECrep}{\ETlc}}
\newcommand{\Ectorder}{\frac{\ECord}{\ETlc}}

\newcommand{\Pf}{p_{F}}
\newcommand{\Pfk}{p_{F,k}}
\newcommand{\Tf}{T_F}

\newcommand{\dt}{\mathrm{d} t}
\newcommand{\dr}{\mathrm{d} r}
\newcommand{\EE}[1]{\mathbb{E}\left[#1\right]}
\newcommand{\trep}{t_{\mathrm{rep}}}
\newcommand{\rrep}{r_{\mathrm{rep}}}
\newcommand{\tord}{t_{\mathrm{ord}}}
\newcommand{\cinv}{c_{\mathrm{inv}}}
\newcommand{\cunav}{c_{\mathrm{unav}}}
\newcommand{\Pro}[1]{\mathrm{Pr}\left(#1\right)}
\newcommand{\mutf}{\mu_{T_{{F}}}}
\newcommand{\Tfi}{T_{F}^{(i)}}
\newcommand{\Pfrrepk}{\Pfk(r_{\mathrm{rep}})}
\newcommand{\Tinvstari}{T_\mathrm{inv}^{(i)^*}}

\title{Leaf It to Renewal: Improved Predictive Maintenance Policies via Renewal Theory and Decision Trees}

\author{Daniel Koutas \\
	ERA Group\\
	TU Munich\\
	Munich, Germany \\
	\texttt{daniel.koutas@tum.de} \\
	\And
	Daniel Straub \\
	ERA Group\\
	TU Munich\\
	Munich, Germany \\
	\texttt{straub@tum.de} \\
}

\renewcommand{\shorttitle}{Leaf It to Renewal}

\newcommand{\settinggraph}{
\begin{figure}[H]
    \def\lx{6.6}
    \def\rx{6.8}
    \def\l{1}
    \def\u{3.2}
    \def\plx{8.5}
    \def\pl{4}
    \def\pw{6.3}
    \def\ph{3}
    \def\boxc{blue!90!black}
    \def\hshift{1.5}
    \def\radius{1.3}
    \def\boundline{thick}
    \def\dnyshift{1.2}
    \def\dnxshift{1.75}

    \tikzset{circle_node/.style={circle, draw, \boundline, minimum size=\radius cm, inner sep=0pt}}
    \centering
    \begin{tikzpicture}
    \node[anchor=south west, inner sep=0] (image) at (0,0) {\includegraphics[width=0.5\textwidth]{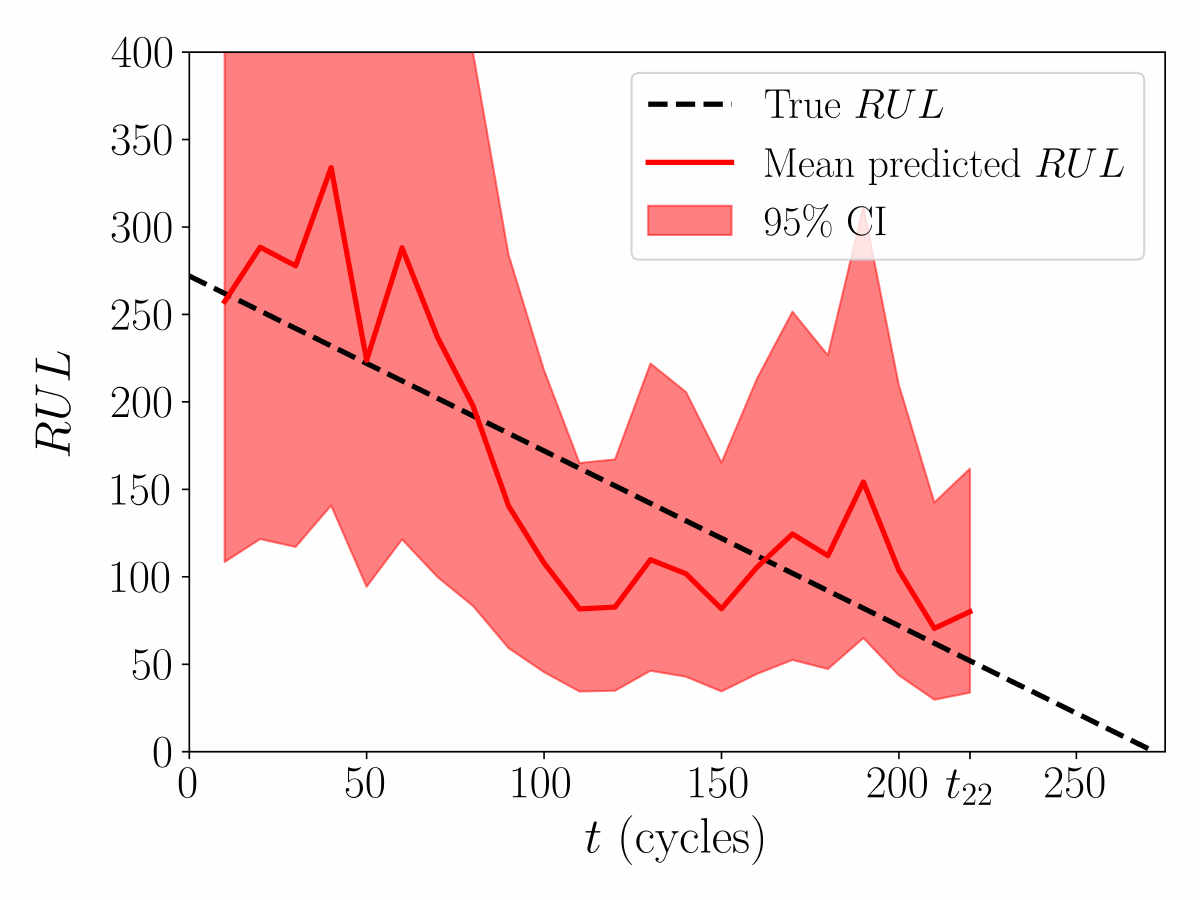}};
        
    \begin{scope}
        \draw[\boxc, thick] ($(\lx,\l)$) rectangle ($(\rx,\u)$); % small rectangle on image
        \draw[-latex, thick] ($(\rx,\u)$) -- ($(\plx,\pl)$);
    \end{scope}
    
    \begin{scope}[xshift=\plx cm, yshift=\pl cm]
        \clip (0,0) rectangle (\pw cm, \ph cm); % clip the magnified area
        \node[anchor=south west, inner sep=-6pt] (pdf) at (0,0) {\includegraphics[width=0.4\textwidth]{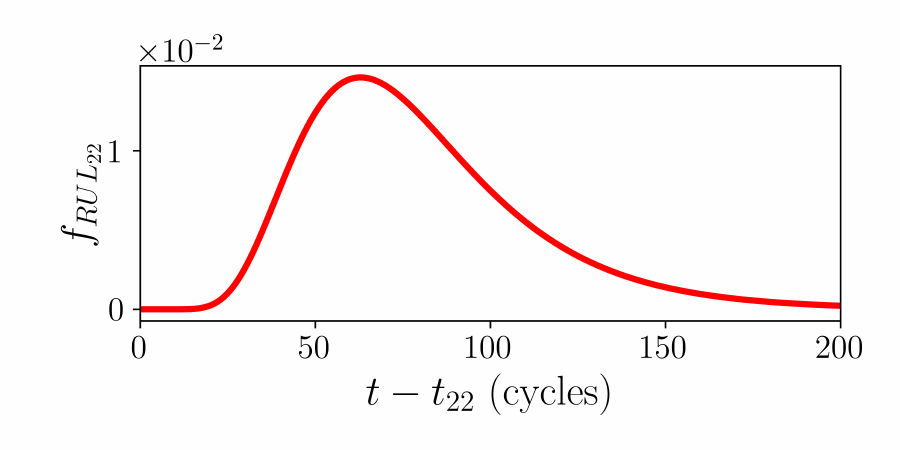}};
        \draw[thick, \boxc] (0,0) rectangle (\pw cm, \ph cm); % border around magnified part
    \end{scope}

    \coordinate[below=\hshift cm of pdf] (headcord) {};
    \node[circle_node] (head) at (headcord) {\Large $\Pi$};
    \draw[-latex, thick] (pdf) -- (head);

    \coordinate[below=\dnyshift cm of head, xshift=-\dnxshift cm] (dncord) {};
    \node[draw, ellipse] (dn) at (dncord) {\large Do nothing};
    \draw[-latex, thick] (head) -- (dn);

    \coordinate[below=\dnyshift cm of head, xshift=\dnxshift cm] (mcord) {};
    \node[draw, ellipse, align=center] (m) at (mcord) {\large Maintenance \\ \large action};
    \draw[-latex, thick] (head) -- (m);
    \end{tikzpicture}
    \caption{General workflow of maintenance decision policies. At every time step, the latest RUL-PDF is fed into a heuristic policy ($\Pi$), which outputs a decision to either do nothing or to perform a maintenance action, e.g., replacement or ordering a spare component. The workflow is exemplarily shown for the time step $t_k=220$.}
    \label{fig:decision_setting}
\end{figure}
}

\newcommand{\repairdecisiontrees}{
\begin{figure}[H]
\def\xdist{2.0}
\def\sydist{1.0}
\def\aydist{1.0}
\def\xtextdist{1.2}
\def\ytextdist{1.0}
\def\tdist{3.8}
\def\width{0.4}
\def\radius{0.3}
\def\termradius{0.2}
\def\capsize{\normalsize}
\def\boundline{very thick}
\def\arrowthick{very thick}
\def\tshift{7}
\def\xenumshift{0.6}
\def\yenumshift{0.25}
\tikzset{circle_node/.style={circle, draw, \boundline, minimum size=\radius cm, inner sep=0pt}}
\tikzset{term_node/.style={circle, draw, fill=black, \boundline, minimum size=\termradius cm, inner sep=0pt}}
\tikzset{square_node/.style={rectangle, draw, fill=cyan, \boundline, minimum size=\width cm, inner sep=0pt}}

    \centering
    \begin{tikzpicture}
    \coordinate (a0cord) at (0,0);
    \node[square_node] (a0) at (a0cord) {};
    \coordinate[right=\xdist cm of a0cord, yshift=\sydist cm] (s1dncord) {};
    \node[circle_node] (s1dn) at (s1dncord) {};
    \coordinate[right=\xdist cm of a0cord, yshift=-\sydist cm] (s1prcord) {};
    \node[term_node] (s1pr) at (s1prcord) {};
    \coordinate[right=\xdist cm of s1dncord, yshift=\aydist cm] (a1cord) {};
    \node[square_node] (a1) at (a1cord) {};
    \coordinate[right=\xdist cm of s1dncord, yshift=-\aydist cm] (s1failedcord) {};
    \node[term_node] (s1failed) at (s1failedcord) {};
    \coordinate[right=\xdist cm of a1cord, yshift=\sydist cm] (s2dncord) {};
    \node[circle_node, dashed] (s2dn) at (s2dncord) {};
    \coordinate[right=\xdist cm of a1cord, yshift=-\sydist cm] (s2prcord) {};
    \node[term_node] (s2pr) at (s2prcord) {};

    \coordinate[right=\xtextdist cm of s2prcord] (c2prcord) {};
    \node[anchor=west] (c2pr) at (c2prcord) {$c_p + c_\infty \cdot (r - \Delta t)$};
    \coordinate (c1failedcord) at ($(s1failedcord) + (\xdist + \xtextdist, 0)$) {};
    \node[anchor=west] (c1failed) at (c1failedcord) {$c_c + c_\infty \cdot (r - \EE{RUL \mid RUL < \Delta t})$};
    \coordinate (c1prcord) at ($(s1prcord) + (2*\xdist + \xtextdist, 0)$) {};
    \node[anchor=west] (c1pr) at (c1prcord) {$c_p + c_\infty \cdot r$};

    \coordinate[left=\xdist cm of a0cord, yshift=-\aydist cm] (sprevcord) {};

    \draw[\arrowthick] (a0) -- (s1dn) node[midway, yshift=\tshift pt, rotate=30] {DN};
    \draw[\arrowthick] (a0) -- (s1pr) node[midway, yshift=\tshift pt, rotate=-30] {PR};
    \draw[\arrowthick] (s1dn) -- (a1) node[midway, yshift=\tshift pt, rotate=30] {\textcolor{green!60!black}{Safe}};
    \node[xshift=6pt, yshift=-5 pt, rotate=30] at ($(s1dncord)!0.5!(a1cord)$) {$1-\Pf$};
    \draw[\arrowthick] (s1dn) -- (s1failed) node[midway, yshift=\tshift pt, rotate=-30] {\textcolor{red!80!black}{Failed}};
    \node[xshift=-5pt, yshift=-\tshift pt, rotate=-30] at ($(s1dncord)!0.5!(s1failedcord)$) {$\Pf$};
    \draw[dashed, \arrowthick] (a1) -- (s2dn) node[midway, yshift=\tshift pt, rotate=30] {DN};
    \draw[\arrowthick] (a1) -- (s2pr)node[midway, yshift=\tshift pt, rotate=-30] {PR};
    \draw[\arrowthick] (s2pr) -- (c2pr);
    \draw[\arrowthick] (s1failed) -- (c1failed);
    \draw[\arrowthick] (s1pr) -- (c1pr);
    \begin{scope}
        \clip (a0) circle (1cm);
        \draw[dashed, \arrowthick] (a0) -- (sprevcord);
    \end{scope}

    \coordinate[right=\xtextdist cm of s2dncord] (c2dncord) {};
    \draw[dashed, \arrowthick] (s2dn) -- (c2dncord);
    
    \coordinate[above=\tdist cm of a0cord] (t0cord) {};
    \coordinate (t1cord) at ($(t0cord) + (2*\xdist, 0)$) {};
    \coordinate (t2cord) at ($(t1cord) + (3*\xdist, 0)$) {};
    \coordinate (tprevcord) at ($(t0cord) + (-3*\xdist, 0)$) {};

    \draw[decorate,decoration={brace,amplitude=10pt}] (t0cord) -- (t1cord) 
        node[midway,yshift=20pt] {\large Step $k$};
    \begin{scope}
        \clip (t1cord) circle (2.1 cm);
        \draw[decorate,decoration={brace,amplitude=10pt}] (t1cord) -- (t2cord) node[midway, xshift=-2cm, yshift=20pt] {\large Step $k+1$};
    \end{scope}
    \begin{scope}
        \clip (t0cord) circle (2.1 cm);
        \draw[decorate,decoration={brace,amplitude=10pt}] (tprevcord) -- (t0cord) node[midway, xshift=2cm, yshift=20pt] {\large Step $k-1$};
    \end{scope}

    \draw[thick, dashed, gray] ($(t0cord) + (0, 1)$) -- (a0) -- ($(a0cord) + (0, -1.5)$){};
    \node[gray] at ($(t0cord) + (0, 1.4)$) {$\begin{array}{l}
        r=0\\
        t=t_k
    \end{array}$};
    \draw[thick, dashed, gray] ($(t1cord) + (0, 1)$) -- (a1) -- (s1failed) -- ($(s1failedcord) + (0, -1.5)$) {};
    \node[gray] at ($(t1cord) + (0, 1.4)$) {$\begin{array}{l}
        r=\Delta t\\
        t=t_{k+1}
    \end{array}$};

    \coordinate[right=\xenumshift cm of s1prcord, yshift=\yenumshift cm] (bacord) {};
    \node at (bacord) {$(a)$};
    \coordinate[right=\xenumshift cm of s1failedcord, yshift=\yenumshift cm] (bbcord) {};
    \node at (bbcord) {$(b)$};
    \coordinate[right=\xenumshift cm of s2prcord, yshift=\yenumshift cm] (bccord) {};
    \node at (bccord) {$(c)$};
    \coordinate[right=\xenumshift cm of s2dncord, 
    yshift=\yenumshift cm
    ] (bdcord) {};
    \node at (bdcord) {$(d)$};
    
    \end{tikzpicture}
    \captionsetup{width=\textwidth}
    \caption{Decision tree for a deteriorating component that can be preventively replaced. Component states are represented by round nodes, and maintenance decisions ($a_k$) as square nodes. The action taken (PR/DN) as well as the condition of the component (Failed/Safe) with the corresponding event probabilities are included in the respective paths, where $\Pf=\Pr(RUL\leq\Delta t)$. The component's end of life (filled black circle) is shown together with the associated expected costs. The Figure is adapted from \cite{bismut2021optimal}.}
    \label{fig:rep_decision_tree}
\end{figure}
}

\newcommand{\orderdecisiontrees}{
\begin{figure}[H]
\captionsetup{aboveskip=20pt}
\def\xdist{3}
\def\sydist{1.8}
\def\aydist{1.1}
\def\xtextdist{1.2}
\def\ytextdist{1.0}
\def\tdist{4.0}
\def\width{0.4}
\def\radius{0.3}
\def\termradius{0.2}
\def\capsize{\normalsize}
\def\boundline{very thick}
\def\arrowthick{very thick}
\def\tshift{7}
\def\otklen{3}
\def\xenumshift{0.6}
\def\yenumshift{0.25}
\tikzset{circle_node/.style={circle, draw, \boundline, minimum size=\radius cm, inner sep=0pt}}
\tikzset{term_node/.style={circle, draw, fill=black, \boundline, minimum size=\termradius cm, inner sep=0pt}}
\tikzset{square_node/.style={rectangle, draw, fill=cyan, \boundline, minimum size=\width cm, inner sep=0pt}}
\tikzset{ghost_node/.style={circle, draw, fill=grey, \boundline, minimum size=0.0 cm, inner sep=0pt}}

    \centering
    \begin{tikzpicture}
    \coordinate (a0cord) at (0,0);
    \node[square_node] (a0) at (a0cord) {};
    \coordinate[right=2*\xdist cm of a0cord, yshift=\sydist cm] (otkcord) {};
    \node[circle_node] (otk) at (otkcord) {};
    \coordinate[right=\xdist cm of a0cord, yshift=-\sydist cm] (olatercord) {};
    \coordinate[right=\xdist cm of olatercord] (otk1cord) {};
    \node[circle_node] (otk1) at (otk1cord) {};
    \coordinate[right=\xdist cm of otkcord, yshift=-\aydist cm] (fotkcord) {};
    \node[term_node] (fotk) at (fotkcord) {};
    \coordinate[right=\xdist cm of otkcord, yshift=\aydist cm] (nfotkcord) {};
    \node[term_node] (nfotk) at (nfotkcord) {};
    \coordinate[right=\xdist cm of otk1cord, yshift=-\aydist cm] (fotk1cord) {};
    \node[term_node] (fotk1) at (fotk1cord) {};
    \coordinate[right=\xdist cm of otk1cord, yshift=+\aydist cm] (nfotk1cord) {};
    \node[term_node] (nfotk1) at (nfotk1cord) {};

    \pgfmathsetmacro\oangtk{atan(\sydist / (2*\xdist))}
    \pgfmathsetmacro\oanglater{atan(\sydist / \xdist)}

    \pgfmathsetmacro\fang{atan(\aydist / \xdist)}

    \coordinate[right=\xtextdist cm of fotkcord] (cfotkcord) {};
    \coordinate[right=\xtextdist cm of nfotkcord] (cnfotkcord) {};
    \coordinate[right=\xtextdist cm of fotk1cord] (cfotk1cord) {};
    \coordinate[right=\xtextdist cm of nfotk1cord] (cnfotk1cord) {};

    \draw[\arrowthick] (a0) -- (otk) node[midway, yshift=\tshift pt, rotate=\oangtk] {Order at $t_k$};
    \draw[\arrowthick] (a0) -- (olatercord) node[midway, yshift=\tshift pt, rotate=-\oanglater] {Order later};
    \draw[\arrowthick] (olatercord) -- (otk1) node[midway, yshift=\tshift pt] {Order at $t_{k+1}$};
    \draw[\arrowthick] (otk) -- (nfotk) node[midway, yshift=\tshift pt, rotate=\fang] {\textcolor{green!60!black}{Safe}};
    \node[xshift=6pt, yshift=-6 pt, rotate=\fang] at ($(otkcord)!0.5!(nfotkcord)$) {$1 - \Pf^{(1)}$};
    \draw[\arrowthick] (otk) -- (fotk) node[midway, yshift=\tshift pt, rotate=-\fang] {\textcolor{red!80!black}{Failed}};
    \node[xshift=-5pt, yshift=-\tshift pt, rotate=-\fang] at ($(otkcord)!0.5!(fotkcord)$) {$\Pf^{(1)}$};
    \draw[\arrowthick] (otk1) -- (nfotk1) node[midway, yshift=\tshift pt, rotate=\fang] {\textcolor{green!60!black}{Safe}};
    \node[xshift=6pt, yshift=-6 pt, rotate=\fang] at ($(otk1cord)!0.5!(nfotk1cord)$) {$1 - \Pf^{(2)}$};
    \draw[\arrowthick] (otk1) -- (fotk1) node[midway, yshift=\tshift pt, rotate=-\fang] {\textcolor{red!80!black}{Failed}};
    \node[xshift=-5pt, yshift=-\tshift pt, rotate=-\fang] at ($(otk1cord)!0.5!(fotk1cord)$) {$\Pf^{(2)}$};
    \draw[\arrowthick] (fotkcord) --(cfotkcord);
    \draw[\arrowthick] (nfotkcord) --(cnfotkcord);
    \draw[\arrowthick] (fotk1cord) --(cfotk1cord);
    \draw[\arrowthick] (nfotk1cord) --(cnfotk1cord);

    \coordinate[right=\xdist cm of olatercord, yshift=-\sydist cm] (otk2cord) {};
    \node[circle_node, dashed] (otk2) at (otk2cord) {};
    \pgfmathsetmacro\pang{atan(\sydist / \xdist)}
    \draw[dashed, \arrowthick] (olatercord) -- (otk2) node[midway, yshift=3 pt, xshift=\tshift pt, rotate=-\pang] {Order at $t_{k+2}$};

    \coordinate[right=\xdist cm of olatercord, yshift=-2*\sydist cm] (otk3cord) {};
    \coordinate (otk3distcord) at ($(olatercord) ! 0.5 ! (otk3cord)$);
    \draw[dashed, \arrowthick] (olatercord) -- (otk3distcord);    

    \coordinate[right=\xenumshift cm of nfotkcord, yshift=\yenumshift cm] (bacord) {};
    \node at (bacord) {$(a)$};
    \coordinate[right=\xenumshift cm of fotkcord, yshift=\yenumshift cm] (bbcord) {};
    \node at (bbcord) {$(b)$};
    \coordinate[right=\xenumshift cm of nfotk1cord, yshift=\yenumshift cm] (bccord) {};
    \node at (bccord) {$(c)$};
    \coordinate[right=\xenumshift cm of fotk1cord, 
    yshift=\yenumshift cm
    ] (bdcord) {};
    \node at (bdcord) {$(d)$};
    
    \end{tikzpicture}
    \captionsetup{width=\textwidth}
    \caption{Decision tree for a deteriorating component, where a spare can be preventively ordered. Component states are represented by round nodes and the ordering decision ($a_k$) as a square node. The action taken, as well as the condition of the component (Failed/Safe) with the corresponding event probabilities are included in the respective paths, where $\Pf^{(1)} = \Pr(RUL\leq l)$ and $\Pf^{(2)}=\Pr(RUL\leq l+\Delta t)$. Here, filled black circles indicate the component's eventual end of life (the ordering decision has no influence on that). The costs for each branch are not shown in the graph to avoid cluttering; instead, they are listed in \Cref{eq:h_ord_tk_fail,eq:h_ord_tk_surv,eq:h_ord_tk+1_fail,eq:h_ord_tk+1_surv}.}
    \label{fig:order_decision_tree}
\end{figure}
}

\pgfmathdeclarefunction{lognormal}{3}{%
\pgfmathparse{%
    (exp(-(ln(#1+0.01)-#2)^2 / (2*#3^2)) / (#1*#3*sqrt(2)+0.01))%
    }%
}

\newcommand{\ordergraphs}{
\begin{figure}[H]
\def\xdist{7.4}
\def\tlabelshift{-0.5}
\def\ytickheight{0.2}
\def\xtk{0.3}
\def\tickdist{1.3}
\def\ylabelshift{1.0}
\def\ylshift{1.2}
\pgfmathsetmacro{\Ldist}{2.7*\tickdist}
\def\Lytickheight{0.1}
\def\Lylabelshift{1.0}
\def\linethick{ultra thick}
\def\tickthick{very thick}
\def\rulcolor{red}
\def\domain{1.2:7.1}
\def\mu{1.4}
\def\sigma{0.4}
\def\Tfone{2.5}
\def\Tftwo{6.6}
\def\Tfonecolor{blue!80!black}
\def\Tftwocolor{green!50!black}
\def\plotyshift{0.1}
\def\plotfactor{2.0}
\def\textxshift{0.0}
\def\textyshift{1.3}
\def\ycshift{0.8}
\def\cunavcolor{\Tfonecolor}
\def\cinvcolor{\Tftwocolor}
\def\ectrcolor{orange!90!black}
\def\colorshift{0.04}
\def\yrshift{0.9}
\def\xrshift{1.0}
\def\rthick{thick}
\def\rtickheight{0.15}
\def\rlabelshift{0.8}

    \centering
    \begin{subfigure}{0.49\textwidth}
        \centering
        \begin{tikzpicture}
            \coordinate (0) at (0,0);
            \coordinate (end) at ($(0) + (\xdist, 0)$);
            \node[yshift=\tlabelshift cm] at (end) {$t$};
            \coordinate (tk) at ($(\xtk,0)$);
            \draw[\tickthick] ($(tk) + (0,\ytickheight)$) -- ($(tk) + 
            (0,-\ytickheight)$) node[anchor=north, yshift=-\ylabelshift]{\large $t_k$};
            \coordinate (tk1) at ($(tk) + (\tickdist,0)$);
            \draw[\tickthick] ($(tk1) + (0,\ytickheight)$) -- ($(tk1) + 
            (0,-\ytickheight)$) node[anchor=north, yshift=-\ylabelshift]{\large $t_{k+1}$};

            \foreach \k in {1,2,4}{
                \draw[\tickthick] ($(tk) + (\k*\tickdist,\ytickheight)$) -- ($(tk) + (\k*\tickdist,-\ytickheight)$);
            } 
            \coordinate (tkw) at ($(tk) + (3*\tickdist,0)$);
            \draw[\tickthick] ($(tkw) + (0,\ytickheight)$) -- 
            ($(tkw) + (0,-\ytickheight)$) node[anchor=north, yshift=-\ylabelshift]{\large $t_k+w$};
    
            \coordinate (Ll) at ($(tk) + (0,-\ylshift)$);
            \coordinate (Lr) at ($(Ll) + (\Ldist,0)$);
            \draw[\tickthick] (Ll) -- (Lr) node[midway, anchor=north, yshift=-\Lylabelshift]{$l$};
            \draw[\tickthick] ($(Ll) + (0,\Lytickheight)$) -- ($(Ll) + (0,-\Lytickheight)$);
            \draw[\tickthick] ($(Lr) + (0,\Lytickheight)$) -- ($(Lr) + (0,-\Lytickheight)$);
    
            \draw[\rulcolor, very thick, domain=\domain, samples=200] 
            plot ({\x}, {\plotfactor*lognormal(\x,\mu,\sigma)+\plotyshift});
            \node[] at ($(tkw) + (\textxshift, \textyshift)$) {$\color{\rulcolor} \frulr$};
            \pgfmathsetmacro\Tfoney{\plotfactor*lognormal(\Tfone,\mu,\sigma)+\plotyshift}
            \draw[\Tfonecolor, dashed, very thick] (\Tfone,0) -- (\Tfone,\Tfoney) 
            node[very near start, yshift=-0.55cm]{\textcolor{\Tfonecolor}{$T_F^{(1)}$}};
    
            \pgfmathsetmacro{\rrend}{\xtk + \Ldist}
            \draw[\linethick] (0) -- (\Tfone,0);
            \draw[dotted, \linethick] (\Tfone,0) -- (\rrend,0);
            \draw[-latex, \linethick] (\rrend,0) -- (end);

            \coordinate (cl) at ($(\Tfone,0) + (0,-\ylshift-\ycshift)$);
            \coordinate (cr) at ($(Lr) + (0,-\ycshift)$);
            \draw[\tickthick, \cunavcolor] (cl) -- (cr) node[midway, anchor=north, yshift=-\Lylabelshift]{$\color{\cunavcolor} \cunav$};
            \draw[\tickthick, \cunavcolor] ($(cl) + (0,\Lytickheight)$) -- ($(cl) + (0,-\Lytickheight)$);
            \draw[\tickthick, \cunavcolor] ($(cr) + (0,\Lytickheight)$) -- ($(cr) + (0,-\Lytickheight)$);
            \coordinate (ectrl) at ($(cr) + (\colorshift,0)$);
            \coordinate (ectrr) at ($(end) + (-\colorshift,0) + (0,-\ylshift-\ycshift)$);
            \draw[\tickthick, \ectrcolor] (ectrl) -- (ectrr) node[midway, anchor=north, yshift=-\Lylabelshift]{$\color{\ectrcolor} c_\infty$};
            \draw[\tickthick, \ectrcolor] ($(ectrl) + (0,\Lytickheight)$) -- ($(ectrl) + (0,-\Lytickheight)$);

            \coordinate (r1) at ($(tk) + (0,\yrshift)$);
            \coordinate (r2) at ($(r1) + (\xrshift, 0)$);
            \draw[-latex, \rthick] (r1) -- (r2) node[pos=0.9, anchor=south, yshift=4*\rlabelshift] {$r$};
            \draw[\rthick] ($(r1) + 
            (0,-\rtickheight)$) -- ($(r1) + (0,\rtickheight)$)  node[anchor=south, yshift=\rlabelshift]{$0$};
        \end{tikzpicture}
        \caption{Order at $t_k$, failure before delivery}
        \label{subfig:ord_at_tk_fail_bef_del}
    \end{subfigure}
    \hfill
    \begin{subfigure}{0.49\textwidth}
        \centering
        \begin{tikzpicture}
            \coordinate (0) at (0,0);
            \coordinate (end) at ($(0) + (\xdist, 0)$);
            \node[yshift=\tlabelshift cm] at (end) {$t$};
            \coordinate (tk) at ($(\xtk,0)$);
            \draw[\tickthick] ($(tk) + (0,\ytickheight)$) -- ($(tk) + 
            (0,-\ytickheight)$) node[anchor=north, yshift=-\ylabelshift]{\large $t_k$};
            \coordinate (tk1) at ($(tk) + (\tickdist,0)$);
            \draw[\tickthick] ($(tk1) + (0,\ytickheight)$) -- ($(tk1) + 
            (0,-\ytickheight)$) node[anchor=north, yshift=-\ylabelshift]{\large $t_{k+1}$};
    
            \foreach \k in {1,2,4}{
                \draw[\tickthick] ($(tk) + (\k*\tickdist,\ytickheight)$) -- ($(tk) + (\k*\tickdist,-\ytickheight)$);
            } 
            \coordinate (tkw) at ($(tk) + (3*\tickdist,0)$);
            \draw[\tickthick] ($(tkw) + (0,\ytickheight)$) -- 
            ($(tkw) + (0,-\ytickheight)$) node[anchor=north, yshift=-\ylabelshift]{\large $t_k+w$};
    
            \coordinate (Ll) at ($(tk) + (0,-\ylshift)$);
            \coordinate (Lr) at ($(Ll) + (\Ldist,0)$);
            \draw[\tickthick] (Ll) -- (Lr) node[midway, anchor=north, yshift=-\Lylabelshift]{$l$};
            \draw[\tickthick] ($(Ll) + (0,\Lytickheight)$) -- ($(Ll) + (0,-\Lytickheight)$);
            \draw[\tickthick] ($(Lr) + (0,\Lytickheight)$) -- ($(Lr) + (0,-\Lytickheight)$);
    
            \draw[\rulcolor, very thick, domain=\domain, samples=200] 
            plot ({\x}, {\plotfactor*lognormal(\x,\mu,\sigma)+\plotyshift});
            \node[] at ($(tkw) + (\textxshift, \textyshift)$) {$\color{\rulcolor} \frulr$};
            \pgfmathsetmacro\Tftwoy{\plotfactor*lognormal(\Tftwo,\mu,\sigma)+\plotyshift}
            \draw[\Tftwocolor, dashed, very thick] (\Tftwo,0) -- (\Tftwo,\Tftwoy) 
            node[very near start, yshift=-0.55cm]{\textcolor{\Tftwocolor}{$T_F^{(2)}$}};
    
            \draw[-latex, \linethick] (0) -- (end);
    
            \coordinate (cl) at ($(Lr) + (0,-\ycshift)$);
            \coordinate (cr) at ($(\Tftwo,0) + (0,-\ylshift-\ycshift)$);
            \draw[\tickthick, \cinvcolor] (cl) -- (cr) node[midway, anchor=north, yshift=-\Lylabelshift]{$\color{\cinvcolor} \cinv$};
            \draw[\tickthick, \cinvcolor] ($(cl) + (0,\Lytickheight)$) -- ($(cl) + (0,-\Lytickheight)$);
            \draw[\tickthick, \cinvcolor] ($(cr) + (0,\Lytickheight)$) -- ($(cr) + (0,-\Lytickheight)$);
            \coordinate (ectrl) at ($(cr) + (\colorshift,0)$);
            \coordinate (ectrr) at ($(end) + (-\colorshift,0) + (0,-\ylshift-\ycshift)$);
            \draw[\tickthick, \ectrcolor] (ectrl) -- (ectrr) node[midway, anchor=north, yshift=-\Lylabelshift]{$\color{\ectrcolor} c_\infty$};
            \draw[\tickthick, \ectrcolor] ($(ectrl) + (0,\Lytickheight)$) -- ($(ectrl) + (0,-\Lytickheight)$);

            \coordinate (r1) at ($(tk) + (0,\yrshift)$);
            \coordinate (r2) at ($(r1) + (\xrshift, 0)$);
            \draw[-latex, \rthick] (r1) -- (r2) node[pos=0.9, anchor=south, yshift=4*\rlabelshift] {$r$};
            \draw[\rthick] ($(r1) + 
            (0,-\rtickheight)$) -- ($(r1) + (0,\rtickheight)$)  node[anchor=south, yshift=\rlabelshift]{$0$};
        \end{tikzpicture}
        \caption{Order at $t_k$, failure after delivery}
    \label{subfig:ord_at_tk_fail_aft_del}
    \end{subfigure} 
    \\[1.5em]
        \begin{subfigure}{0.49\textwidth}
        \centering
        \begin{tikzpicture}
            \coordinate (0) at (0,0);
            \coordinate (end) at ($(0) + (\xdist, 0)$);
            \node[yshift=\tlabelshift cm] at (end) {$t$};
            \coordinate (tk) at ($(\xtk,0)$);
            \draw[\tickthick] ($(tk) + (0,\ytickheight)$) -- ($(tk) + 
            (0,-\ytickheight)$) node[anchor=north, yshift=-\ylabelshift]{\large $t_k$};
            \coordinate (tk1) at ($(tk) + (\tickdist,0)$);
            \draw[\tickthick] ($(tk1) + (0,\ytickheight)$) -- ($(tk1) + 
            (0,-\ytickheight)$) node[anchor=north, yshift=-\ylabelshift]{\large $t_{k+1}$};
    
            \foreach \k in {2,3,4}{
                \draw[\tickthick] ($(tk) + (\k*\tickdist,\ytickheight)$) -- ($(tk) + (\k*\tickdist,-\ytickheight)$);
            } 
            \coordinate (tkw) at ($(tk1) + (3*\tickdist,0)$);
            \draw[\tickthick] ($(tkw) + (0,\ytickheight)$) -- 
            ($(tkw) + (0,-\ytickheight)$) node[anchor=north, yshift=-\ylabelshift]{\large $t_{k+1}+w$};
    
            \coordinate (Ll) at ($(tk1) + (0,-\ylshift)$);
            \coordinate (Lr) at ($(Ll) + (\Ldist,0)$);
            \draw[\tickthick] (Ll) -- (Lr) node[midway, anchor=north, yshift=-\Lylabelshift]{$l$};
            \draw[\tickthick] ($(Ll) + (0,\Lytickheight)$) -- ($(Ll) + (0,-\Lytickheight)$);
            \draw[\tickthick] ($(Lr) + (0,\Lytickheight)$) -- ($(Lr) + (0,-\Lytickheight)$);
    
            \draw[\rulcolor, very thick, domain=\domain, samples=200] 
            plot ({\x}, {\plotfactor*lognormal(\x,\mu,\sigma)+\plotyshift});
            \node[] at ($(tkw) + (\textxshift, \textyshift)$) {$\color{\rulcolor} \frulr$};
            \pgfmathsetmacro\Tfoney{\plotfactor*lognormal(\Tfone,\mu,\sigma)+\plotyshift}
            \draw[\Tfonecolor, dashed, very thick] (\Tfone,0) -- (\Tfone,\Tfoney) 
            node[very near start, yshift=-0.55cm]{\textcolor{\Tfonecolor}{$T_F^{(1)}$}};
    
            \pgfmathsetmacro{\rrend}{\xtk + \tickdist + \Ldist}
            \draw[\linethick] (0) -- (\Tfone,0);
            \draw[dotted, \linethick] (\Tfone,0) -- (\rrend,0);
            \draw[-latex, \linethick] (\rrend,0) -- (end);
    
            \coordinate (cl) at ($(\Tfone,0) + (0,-\ylshift-\ycshift)$);
            \coordinate (cr) at ($(Lr) + (0,-\ycshift)$);
            \draw[\tickthick, \cunavcolor] (cl) -- (cr) node[midway, anchor=north, yshift=-\Lylabelshift]{$\color{\cunavcolor} \cunav$};
            \draw[\tickthick, \cunavcolor] ($(cl) + (0,\Lytickheight)$) -- ($(cl) + (0,-\Lytickheight)$);
            \draw[\tickthick, \cunavcolor] ($(cr) + (0,\Lytickheight)$) -- ($(cr) + (0,-\Lytickheight)$);
            \coordinate (ectrl) at ($(cr) + (\colorshift,0)$);
            \coordinate (ectrr) at ($(end) + (-\colorshift,0) + (0,-\ylshift-\ycshift)$);
            \draw[\tickthick, \ectrcolor] (ectrl) -- (ectrr) node[midway, anchor=north, yshift=-\Lylabelshift]{$\color{\ectrcolor} c_\infty$};
            \draw[\tickthick, \ectrcolor] ($(ectrl) + (0,\Lytickheight)$) -- ($(ectrl) + (0,-\Lytickheight)$);
            
            \coordinate (r1) at ($(tk) + (0,\yrshift)$);
            \coordinate (r2) at ($(r1) + (\xrshift, 0)$);
            \draw[-latex, \rthick] (r1) -- (r2) node[pos=0.9, anchor=south, yshift=4*\rlabelshift] {$r$};
            \draw[\rthick] ($(r1) + 
            (0,-\rtickheight)$) -- ($(r1) + (0,\rtickheight)$)  node[anchor=south, yshift=\rlabelshift]{$0$};
        \end{tikzpicture}
        \caption{Order at $t_{k+1}$, failure before delivery}
        \label{subfig:ord_at_tk1_fail_bef_del}
    \end{subfigure}
    \hfill
    \begin{subfigure}{0.49\textwidth}
        \centering
        \begin{tikzpicture}
            \coordinate (0) at (0,0);
            \coordinate (end) at ($(0) + (\xdist, 0)$);
            \node[yshift=\tlabelshift cm] at (end) {$t$};
            \coordinate (tk) at ($(\xtk,0)$);
            \draw[\tickthick] ($(tk) + (0,\ytickheight)$) -- ($(tk) + 
            (0,-\ytickheight)$) node[anchor=north, yshift=-\ylabelshift]{\large $t_k$};
            \coordinate (tk1) at ($(tk) + (\tickdist,0)$);
            \draw[\tickthick] ($(tk1) + (0,\ytickheight)$) -- ($(tk1) + 
            (0,-\ytickheight)$) node[anchor=north, yshift=-\ylabelshift]{\large $t_{k+1}$};
    
            \foreach \k in {2,3,4}{
                \draw[\tickthick] ($(tk) + (\k*\tickdist,\ytickheight)$) -- ($(tk) + (\k*\tickdist,-\ytickheight)$);
            } 
            \coordinate (tkw) at ($(tk1) + (3*\tickdist,0)$);
            \draw[\tickthick] ($(tkw) + (0,\ytickheight)$) -- 
            ($(tkw) + (0,-\ytickheight)$) node[anchor=north, yshift=-\ylabelshift]{\large $t_{k+1}+w$};
    
            \coordinate (Ll) at ($(tk1) + (0,-\ylshift)$);
            \coordinate (Lr) at ($(Ll) + (\Ldist,0)$);
            \draw[\tickthick] (Ll) -- (Lr) node[midway, anchor=north, yshift=-\Lylabelshift]{$l$};
            \draw[\tickthick] ($(Ll) + (0,\Lytickheight)$) -- ($(Ll) + (0,-\Lytickheight)$);
            \draw[\tickthick] ($(Lr) + (0,\Lytickheight)$) -- ($(Lr) + (0,-\Lytickheight)$);
    
            \draw[\rulcolor, very thick, domain=\domain, samples=200] 
            plot ({\x}, {\plotfactor*lognormal(\x,\mu,\sigma)+\plotyshift});
            \node[] at ($(tkw) + (\textxshift, \textyshift)$) {$\color{\rulcolor} \frulr$};
            \pgfmathsetmacro\Tftwoy{\plotfactor*lognormal(\Tftwo,\mu,\sigma)+\plotyshift}
            \draw[\Tftwocolor, dashed, very thick] (\Tftwo,0) -- (\Tftwo,\Tftwoy) 
            node[very near start, yshift=-0.55cm]{\textcolor{\Tftwocolor}{$T_F^{(2)}$}};
    
            \draw[-latex, \linethick] (0) -- (end);
            
            \coordinate (cl) at ($(Lr) + (0,-\ycshift)$);
            \coordinate (cr) at ($(\Tftwo,0) + (0,-\ylshift-\ycshift)$);
            \draw[\tickthick, \cinvcolor] (cl) -- (cr) node[midway, anchor=north, yshift=-\Lylabelshift]{$\color{\cinvcolor} \cinv$};
            \draw[\tickthick, \cinvcolor] ($(cl) + (0,\Lytickheight)$) -- ($(cl) + (0,-\Lytickheight)$);
            \draw[\tickthick, \cinvcolor] ($(cr) + (0,\Lytickheight)$) -- ($(cr) + (0,-\Lytickheight)$);
            \coordinate (ectrl) at ($(cr) + (\colorshift,0)$);
            \coordinate (ectrr) at ($(end) + (-\colorshift,0) + (0,-\ylshift-\ycshift)$);
            \draw[\tickthick, \ectrcolor] (ectrl) -- (ectrr) node[midway, anchor=north, yshift=-\Lylabelshift]{$\color{\ectrcolor} c_\infty$};
            \draw[\tickthick, \ectrcolor] ($(ectrl) + (0,\Lytickheight)$) -- ($(ectrl) + (0,-\Lytickheight)$);

            \coordinate (r1) at ($(tk) + (0,\yrshift)$);
            \coordinate (r2) at ($(r1) + (\xrshift, 0)$);
            \draw[-latex, \rthick] (r1) -- (r2) node[pos=0.9, anchor=south, yshift=4*\rlabelshift] {$r$};
            \draw[\rthick] ($(r1) + 
            (0,-\rtickheight)$) -- ($(r1) + (0,\rtickheight)$)  node[anchor=south, yshift=\rlabelshift]{$0$};
    
        \end{tikzpicture}
        \caption{Order at $t_{k+1}$, failure after delivery}
    \label{subfig:ord_at_tk1_fail_aft_del}
    \end{subfigure}
    \caption{Exemplary incurred costs for each of the considered branches of the decision tree in \Cref{fig:order_decision_tree}. The dotted black line indicates interruption of the component's function; hence, no costs are accumulated due to the continuation of the renewal-reward process.}
    \label{fig:order_graphs}
\end{figure}
}

\newcommand{\repgraphs}{
\begin{figure}[H]
\def\xdist{7.4}
\def\tlabelshift{-0.5}
\def\ytickheight{0.2}
\def\xtk{0.3}
\def\tickdist{4}
\def\ylabelshift{1}
\def\linethick{ultra thick}
\def\tickthick{very thick}
\def\rulcolor{red}
\def\domain{1.2:7.1}
\def\mu{1.4}
\def\sigma{0.4}
\def\Tfone{2.5}
\def\Tftwo{5.6}
\def\Tfonecolor{blue!80!black}
\def\Tftwocolor{green!50!black}
\def\plotyshift{0.1}
\def\plotfactor{2}
\def\textxshift{0}
\def\textyshift{1.3}
\def\ctickheight{0.1}
\def\yclabelshift{1}
\def\ycshift{1.2}
\def\cccolor{\Tfonecolor}
\def\cpcolor{\Tftwocolor}
\def\ectrcolor{orange!90!black}
\def\colorshift{0.04}
\def\yrshift{0.9}
\def\xrshift{1.0}
\def\rthick{thick}
\def\rtickheight{0.15}
\def\rlabelshift{0.8}

\centering
\begin{subfigure}{0.49\textwidth}
    \centering
    \begin{tikzpicture}
        \coordinate (0) at (0,0);
        \coordinate (end) at ($(0) + (\xdist, 0)$);
        \node[yshift=\tlabelshift cm] at (end) {$t$};
        \coordinate (tk) at ($(\xtk,0)$);
        \draw[\tickthick] ($(tk) + (0,\ytickheight)$) -- ($(tk) + 
        (0,-\ytickheight)$) node[anchor=north, yshift=-\ylabelshift]{\large $t_k$};
        \coordinate (tk1) at ($(tk) + (\tickdist,0)$);
        \draw[\tickthick] ($(tk1) + (0,\ytickheight)$) -- ($(tk1) + 
        (0,-\ytickheight)$) node[anchor=north, yshift=-\ylabelshift]{\large $t_{k+1}$};
    
        \draw[\rulcolor, very thick, domain=\domain, samples=200] 
        plot ({\x}, {\plotfactor*lognormal(\x,\mu,\sigma)+\plotyshift});
        \node[] at ($(tk1) + (\textxshift, \textyshift)$) {$\color{\rulcolor} \frulr$};
    
        \draw[-latex, \linethick] (0) -- (end);
    
        \coordinate (cp) at ($(tk) + (0,-\ycshift)$);
        
        \draw[dotted, \tickthick, \cpcolor] (tk) -- ($(tk) + (0,-0.4)$);
        \draw[dotted, \tickthick, \cpcolor] ($(cp)+(0,0.4)$) -- (cp) node[pos=1, anchor=north, yshift=-\yclabelshift]{$\color{\cpcolor} c_p$};
        \draw[\tickthick, \cpcolor] ($(cp) + (0,\ctickheight)$) -- ($(cp) + (0,-\ctickheight)$);
        \draw[\tickthick, \cpcolor] ($(cp) + (0,\ctickheight)$) -- ($(cp) + (0,-\ctickheight)$);
        \coordinate (ectrl) at ($(cp) + (\colorshift,0)$);
        \coordinate (ectrr) at ($(end) + (-\colorshift,0) + (0,-\ycshift)$);
        \draw[\tickthick, \ectrcolor] (ectrl) -- (ectrr) node[midway, anchor=north, yshift=-\yclabelshift]{$\color{\ectrcolor} c_\infty$};
        \draw[\tickthick, \ectrcolor] ($(ectrl) + (0,\ctickheight)$) -- ($(ectrl) + (0,-\ctickheight)$);

        \coordinate (r1) at ($(tk) + (0,\yrshift)$);
        \coordinate (r2) at ($(r1) + (\xrshift, 0)$);
        \draw[-latex, \rthick] (r1) -- (r2) node[pos=0.9, anchor=south, yshift=4*\rlabelshift] {$r$};
        \draw[\rthick] ($(r1) + 
        (0,-\rtickheight)$) -- ($(r1) + (0,\rtickheight)$)  node[anchor=south, yshift=\rlabelshift]{$0$};
    \end{tikzpicture}
\caption{Preventive replacement at $t_k$}
\label{subfig:rep_at_tk}
\end{subfigure}
\\[1.5em]
\begin{subfigure}{0.49\textwidth}
    \centering
    \begin{tikzpicture}
    \coordinate (0) at (0,0);
        \coordinate (end) at ($(0) + (\xdist, 0)$);
        \node[yshift=\tlabelshift cm] at (end) {$t$};
        \coordinate (tk) at ($(\xtk,0)$);
        \draw[\tickthick] ($(tk) + (0,\ytickheight)$) -- ($(tk) + 
        (0,-\ytickheight)$) node[anchor=north, yshift=-\ylabelshift]{\large $t_k$};
        \coordinate (tk1) at ($(tk) + (\tickdist,0)$);
        \draw[\tickthick] ($(tk1) + (0,\ytickheight)$) -- ($(tk1) + 
        (0,-\ytickheight)$) node[anchor=north, yshift=-\ylabelshift]{\large $t_{k+1}$};
    
        \draw[\rulcolor, very thick, domain=\domain, samples=200] 
        plot ({\x}, {\plotfactor*lognormal(\x,\mu,\sigma)+\plotyshift});
        \node[] at ($(tk1) + (\textxshift, \textyshift)$) {$\color{\rulcolor} \frulr$};

        \pgfmathsetmacro\Tfoney{%
            \plotfactor*lognormal(\Tfone,\mu,\sigma)+\plotyshift %
        }
        \draw[\Tfonecolor, dashed, very thick] (\Tfone,0) -- (\Tfone,\Tfoney) 
        node[very near start, yshift=-0.55cm]{\textcolor{\Tfonecolor}{$T_F^{(1)}$}};
    
        \draw[-latex, \linethick] (0) -- (end);
    
        \coordinate (cc) at ($(\Tfone,0) + (0,-\ycshift)$);
        \draw[dotted, \tickthick, \cccolor]  ($(\Tfone,0)$) -- ($(\Tfone,0)+(0,-0.2)$);
        \draw[dotted, \tickthick, \cccolor]  ($(cc)+(0,0.4)$) -- (cc) node[pos=1, anchor=north, yshift=-\yclabelshift]{$\color{\cccolor} c_c$};
        \draw[\tickthick, \cccolor] ($(cc) + (0,\ctickheight)$) -- ($(cc) + (0,-\ctickheight)$);
        \draw[\tickthick, \cccolor] ($(cc) + (0,\ctickheight)$) -- ($(cc) + (0,-\ctickheight)$);
        \coordinate (ectrl) at ($(cc) + (\colorshift,0)$);
        \coordinate (ectrr) at ($(end) + (-\colorshift,0) + (0,-\ycshift)$);
        \draw[\tickthick, \ectrcolor] (ectrl) -- (ectrr) node[midway, anchor=north, yshift=-\yclabelshift]{$\color{\ectrcolor} c_\infty$};
        \draw[\tickthick, \ectrcolor] ($(ectrl) + (0,\ctickheight)$) -- ($(ectrl) + (0,-\ctickheight)$);
        
        \coordinate (r1) at ($(tk) + (0,\yrshift)$);
        \coordinate (r2) at ($(r1) + (\xrshift, 0)$);
        \draw[-latex, \rthick] (r1) -- (r2) node[pos=0.9, anchor=south, yshift=4*\rlabelshift] {$r$};
        \draw[\rthick] ($(r1) + 
        (0,-\rtickheight)$) -- ($(r1) + (0,\rtickheight)$)  node[anchor=south, yshift=\rlabelshift]{$0$};
    \end{tikzpicture}
\caption{Do nothing at $t_k$, failure before $t_{k+1}$}
\label{subfig:rep_at_tk1_fail_bef_dt}
\end{subfigure}
\hfill
\begin{subfigure}{0.49\textwidth}
    \centering
    \begin{tikzpicture}
        \coordinate (0) at (0,0);
        \coordinate (end) at ($(0) + (\xdist, 0)$);
        \node[yshift=\tlabelshift cm] at (end) {$t$};
        \coordinate (tk) at ($(\xtk,0)$);
        \draw[\tickthick] ($(tk) + (0,\ytickheight)$) -- ($(tk) + 
        (0,-\ytickheight)$) node[anchor=north, yshift=-\ylabelshift]{\large $t_k$};
        \coordinate (tk1) at ($(tk) + (\tickdist,0)$);
        \draw[\tickthick] ($(tk1) + (0,\ytickheight)$) -- ($(tk1) + 
        (0,-\ytickheight)$) node[anchor=north, yshift=-\ylabelshift]{\large $t_{k+1}$};
    
        \draw[\rulcolor, very thick, domain=\domain, samples=200] 
        plot ({\x}, {\plotfactor*lognormal(\x,\mu,\sigma)+\plotyshift});
        \node[] at ($(tk1) + (\textxshift, \textyshift)$) {$\color{\rulcolor} \frulr$};

        \pgfmathsetmacro\Tftwoy{%
            \plotfactor*lognormal(\Tftwo,\mu,\sigma)+\plotyshift %
        }
        \draw[\Tftwocolor, dashed, very thick] (\Tftwo,0) -- (\Tftwo,\Tftwoy) 
        node[very near start, yshift=-0.55cm]{\textcolor{\Tftwocolor}{$T_F^{(2)}$}};
    
        \draw[-latex, \linethick] (0) -- (end);
    
        \coordinate (cp) at ($(tk1) + (0,-\ycshift)$);

        \draw[dotted, \tickthick, \cpcolor] (tk1) -- ($(tk1) + (0,-0.4)$);
        \draw[dotted, \tickthick, \cpcolor]  ($(\Tftwo,0)$) -- ($(\Tftwo,0)+(0,-0.2)$);
        \draw[dotted, \tickthick, \cpcolor]  ($(cp)+(0,0.4)$) -- (cp) node[pos=1, anchor=north, yshift=-\yclabelshift]{$\color{\cpcolor} c_p$};
        \draw[\tickthick, \cpcolor] ($(cp) + (0,\ctickheight)$) -- ($(cp) + (0,-\ctickheight)$);
        \draw[\tickthick, \cpcolor] ($(cp) + (0,\ctickheight)$) -- ($(cp) + (0,-\ctickheight)$);
        \coordinate (ectrl) at ($(cp) + (\colorshift,0)$);
        \coordinate (ectrr) at ($(end) + (-\colorshift,0) + (0,-\ycshift)$);
        \draw[\tickthick, \ectrcolor] (ectrl) -- (ectrr) node[midway, anchor=north, yshift=-\yclabelshift]{$\color{\ectrcolor} c_\infty$};
        \draw[\tickthick, \ectrcolor] ($(ectrl) + (0,\ctickheight)$) -- ($(ectrl) + (0,-\ctickheight)$);

        \coordinate (r1) at ($(tk) + (0,\yrshift)$);
        \coordinate (r2) at ($(r1) + (\xrshift, 0)$);
        \draw[-latex, \rthick] (r1) -- (r2) node[pos=0.9, anchor=south, yshift=4*\rlabelshift] {$r$};
        \draw[\rthick] ($(r1) + 
        (0,-\rtickheight)$) -- ($(r1) + (0,\rtickheight)$)  node[anchor=south, yshift=\rlabelshift]{$0$};
    \end{tikzpicture}
\caption{Do nothing at $t_k$, preventive replacement at $t_{k+1}$}
\label{subfig:rep_at_tk1}
\end{subfigure}
    
\caption{Exemplary incurred costs for each of the considered branches (a) to (c) of the decision tree in \Cref{fig:rep_decision_tree}.}
\label{fig:rep_graphs}
\end{figure}
}

\newcommand{\nnarchitecture}{
\begin{figure}[H]
\vspace{-1.5cm}
\def\sdist{1.9}
\def\sxdist{2.5}
\def\sydist{1}
\def\swidth{1.5}
\def\sheight{1}
\pgfmathsetmacro\sensorwidth{0.98*\swidth}
\pgfmathsetmacro\iheight{3.2*\sheight}
\pgfmathsetmacro\iwidth{1.1*\swidth}
\def\layerwidth{0.6}
\def\ldist{1.2}
\def\incolor{gray!70!white}
\def\fccolor{cyan!40!white}
\def\lstmcolor{lime!60}
\def\actcolor{orange!50!white}
\def\dropcolor{pink!90!black!}

\def\boxtextyshift{0.3}

\def\radius{0.3}
\def\termradius{0.2}
\def\boundline{very thick}
\def\arrowthick{very thick}
\def\dashthick{thick}
\tikzset{circle_node/.style={circle, draw, \boundline, minimum size=\radius cm, inner sep=0pt}}
\tikzset{term_node/.style={circle, draw, fill=black, \boundline, minimum size=\termradius cm, inner sep=0pt}}
\tikzset{square_node/.style={rectangle, draw, align=center, \boundline, inner sep=0pt}}
\tikzset{sens_node/.style={rectangle, align=center, inner sep=0pt}}

\tikzset{ghost_node/.style={circle, minimum size=0 cm, inner sep=0pt}}
    \centering
    \begin{tikzpicture}
    \coordinate (s0cord) at (0,0);
    \node[square_node, minimum height=\layerwidth cm, minimum width=\iheight cm, rotate=90, fill=\incolor] (s0) at (s0cord) {\Large Input Windows};
    \coordinate[right=\sdist cm of s0cord] (lstm1cord);
    \node[square_node, minimum height=\layerwidth cm, minimum width=\iheight cm, rotate=90, fill=\lstmcolor] (lstm1) at (lstm1cord) {\Large LSTM};
    \coordinate[right=\ldist cm of lstm1cord] (drop1cord);
    \node[square_node, minimum height=\layerwidth cm, minimum width=\iheight cm, rotate=90, fill=\dropcolor] (drop1) at (drop1cord) {\Large Dropout};
    \coordinate[right=\ldist cm of drop1cord] (lstm2cord);
    \node[square_node, minimum height=\layerwidth cm, minimum width=\iheight cm, rotate=90, fill=\lstmcolor] (lstm2) at (lstm2cord) {\Large LSTM};
    \coordinate[right=\sdist cm of lstm2cord] (gelucord);
    \node[square_node, minimum height=\layerwidth cm, minimum width=\iheight cm, rotate=90, fill=\actcolor] (gelu) at (gelucord) {\Large GELU};
    \coordinate[right=\ldist cm of gelucord] (drop2cord);
    \node[square_node, minimum height=\layerwidth cm, minimum width=\iheight cm, rotate=90, fill=\dropcolor] (drop2) at (drop2cord) {\Large Dropout};
    \coordinate[right=\ldist cm of drop2cord] (fccord);
    \node[square_node, minimum height=\layerwidth cm, minimum width=\iheight cm, rotate=90, fill=\fccolor] (fc) at (fccord) {\Large FC};
    \coordinate[right=\ldist cm of fccord] (outcord);
    \node[square_node, minimum height=\layerwidth cm, minimum width=\iheight cm, rotate=90, fill=\incolor] (out) at (outcord) {\Large RUL Classes};

    \draw[-latex, \arrowthick] (s0) -- (lstm1);
    \draw[-latex, \arrowthick] (lstm1) -- (drop1);
    \draw[-latex, \arrowthick] (drop1) -- (lstm2);
    \draw[-latex, \arrowthick] (lstm2) -- (gelu);
    \draw[-latex, \arrowthick] (gelu) -- (drop2);
    \draw[-latex, \arrowthick] (drop2) -- (fc);
    \draw[-latex, \arrowthick] (fc) -- (out);

    \def\dashxdist{0.55}
    \def\dashydist{1.8}
    \coordinate[left=\dashxdist cm of lstm1cord] (bllcord);
    \coordinate[above=\dashydist cm of drop1cord] (bltcord);
    \coordinate[right=\dashxdist cm of lstm2cord] (blrcord);
    \coordinate[below=\dashydist cm of lstm1cord] (blbcord);
    \draw[dashed, \dashthick] (bllcord) |- (bltcord) -| (blrcord) |- (blbcord) -| (bllcord);
    \coordinate[above=\boxtextyshift cm of bltcord] (bltextcord);
    \node[ghost_node] (bltext) at (bltextcord) {1-2 LSTM-Dropout layers};
    \end{tikzpicture}
    \captionsetup{width=\textwidth}
    \caption{NN RUL-classifier architecture with input and output (gray), lstm (lime), dropout (dark pink), and fully connected (cyan) layers, as well as GELU activation function (orange).}
    \label{fig:NN_architecture}
\end{figure}
}

\newcommand{\evaldiagram}{
\begin{figure}[H]
\def\xdist{3.2}
\def\ydist{2}
\def\y2dist{1.5}
\def\dswidth{1.85}
\def\dsheight{1}
\def\evalxdist{2.4}
\pgfmathsetmacro\yeval{\ydist + \y2dist/2}
\def\pad{0.4}
\def\evalcolor{gray!70!white}
\def\dscolor{cyan!40!white}
\def\modelcolor{lime!60}
\def\heurcolor{orange!50!white}
\def\boundline{very thick}
\def\arrowthick{very thick}
\def\dashthick{thick}
\tikzset{square_node/.style={rectangle, draw, align=center, \boundline, inner sep=0pt}}

    \centering
    \begin{tikzpicture}
    \coordinate (dscord) at (0,0);
    \node[square_node, minimum height=\dsheight cm, minimum width=\dswidth cm, fill=\dscolor] (ds) at (dscord) {Dataset};
    \coordinate[right=\xdist cm of dscord] (trdscord);
    \node[square_node, minimum height=\dsheight cm, minimum width=\dswidth cm, fill=\dscolor] (trds) at (trdscord) {Train \& val.\\set};
    \coordinate[above=\y2dist cm of trdscord] (tedscord);
    \node[square_node, minimum height=\dsheight cm, minimum width=\dswidth cm, fill=\dscolor] (teds) at (tedscord) {Test set};

    \coordinate[below=\ydist cm of dscord] (pmcord);
    \node[square_node, minimum height=\dsheight cm, minimum width=\dswidth cm, fill=\modelcolor] (pm) at (pmcord) {Prognostic\\model};
    \coordinate[right=\xdist cm of pmcord] (mtcord);
    \node[square_node, minimum height=\dsheight cm, minimum width=\dswidth cm, fill=\modelcolor] (mt) at (mtcord) {Hyperparam.\\search};
    \coordinate[right=\xdist cm of mtcord] (bmcord);
    \node[square_node, minimum height=\dsheight cm, minimum width=\dswidth cm, fill=\modelcolor] (bm) at (bmcord) {Best model};

    \coordinate[below=\ydist cm of pmcord] (dhcord);
    \node[square_node, minimum height=\dsheight cm, minimum width=\dswidth cm, fill=\heurcolor] (dh) at (dhcord) {Heuristic\\policy};
    \coordinate[below=\ydist cm of bmcord] (piocord);
    \node[square_node, minimum height=\dsheight cm, minimum width=\dswidth cm, fill=\heurcolor] (pio) at (piocord) {Parameter\\init./opt.};
    \coordinate[right=\xdist cm of piocord] (bhcord);
    \node[square_node, minimum height=\dsheight cm, minimum width=\dswidth cm, fill=\heurcolor] (bh) at (bhcord) {Configured\\policy};

    \coordinate[right=\evalxdist cm of bhcord, yshift=\yeval cm] (evalcord);
    \node[square_node, minimum height=\dsheight cm, minimum width=\dswidth cm, fill=\evalcolor] (eval) at (evalcord) {Policy\\evaluation};

    \draw[-latex, \arrowthick] (ds) -- (trds);
    \draw[-latex, \arrowthick] (trds) -- (mt);
    \draw[-latex, \arrowthick] (trds) -- (pio);
    \draw[-latex, \arrowthick] (ds) |- (teds);
    \draw[-latex, \arrowthick] (teds) -- (eval);
    \draw[-latex, \arrowthick] (pm) -- (mt);
    \draw[-latex, \arrowthick] (mt) -- (bm);
    \draw[-latex, \arrowthick] (bm) -- (pio);
    \draw[-latex, \arrowthick] (bm) -- (eval);
    \draw[-latex, \arrowthick] (dh) -- (pio);
    \draw[-latex, \arrowthick] (pio) -- (bh);
    \draw[-latex, \arrowthick] (bh) -- (eval);

    \coordinate (xalign) at ($(ds)!0.5!(trds)$);
    \node[above] at (xalign |- trds) {$80\%$};
    \node[above] at (xalign |- teds) {$20\%$};

    \begin{pgfonlayer}{background}
        \node[fit=(ds) (teds) (eval) (bh) (dh), fill=pink!50!white, inner sep=20pt, rounded corners] (pinkbox) {};

        \fill[teal!20!white, rounded corners]
        ($(dh.north west) + 
        (-\pad,\pad)$) --
        ($(bh.north west) + 
        (-\pad,\pad)$) --
        ($(eval.north west) + 
        (-\pad,\pad)$) --
        ($(eval.north east) + (\pad,\pad)$) --
        ($(eval.south east) + 
        (\pad,-1.5*\pad)$) -- 
        ($(bh.south east) + 
        (\pad,-1.5*\pad)$) -- 
        ($(dh.south west) + 
        (-\pad,-1.5*\pad)$) --
        ($(dh.north west) + 
        (-\pad,\pad)$);
    \end{pgfonlayer}
    \node[anchor=north west, yshift=0pt] at (pinkbox.north west)
    {{Repeat 100 times}};

    \node[anchor=south west, yshift=0pt] at ($(dh.south west) + (-\pad,-1.5*\pad)$)
    {{Repeat for each heuristic policy}};

    \end{tikzpicture}
    \captionsetup{width=\textwidth}
    \caption{Evaluation procedure for comparing different heuristic policies.}
    \label{fig:cmapss_evaluation_process}
\end{figure}
}

\begin{document}
\maketitle

\begin{abstract}
We propose a hybrid planning method for deriving prognostics-based predictive maintenance policies. The method accounts for the available decision options, the information on the future state of the system provided by a prognostic model, and the costs of the underlying renewal-reward process. It results in policies defined by only a few parameters, which can be determined based on theoretical considerations or by optimization from run-to-failure data. We demonstrate the potential of the method in two separate predictive maintenance decision settings: preventive replacement and preventive ordering. Our numerical investigations on a virtual remaining useful life simulator and the C-MAPSS dataset show that the derived policies rival the performance of optimized benchmark policies, while being significantly more robust against overfitting.
\end{abstract}

% keywords can be removed

%
%
%%%%%%%%%%%%%%%%%%%%%%%%%%%%%%%%%%%%%%%%%%%%%%%%%%%%%%%%%%%%%%%%%
%%%%%%%%%%%%%%%%%%%%%%%%%%%%%%%%%%%%%%%%%%%%%%%%%%%%%%%%%%%%%%%%%
%%%%%%%%%%%%%%%%%%%%%%%%%%%%%%%%%%%%%%%%%%%%%%%%%%%%%%%%%%%%%%%%%
%
%
\section{Introduction}
\label{sec:Introduction}
\subsection{Motivation and background}
\label{subsec:background}
Prognostics is concerned with predicting the future health of a system \cite{pecht2009prognostics}. Most often, one is interested in estimating the remaining useful life (RUL), i.e., the time left until a specific component or system fails and ceases to perform its intended function \cite{kim2017prognostics}. To quantify the uncertainty about the true underlying system state (not available due to, e.g., noisy sensors or hidden faults) as well as the uncertainties associated with making predictions (e.g., aleatory future loads), one typically resorts to a probability distribution over the RUL \cite{kim2017prognostics,si2011remaining}.
Prognostics is complemented by health management, which concerns making optimal maintenance decisions based on information obtained from the prognostic model \cite{kim2017prognostics,nguyen2019new}. 
This prognostic-based maintenance is summarized with the umbrella term \emph{Prognostics and Health Management} (PHM), or alternatively \emph{predictive maintenance} (PdM) \cite{ran2019survey}. 

This contribution focuses on health management, assuming a prognostic model is already available. There are two main options for obtaining maintenance decisions: black box models and heuristics. Black box decision models, such as Deep Reinforcement Learning (DRL) approaches  \cite{koutas2024investigation,lee2023deep,abbas2022interpretable,jimenez2024maintenance} often achieve state-of-the-art performance.
However, training these models either requires a large corpus of training data or a model that accurately represents the underlying deterioration process. In PHM applications, usually, only limited amounts of training data are available, and there is still a substantial gap between the developed models and the underlying reality \cite{chao2022fusing}.  
Another main drawback of DRL approaches lies in the limited transparency and interpretability of the resulting decision policies. In real-life applications, especially for safety-critical systems, this lack of insight is often critical. Therefore, the preferred approach in practice is the use of heuristic policies, i.e., simple rules on how to use the information at hand to choose maintenance actions. Heuristics are guided by intuition and, thus, are easily understood by engineers and operators. Examples for such heuristic policies include probability thresholds and fixed intervals \cite{kamariotis2024metric,luque2019risk,hodkiewicz2018optimizing,he2023condition}.

We propose a method to derive heuristic policies by relying on insights from renewal theory \cite{doob1948renewal, smith1958renewal}, which generalizes the Poisson process to (almost) arbitrary inter-arrival distributions. By regarding the inter-arrival times as lifetimes of individual components, and the associated failure costs as the rewards, the optimal maintenance problem can be seen as an application of \emph{renewal-reward processes}. Herein, the overarching goal is to find a policy that minimizes the long-running maintenance cost per unit time $c_\infty$.

%
%
%%%%%%%%%%%%%%%%%%%%%%%%%%%%%%%%%%%%%%%%
%
%
\subsection{Related work}
\label{subsec:related_work}
Many prognostics-based heuristic maintenance policies have been developed without the use of renewal theory \cite[e.g.,][]{lei2016phm,asgarpour2018bayesian,nguyen2019new,wang2024dynamic,wang2024predictive,kamariotis2024metric}. In these works, the parameters of the proposed heuristics are either already fully defined a-priori, or optimized based on available training data. Most commonly, a threshold approach is chosen. This threshold can be defined on various quantities, e.g., on the predicted RUL of a deterministic prognostic model, or on the predicted failure probability within a certain time period. 

Previous works \cite[e.g.,][]{zhuang2023prognostic,nguyen2019new,dong2025data,wang2025uncertainty,chen2022data} formulate the maintenance optimization problem as the minimization of the cost rate of the monitored components. At first glance, this seems similar to the $c_\infty$ measure of renewal-reward processes. But this cost rate only takes into account the expected costs and lifetime of the monitored component and \emph{not} of the system lifetime and overall component population. In addition, some of the works compute the target as an expectation over individual cost rates, which is \emph{not} equal to the ratio of expected cost over expected lifetime. There are also studies that use $c_\infty$ as a target for making maintenance decisions, such as \cite{fauriat2020optimization,mitici2023dynamic,yan2022online,prakash2017condition}. However, these works formulate the long-running maintenance cost per unit time based on the predicted RUL distribution of the single monitored component. By contrast, renewal theory dictates the use of the time-to-failure (TTF) distribution of the whole population, which is evidently different. 
The use of the component-specific instead of the population-based TTF distribution has been discussed by \citet{kamariotis2024metric}, who show the negative implications of this choice on the effectiveness of the resulting maintenance policies.
Examples for the application of renewal theory to heuristic maintenance strategies on the basis of the whole population of components can be found in the context of optimizing age-based replacement without monitoring, e.g., finding an a-priori optimal age-based replacement time of components \cite{jiang2009accurate}. This approach can be expanded in various ways, e.g., to settings with unknown starting component age \cite{sidibe2017preventive} or with the inclusion of shocks \cite{cha2017preventive}. Other age-based decision settings are also considered, e.g., for imperfect repairs \cite{finkelstein2015optimal}, combinations of repair and replacement (warranties) \cite{jack2009repair}, or ordering \cite{chelbi2001spare}.

Renewal theory has also been used in continuous monitoring settings in which the degradation state of the system can be perfectly observed. In this setting, the asymptotic failure rate of the system subject to a particular stochastic degradation model (e.g., a Gamma process) can be derived, which is then used to find the optimal policy. Examples include the investigation of delayed maintenance operations \cite{grall2006asymptotic}, combinations of imperfect as well as perfect repairs \cite{chien2012optimal,mercier2013modelling}, or the inclusion of multiple deterioration modes \cite{saassouh2007online}. Again, these studies do not rely on information obtained from a prognostic model. 

Similar to continuous monitoring, renewal theory can also be used to derive maintenance decisions based on degradation information obtained through inspections. The majority of these studies consider periodic inspections \cite{de2020review}. Examples include the consideration of imperfect repairs with unavailability constraints \cite{guo2013maintenance}, the optimization of inspection intervals for a given degradation threshold \cite{wu2015cost}, or the inclusion of imperfect inspections \cite{wang2015preventive}. The difference between these works and our contribution lies in their limitation to a specific stochastic degradation behaviour (e.g., Wiener or Gamma processes). Moreover, these approaches find an optimal a-priori deterioration threshold or inspection interval, instead of using the most recent RUL distribution for dynamic online decision optimization.

For preventive ordering (our second investigated decision setting), there are some works that perform prognostics-based cost rate minimization to find the optimal ordering time for spare parts, such as in \cite{elwany2008sensor}, or with an extension to stochastic lead times \cite{wang2015prognostics}. Similar to the discussed works on maintenance decisions, the optimization target of these works is the cost rate of the single monitored component and not the overall ordering cost rate when this policy is applied to the whole population of components. 

In addition, there are some works that consider renewal theory-based maintenance for multi-component system. Again, these works do not take into account the information available from a prognostics model. For a more detailed overview, the reader is referred to \cite{de2020review}. 

%
%
%%%%%%%%%%%%%%%%%%%%%%%%%%%%%%%%%%%%%%%%
%
%
\subsection{Contributions and overview}
\label{subsec:contributions_and_overview}

The objective of this paper is to exploit renewal theory to derive prognostics-based heuristic maintenance policies. Our specific contributions are summarized as follows.

\begin{itemize}
    \item We provide a general method for deriving predictive maintenance policies that incorporate the costs of the underlying renewal process. It consists of constructing a one-step decision tree at each decision time step, where the branch probabilities are computed from the latest RUL distribution, and the branch costs are obtained as the sum of the direct component costs and the costs associated with continuing the underlying renewal-reward process.
    \item We implement the proposed method in two decision settings, namely preventive replacement and preventive ordering, for which we provide step-by-step derivations of the resulting maintenance policies.
    \item We compare the derived maintenance policies against a set of benchmark policies using a virtual RUL simulator and the C-MAPSS dataset. Our derived policies rival the performance of optimized policies, without relying on optimization, and are hence significantly more robust against overfitting in low data regimes. 
    \item We release the RUL simulator as a common reference problem for evaluating and comparing maintenance policies in future studies. The corresponding source code for this work can be found at \url{github.com/Dakout/leaf_it_to_renewal}.
\end{itemize}

Besides its improved performance in the low-data regime, the proposed method offers additional advantages over the state of the art. The method merely requires the availability of a prognostic model that outputs RUL distributions. Hence, this work is not constrained to a specific stochastic degradation process. Furthermore, it is applicable across a range of decision settings; additional settings beyond those investigated in the paper are listed in the discussion section.

The remainder of this paper is structured as follows. \Cref{sec:Notation} introduces the general framework of predictive maintenance and the two decision settings investigated in this work. \Cref{sec:pdm_policies} proposes a new method for deriving predictive maintenance heuristics and applies it to the introduced decision settings, together with a set of benchmark heuristic policies. \Cref{sec:Numerical_investigations_RUL_simulator} conducts a detailed investigation of the proposed policies and benchmarks using a virtual RUL simulator. \Cref{sec:validation_on_cmapss_dataset} validates the obtained results on the prominently used C-MAPSS dataset. \Cref{sec:discussion} provides a thorough discussion on the limitations of the derived heuristics, as well as suggestions for further applications of the proposed method. \Cref{sec:conclusion} concludes this paper with a short summary.

\newpage
%
%
%%%%%%%%%%%%%%%%%%%%%%%%%%%%%%%%%%%%%%%%%%%%%%%%%%%%%%%%%%%%%%%%%
%%%%%%%%%%%%%%%%%%%%%%%%%%%%%%%%%%%%%%%%%%%%%%%%%%%%%%%%%%%%%%%%%
%%%%%%%%%%%%%%%%%%%%%%%%%%%%%%%%%%%%%%%%%%%%%%%%%%%%%%%%%%%%%%%%%
%
%
\section{Decision setting and notation}
\label{sec:Notation}
We adopt large parts of the setting from \cite{kamariotis2024metric}, i.e., we consider a single component subject to continuous monitoring, which is used to periodically make a probabilistic RUL prediction. This prediction is then fed into a maintenance policy, which outputs a decision to either do nothing or to perform a maintenance action\footnote{Here, a maintenance action refers to any action that helps to keep the system in a working condition, and, therefore, includes replacement as well as ordering.}. The workflow is illustrated in \Cref{fig:decision_setting}. 

\settinggraph

In this paper, we consider and discuss two PdM decision settings:

\begin{enumerate}
    \item[1.] Pure replacement: one decides when to replace a component subject to two possible types of replacement, namely preventive replacement with cost $c_p$ and corrective replacement upon failure with cost $c_c$, where $c_c > c_p$ to represent higher damages or longer downtime. It is assumed that a spare component is available.
        
    \item[2.] Pure ordering: one decides when to order a new component, which arrives after a deterministic lead time $l$, considering unavailability $\cunav$ and inventory holding $\cinv$ cost rates, where $\cunav > \cinv$. The component is replaced following failure as soon as the spare is available. In this paper, for simplicity, it is assumed that replacement of the failed component incurs no additional costs, but this assumption can easily be removed.
\end{enumerate}

$t$ is the time since a component has been installed. The set of discrete decision times for ordering/replacement is denoted with $\{t_k=k\cdot \Delta t,k=1,2,\ldots\}$. Further assumptions are that failures are self-announcing, replacements are perfect, instantaneous (take 0 time) and are performed immediately upon failure of the component. Furthermore, the lifetimes of all components are independent and identically distributed (with finite mean). Hence, we can use renewal-reward theory, where a renewal cycle is defined by the time interval between two successive replacements. Finally, we do not consider discounting of future costs, assuming that component lifetimes are comparably short.

Preventive replacement, order and failure times of a component are denoted with $T_{\mathrm{rep}}$, $T_{\mathrm{ord}}$,  and $\Tf$, respectively. If reference to a specific sample is needed, we explicitly refer to it via a sample superscript for component $i$, e.g., $\Tfi$; otherwise, the sample-based notation is avoided wherever possible for the sake of readability. We further introduce relative time $r_k=t-t_k$, which follows the definition of the remaining useful life $RUL(t)=\Tf - t$, where $RUL_k\coloneqq RUL(t_k)$. For the derivations and definitions of the proposed heuristic maintenance policies in this work, both notations ($r$ vs. $t$ \& $RUL$ vs. $\Tf$) will be used to reduce verbosity of the formulas as well as to facilitate intuition.

To evaluate the average performance of a certain maintenance policy, it needs to be repeatedly applied to numerous components that follow some failure time distribution. Upon failure/replacement, the component is assumed to be as good as new, therefore, we can use the elementary renewal-reward theory in both decision settings to compute the policy's long-running maintenance cost per unit time, $c_\infty$, as \cite{grimmett2020probability}
\begin{equation}
\label{eq:ectr}
    c_\infty = \lim_{t \rightarrow \infty} \frac{C(t)}{t} = \Ectratio.
\end{equation}
The expectations in the numerator and denominator of \Cref{eq:ectr} are both taken based on the resulting renewal cycles of the component, and $C$ and $T_{\mathrm{lc}}$ are determined by the specific policy of the respective decision settings. For the replacement setting, $C$ is either $c_p$ or $c_c$, and $T_{\mathrm{lc}}$ denotes either the point in time when the component is preventively replaced or fails, whichever happens first. For the ordering setting, $C$ is either the accumulated inventory or unavailability cost, and $T_{\mathrm{lc}}$ is always the component's failure time. Hence, $c_\infty$ refers to $\Ectrep$ or $\Ectorder$, depending on the decision context. The overarching goal is to find policies which minimize $c_\infty$.

The investigated policies are evaluated by computing a Monte Carlo (MC) estimate of the relative deviation between their cost rate and the optimal policy’s cost rate. This metric was introduced as $\hat{M}$ in \cite{kamariotis2024metric}:
\begin{equation}
    \label{eq:M}
    \hat{M} = \frac{\hat{C}_\infty - \hat{C}_{\infty}^{*}}{\hat{C}_{\infty}^{*}}, \qquad 0\leq\hat{M}.
\end{equation}
A higher value of $\hat{M}$ indicates a greater deviation from the optimal policy, and, consequently, poorer performance. $\hat{C}_{\infty}$ and $\hat{C}_{\infty}^{*}$ are MC estimates of the long-running cost rate for an investigated policy and the perfect policy, respectively. They are defined as \cite{kamariotis2024metric}
\begin{equation}
    \label{eq:R}
    \hat{C}_{\infty} = \frac{\frac{1}{n}\sum_{i=1}^n C^{(i)}}{\frac{1}{n}\sum_{i=1}^n T_{\mathrm{lc}}^{(i)}},
\end{equation}
\begin{equation}
    \label{eq:R_perfect}
    \hat{C}_{\infty}^{*} = \frac{\frac{1}{n}\sum_{i=1}^n C^{(i)^*}}{\frac{1}{n}\sum_{i=1}^n T_{\mathrm{lc}}^{(i)^*}}.
\end{equation}
$C^{*}$ and $T_{\mathrm{lc}}^{*}$ are the maintenance costs and component lifetimes resulting from a policy that has perfect foresight of each component's failure time and therefore can act optimally. For a given decision setting, they can be computed as the arguments that minimize a component's cost rate. For example, in the preventive replacement setting, the two options are 1) let the component fail or 2) replace at the last decision time before failure:
\begin{equation}
    \label{eq:cper_and_tperf_rep}
    \frac{C^{(i)^*}}{T_{\mathrm{lc}}^{(i)^*}} = 
    \min \Bigg( \frac{c_c}{\Tfi}, ~\frac{c_p}{\floor*{\frac{\Tfi}{\Delta t}} \cdot \Delta t}\Bigg).
\end{equation}
Likewise, in the preventive ordering setting, the lifetimes are equal to the failure times: $T_{\mathrm{lc}}^{(i)^*}\coloneqq \Tfi$, and the ordering times are chosen such that the inventory/unavailability costs are minimized:
\begin{equation}
    C^{(i)^{*}}=\min \left(
    \cinv \cdot \left[
    \Tfi - \Tinvstari
    \right],~
    \cunav \cdot\left[
    \Tinvstari+\Delta t +l - \Tfi
    \right]
    \right), \qquad \text{with}~ 
    \Tinvstari=\floor*{\frac{\Tfi-l}{\Delta t}} \cdot \Delta t,
\end{equation}
where $\floor* \cdot $ and $\ceil* \cdot $ represent the floor and ceiling functions.

In the considered setting, one obtains from a prognostic model at each $t_k$ a predicted RUL probability density function (PDF) denoted with $\frulr$. This prediction is based on all monitoring information available up to time $t_k$. This distribution is the input to the decision policy, as illustrated in \Cref{fig:decision_setting}. Throughout this work, $\frulr$ is used as the reference distribution to compute cumulative probabilities and expectations. In particular, the probability of the RUL being in the interval $(a,b]$, \Cref{eq:prob_RUL_betweeen_a_b}, and the expected RUL given that it lies in the interval $(a,b]$, \Cref{eq:exp_RUL_betweeen_a_b} are obtained as

\begin{equation}
    \label{eq:prob_RUL_betweeen_a_b}
    \Pr(a<RUL_{k}\leq b) \coloneqq \int_a^b \frulr~ \dr    
\end{equation}

\begin{equation}
    \EE{RUL_{k} \mid a<RUL_{k}\leq b} \coloneqq \frac{1}{\Pr(a<RUL_{k}\leq b)} 
    \cdot \int_a^b r \cdot \frulr ~\dr.
    \label{eq:exp_RUL_betweeen_a_b}
\end{equation}

The decision interval $\Delta t$ is highly problem-dependent. As an example, for aero-engines, \citet{kamariotis2024metric} suggest that $\Delta t=5-10$ flight cycles is realistic.

The available actions in the two decision settings, ``do nothing'', ``preventively replace'', and ``preventively order'', are abbreviated as DN, PR, and PO, respectively. 

%
%
%%%%%%%%%%%%%%%%%%%%%%%%%%%%%%%%%%%%%%%%%%%%%%%%%%%%%%%%%%%%%%%%%
%%%%%%%%%%%%%%%%%%%%%%%%%%%%%%%%%%%%%%%%%%%%%%%%%%%%%%%%%%%%%%%%%
%%%%%%%%%%%%%%%%%%%%%%%%%%%%%%%%%%%%%%%%%%%%%%%%%%%%%%%%%%%%%%%%%
%
%
\section{PdM policies}
\label{sec:pdm_policies}
In a general maintenance setting, individual decisions affect the state of the underlying system; they also impact future decision alternatives. Hence, one must consider the full sequence of decisions over the entire structure's life to quantify the effect of individual decisions. 
Because of uncertainty in the past, present, and future states and available information, this setting corresponds to a stochastic sequential decision problem (SDP) \cite{raiffa2000applied,kochenderfer2015decision,bismut2021optimal}. The full decision tree grows polynomially/(super-)exponentially in the set of the number of system states, observations, actions, components and considered time horizon, which is referred to as \emph{curse of dimensionality} \cite{bellman1966dynamic} and \emph{curse of history} \cite{pineau2006anytime}.

The decision on which action to take at time step $t_k$ is determined by a function called \emph{policy} $\pi_{t_k}$. For its definition, we adopt the notation from Markov decision processes (MDPs) \cite{sutton1998reinforcement}. Given a set of component states $s \in \mathcal{S}$ and available actions $a \in \mathcal{A}$, a policy assigns a probability distribution over $\mathcal{A}$, denoted as $\pi_{t_k}(a\mid s)=\Pr\left(a \mid S_{t_k}=s\right)$\footnote{This notation already assumes that the state $s$ is Markovian \cite{ethier2009markov}; otherwise, the policy depends on the full history of states and actions: $\pi_{t_k}(s_{0:t_k},a_{0:t_{k-1}})$ \cite{kochenderfer2015decision}.}. Thus, policies are generally stochastic. However, in this paper, we consider only deterministic policies, denoted as $\pi_{t_k}(s)$, which prescribe an action $a$  to any state $s$ \cite{kaelbling1996reinforcement}.

In practical cases, the true underlying state of the component is typically unknown, i.e., we have state uncertainty. Following the notation of \emph{partially observable} MDPs (POMDPs), this state uncertainty is represented by \emph{belief states} $b \in \mathcal{B}$, which are probability distributions over $\mathcal{S}$. If $b_{t_k}(s)=\Pr(S_{t_k}=s)$ denotes the probability of the system being in state $s$ at time $t$, then by the axioms of probability theory we obtain the following constraints: $0 \leq b_{t_k}(s) \leq 1, ~\forall s \in \mathcal{S}$ and $\sum_{s \in \mathcal{S}} b_{t_k}(s) = 1$. In POMDPs, the policy can be reformulated to be a function of the belief state: $\pi_{t_k}(b)$ \cite{kaelbling1998planning}. This belief is updated with Bayes' rule, taking into account the underlying system dynamics. 

For our predictive maintenance formulation, we define the true remaining useful life as the state of the system, and the RUL-PDF obtained by the prognostic model as the belief state:
\begin{align}
    \label{eq:state_as_rul}
    S_{t_k}&:=RUL_k \\
    b_{t_k}&:=f_{RUL_{k}}~.
    \label{eq_belief_as_rul_pdf}
\end{align}
This RUL-PDF is the input to all considered policies, which is why it is omitted in the following; we simply write $\pi_{t_k}=\pi_{t_k}(b)$. Note that the RUL-PDF is provided by a prognostic model, whose choice and training are deliberately left to the user. Thus, the resulting belief representation does not necessarily coincide with the Bayesian belief over the underlying degradation state, and the final maintenance performance is inherently conditioned on the predictive quality of the employed prognostic model. 

A \emph{strategy} $ \Pi=\{\pi_{t_1}, \pi_{t_2}, \dots\}$ denotes the set of policies for all time steps, where $\mathcal{P}$ is the space of all strategies. The optimal strategy $\Pi^*$ in the context of renewal-reward theory is then the one which minimizes the long-running maintenance cost per unit time:
\begin{equation}
    \label{eq:optimal_strategy}
    \Pi^*= \underset{\Pi \in \mathcal{P}}{\arg \min }~ \left[c_\infty \mid \Pi\right].
\end{equation}
Examining all possible $\Pi$ is computationally unfeasible; hence, one needs a way to reduce the strategy space $\mathcal{P}$. We do this by using a \emph{heuristic}, which is a parametrized strategy $\Pi_{\mathbf{w}}$ with parameters $\mathbf{w}=\{w_1, w_2, \dots, w_h\}$. Instead of \Cref{eq:optimal_strategy}, we thus solve the following optimization problem \cite{bismut2021optimal}:
\begin{equation}
    \label{eq:optimal_heuristic_strategy}
    \mathbf{w}^*= \underset{\mathbf{w}}{\arg \min }~ \left[c_\infty \mid \Pi_{\mathbf{w}}\right].
\end{equation}
Since $\Pi_{\mathbf{w}}$ forms a subspace of $\mathcal{P}$, solving \Cref{eq:optimal_heuristic_strategy} does not necessarily coincide with the global optimum obtained by solving \Cref{eq:optimal_strategy}. 

We consider an infinite time horizon, assuming that components are continuously replaced. In reality, the system will have a finite lifetime, but if the component life is significantly shorter than the system life, the infinite horizon is a good approximation \cite{hinderer1999approximate}. 
Because of the infinite horizon, the optimal strategies are \emph{stationary}, i.e., the policies at each time step are identical. We thus further simplify the policy notation by dropping the time index $\Pi = \pi = \pi_{t_1} = \pi_{t_2} = \dots$ \cite{bismut2021optimal}. For this reason, in this work, we use the terms policy and strategy interchangeably.

%
%
%%%%%%%%%%%%%%%%%%%%%%%%%%%%%%%%%%%%%%%%%%%%%%%%%%%%%%%%%%%%%%%%%%%%%%%%
%%%%%%%%%%%%%%%%%%%%%%%%%%%%%%%%%%%%%%%%%%%%%%%%%%%%%%%%%%%%%%%%%%%%%%%%
%
%
\subsection{Proposed method for deriving PdM policies}
\label{subsec:proposed_method}
Heuristic policies can be viewed as a way to simplify the underlying complexity of the full decision tree. In this paper, we propose a specific method to derive predictive maintenance policies by considering only the one-step decision tree. This corresponds to an online policy planning with a receding horizon of depth $1$. Such a one-step formulation can still be quite effective, as frequent replanning can compensate for the lack of longer-term modeling \cite{kochenderfer2022algorithms}. 

The motivation behind this originates in real-life settings, where maintenance actions (such as order, inspect, repair, replace) cannot be performed at arbitrary times, but are limited to a set of available options. As an example, for the inspection of an aircraft, the restrictions resulting from the planned operation time of the aircraft, as well as the availability of the crew, must be taken into account. In our setting, this restriction is represented by the decision interval $\Delta t$. The steps for our proposed method are enumerated as follows.
\begin{enumerate}[itemsep=1.8pt]
    \item[i.] Construct a simplified decision tree based on the available options and potential consequences.
    \item[ii.] Assign probabilities to each branch based on the RUL-PDF obtained by the prognostic model.
    \item[iii.] Assign to each branch the expected costs resulting from
    \begin{enumerate}[topsep=-5pt, itemsep=0pt]
        \item[a)] the monitored component,
        \item[b)] the underlying renewal-reward process.
    \end{enumerate}
    \item[iv.] Simplify the cost expressions and derive the policy that chooses the action with minimal expected cost according to the simplified decision tree.
    \item[v.] Optional: Perform data-based optimization of selected parameters of the derived policy.
\end{enumerate}
In step iii, renewal theory is used to quantify the expected benefits associated with keeping the component in service longer.
Specifically, we use the long-running cost per unit time of the respective maintenance policy ($c_\infty$). 
If the component is replaced by a new component at time $t_k$ instead of $t_k+\Delta t$, this will incur an expected cost that is equal to $c_\infty$ times $\Delta t$. 
This choice can be justified in two different ways. Firstly, the time-to-failure and maintenance costs of the future components are unknown; $c_\infty$ serves here as an expected value. Secondly, the assumption for the use of renewal theory is that the underlying renewal-reward process continues in the limit until infinity. Herein, $c_\infty$ represents the ``true'' cost rate of the maintenance policy.

The core of the proposed heuristic policies is the simplification of the full sequential decision problem into the binary one-step decision of \emph{perform an action now or perform an action later}. As a representation for \emph{act later}, we choose the action at the next possible decision time step $t_{k+1}$. Our reasoning behind this is that if it is optimal to perform the action at some later stage $t_{k+\tau}$, then we assume that it is better to act at $t_{k+1}$ as opposed to acting at $t_k$. Thus, ultimately, we consider the decision problem of acting immediately at $t_k$ versus at the next possible decision time $t_{k+1}$. The validity of this assumption is discussed in more detail in \Cref{sec:discussion}.

If the heuristic decision rule yields that it is better to act at $t_k$, then the action is performed at $t_k$. On the other hand, if the heuristic favours an action at $t_{k+1}$, then this does not necessarily mean that the action is actually performed at $t_{k+1}$; instead, the decision is revisited at $t_{k+1}$ following the online planning paradigm. At $t_{k+1}$, one obtains new information from the prognostic model and repeats the evaluation of the two decision options.
This simplified binary decision can be represented by a one-step decision tree, which will be used throughout this paper. 

We refer to the heuristic policies derived with the proposed method as \emph{discrete option assessment} policies and correspondingly use $doa$ to refer to these policies. In the following \Cref{subsec:heuristic_replacement_policies,subsec:Heuristic_ordering_policies}, we derive $doa$ policies together with different benchmarks for the preventive replacement and preventive ordering settings.

%
%
%%%%%%%%%%%%%%%%%%%%%%%%%%%%%%%%%%%%%%%%%%%%%%%%%%%%%%%%%%%%%%%%%%%%%%%%
%%%%%%%%%%%%%%%%%%%%%%%%%%%%%%%%%%%%%%%%%%%%%%%%%%%%%%%%%%%%%%%%%%%%%%%%
%
%
\subsection{Heuristic replacement policies}
\label{subsec:heuristic_replacement_policies}
In this section, we consider the pure replacement setting.
First, two benchmark heuristic policies are introduced in \Cref{subsubsec:rep_threshold,subsubsec:rep_optimization}, which are taken from \cite{kamariotis2024metric}. Afterwards, we present a new replacement heuristic based on the one-step decision tree in \Cref{subsubsec:rep_disc_opt_ass}.

%
%
%%%%%%%%%%%%%%%%%%%%%%%%%%%%%%%%%%%%%%%%%%%%%%%%%%%%%%%%%%%%%%%%%%%%%%%%
%%%%%%%%%%%%%%%%%%%%%%%%%%%%%%%%%%%%%%%%%%%%%%%%%%%%%%%%%%%%%%%%%%%%%%%%
%
%
\subsubsection{Reference heuristic: simple threshold}
\label{subsubsec:rep_threshold}
The first heuristic replacement policy, termed $rh1$, is a standard strategy used frequently in the literature for risk-based inspection planning \cite[e.g.,][]{luque2019risk,nielsen2018computational}) and PHM \cite[e.g.,][]{kamariotis2024metric,nguyen2019new,he2023condition}). $rh1$ performs a replacement action whenever the predicted probability of failure in the next time interval exceeds a threshold value $p_{\mathrm{thres}}$ (this threshold is the sole parameter $w$ in this heuristic):
\begin{equation}
    \label{eq:rep_threshold}
    a_{k} = \Pi^{rh1}_{p_{\mathrm{thres}}} = 
    \begin{cases}
        \mathrm{PR}, & \text { if } 
        p_{\mathrm{thres}} < \Pr\left(RUL_k \leq \Delta t\right) \\ 
        \mathrm{DN} & \text { else}.
    \end{cases}
\end{equation}
A default choice for the probability threshold is $p_{\mathrm{thres}}=c_p/c_c$\footnote{This threshold can be derived by comparing the preventive maintenance costs $c_p$ with the direct expected failure costs until the next decision time $\Pr\left(RUL_{k} \leq \Delta t\right)\cdot c_c$.}. However, this choice is suboptimal and can perform poorly \cite{kamariotis2024metric}. If training data is available, it is preferable to find $p_{\mathrm{thres}}^*$ by solving the optimization problem in \Cref{eq:optimal_heuristic_strategy}.

%
%
%%%%%%%%%%%%%%%%%%%%%%%%%%%%%%%%%%%%%%%%%%%%%%%%%%%%%%%%%%%%%%%%%%%%%%%%
%%%%%%%%%%%%%%%%%%%%%%%%%%%%%%%%%%%%%%%%%%%%%%%%%%%%%%%%%%%%%%%%%%%%%%%%
%
%
\subsubsection{Reference heuristic: optimal replacement time based on full RUL-PDF}
\label{subsubsec:rep_optimization}
As another benchmark, we consider \emph{PdM policy 3} by \citet{kamariotis2024metric}, which we denote by $rh2$. At each decision step, an optimization problem is solved to determine the optimal replacement time. The corresponding decision rule is
\begin{equation}
    \label{eq:rep_optimization}
    a_{k} = \Pi^{rh2}_{c_\infty} = 
    \begin{cases}
        \mathrm{PR}, & \text { if } r_{\mathrm{rep}, k}^* < \Delta t\\ 
        \mathrm{DN} & \text { else},
    \end{cases}
\end{equation}

where the optimal replacement time $r_{\mathrm{rep}, k}^*$ is defined as

\begin{equation}
    \label{eq:rep_optimization_obj_fun}
    r_{\mathrm{rep}, k}^* = 
    \underset{r_{\mathrm{rep}}}{\arg \min }~ 
    \Pfrrepk \cdot c_c + (1 - \Pfrrepk) \cdot c_p + c_\infty \cdot \int_{r_{\mathrm{rep}}}^{\infty}(r-\rrep) \cdot \frulr \dr,
\end{equation}
with $\Pfrrepk=\Pr\left(RUL_k \leq r_{\mathrm{rep}}\right)$. The idea behind $rh2$ is to not only consider the costs resulting from preventive/ corrective replacement of the monitored component, but also to include the costs incurred from replacing the component too early. The cost rate for this ``opportunity cost'' is set to $c_\infty$, following the logic explained in \Cref{subsec:proposed_method}. However, the true value of $c_\infty$ is unknown, since the final policy is the accumulation of individual maintenance decisions that are based on $c_\infty$ themselves. To escape this implicit definition, we use an estimate for this value, denoted as $c_{\infty,0}$. Three options to obtain such an estimate are presented in \cite{kamariotis2024metric}. In this work, we only use the first option, which considers a situation without monitoring (NM), where renewal theory can be used to find the optimal replacement time (component age). The objective function for a deterministic replacement time $\trep$ is formulated as:
\begin{align}
\label{eq:obj_func_rep_nm}
    g^{NM}(\trep) = \frac{\Pro{\Tf\leq\trep}\cdot c_c + (1-\Pro{\Tf\leq\trep})\cdot c_p}{\int_0^{\trep}t\cdot f_{\Tf}(t) \dt + (1-\Pro{\Tf\leq\trep})\cdot \trep}, \qquad \text{with~} \Pro{\Tf\leq\trep}=\int_0^{\trep}f_{\Tf}(t) \dt
\end{align}
and the corresponding initial estimate for the cost rate is found by solving the optimization problem:
\begin{equation}
    \label{eq:opt_prob_rep_nm}
    c_{\infty,0} = \min_{\trep} g^{NM}(\trep), \quad \trep \in (0, \infty).
\end{equation}
We have also performed detailed investigations on the performance of alternative versions of $rh2$, the results of which are given in \Cref{sec:app_further_investigations_rh2}.

%
%
%%%%%%%%%%%%%%%%%%%%%%%%%%%%%%%%%%%%%%%%%%%%%%%%%%%%%%%%%%%%%%%%%%%%%%%%
%%%%%%%%%%%%%%%%%%%%%%%%%%%%%%%%%%%%%%%%%%%%%%%%%%%%%%%%%%%%%%%%%%%%%%%%
%
%

\subsubsection{Proposed policy}
\label{subsubsec:rep_disc_opt_ass}

We derive a new $doa$ preventive replacement heuristic with the proposed method of \Cref{subsec:proposed_method}. Since we show all steps for a generic time step $t_k$, we omit the subscript $k$ in the derivation for ease of notation, and subsequently reintroduce it in the final policy expression.
In step i, we construct a one-step decision tree for the time horizon of $\Delta t$, depicted in \Cref{fig:rep_decision_tree}. 

\repairdecisiontrees

There are two decision alternatives (DN: do nothing; PR: preventive replacement) at the current time. Following PR, a new component will be in place during $\Delta t$. This new component is not explicitly modeled. If the do-nothing action is chosen, the component can fail within $\Delta t$ or be safe. In total this leads to four branches for this one-step decision tree, which are labelled (a) to (d) in \Cref{fig:rep_decision_tree}. The probabilities associated with the two branches following the DN action are readily obtained (step ii); the probability of failure during $\Delta t$ is $\Pf=\Pro{RUL\leq\Delta t}$. 

\Cref{fig:rep_graphs} illustrates the relation between the predicted RUL distribution, the failure time, and the resulting cost for branches (a) to (c). 
Following the proposed method, we do not assign costs to branch (d), since we restrict the maintenance options to $t_k$ and $t_{k+1}$, and do not consider the replacement options at $t_{k+\tau}, \tau>1$ when establishing the policy for time $t_k$. In the following, \Cref{fig:rep_graphs} forms the basis for establishing the expected costs associated with branches (a) to (c).\vspace{1em}

\repgraphs

Following \Cref{subsec:proposed_method}, we use $c_{\infty}$ for quantifying the costs of the continuation of the renewal-reward process upon failure or replacement of the component. If $c_{M@t_k}$ denotes the expected costs of performing maintenance action $M$ at time $t_k$, then the expected costs for preventively replacing the component immediately, i.e., branch/subfigure (a), is obtained as
\begin{equation}
\label{eq:h4_c_PR_tk}
    c_{PR@t_k}(r) = c_p + c_{\infty} \cdot r.
\end{equation}
Here, $r$ denotes the considered time horizon. 

The expected cost of planning to replace at the next decision time is a function of expected failure costs, i.e., branch (b),
\begin{equation}
\label{eq:h4_c_dn_fail}
    c_{DN@t_k,\Tf\leq t_{k+1}}(r) = c_c  + c_{\infty} \cdot \left(r- \EE{RUL \mid RUL < \Delta t} \right),
\end{equation}
as well as the expected preventive replacement costs, i.e., branch (c),
\begin{equation}
\label{eq:h4_c_dn_pr}
    c_{DN@t_k,PR@t_{k+1}}(r) =  c_p  + c_{\infty} \left(r - \Delta t \right).
\end{equation}
Since replacement at other points in time is not considered (branch d), the total expected cost of doing nothing at the current time step (planning to replace at $t_{k+1}$) is

\begin{equation}
    \label{eq:h4_c_PR_tk+1}
    c_{PR@t_{k+1}}(r) = \Pro{RUL \leq \Delta t} \cdot c_{DN@t_k,\Tf\leq t_{k+1}}(r) + \left(1-\Pro{RUL \leq \Delta t}\right) \cdot c_{DN@t_k,PR@t_{k+1}} (r).
\end{equation}
By comparing the expected costs for the two available decision options \Cref{eq:h4_c_PR_tk,eq:h4_c_PR_tk+1}, we derive the heuristic replacement rule by choosing the action with minimal expected cost. This is expressed as
\begin{equation}
\label{eq:a_rep_k_h4_gen}
    a_k =
    \Pi^{\mathrm{rep},doa}_{c_{\infty}, r}= 
        \begin{cases}
        \mathrm{PR}, & \text { if } c_{PR@t_{k}}(r)<c_{PR@t_{k+1}}(r)
        \\ \mathrm{DN} & \text { else}.
        \end{cases}
\end{equation}
Inserting \Cref{eq:h4_c_PR_tk,eq:h4_c_dn_fail,eq:h4_c_dn_pr,eq:h4_c_PR_tk+1} into the decision heuristic in \Cref{eq:a_rep_k_h4_gen}, results - after some algebra - in the final $doa$ replacement decision heuristic:
\begin{equation}
\label{eq:a_rep_k_h4}
    a_k = \Pi^{\mathrm{rep},doa}_{c_{\infty}} =
        \begin{cases}
        \mathrm{PR}, & \text { if } c_{\infty} \cdot \Delta t < \Pfk \left[c_c - c_p + c_{\infty}\left(\Delta t -  \EE{RUL_k \mid RUL_k \leq \Delta t}\right) \right]
        \\ \mathrm{DN} & \text { else}.
        \end{cases}
\end{equation}
The full derivation of \Cref{eq:a_rep_k_h4} is given in \Cref{sec:app_derivation_doa_final_costs}. Note that this policy does not depend on $r$, hence there is no need to choose this value.

\Cref{eq:a_rep_k_h4} offers a nice, intuitive interpretation: One should preventively replace the component only if the guaranteed costs incurred from the continuation of the renewal-reward process ($c_{\infty}\cdot \Delta t)$ are smaller than the expected costs of failure in the interval $[t_k,~t_{k+1}]$. These failure costs are composed of the additional replacement costs ($c_c-c_p$) as well as the costs associated with the continuation of the renewal-reward process, reduced by the conditional expected remaining useful life.

Following $rh2$, we again use the estimate $c_{\infty,0}\approx c_{\infty}$ presented in \Cref{subsubsec:rep_optimization} for the heuristic policy. Instead of this estimate, one can also regard $c_{\infty}$ as a mere parameter, and perform data-driven optimization by solving \Cref{eq:optimal_heuristic_strategy} with respect to the final policy performance $\hat{C}_{\infty}$, following step v of the general method. We investigate both options numerically in \Cref{sec:Numerical_investigations_RUL_simulator}.

%
%
%%%%%%%%%%%%%%%%%%%%%%%%%%%%%%%%%%%%%%%%%%%%%%%%%%%%%%%%%%%%%%%%%
%%%%%%%%%%%%%%%%%%%%%%%%%%%%%%%%%%%%%%%%%%%%%%%%%%%%%%%%%%%%%%%%%
%%%%%%%%%%%%%%%%%%%%%%%%%%%%%%%%%%%%%%%%%%%%%%%%%%%%%%%%%%%%%%%%%
%
%
\subsection{Heuristic ordering policies}
\label{subsec:Heuristic_ordering_policies}
In this section, we use the proposed method to develop a policy for preventive ordering. We assume here that if a spare component is in stock, the component can be replaced immediately upon failure. For holding a component in stock, one pays the inventory holding cost:
\begin{equation}
    \label{eq:inventory_cost}
        C_{\mathrm{stock}}
        = \max \left(\Tf - (T_{\mathrm{ord}} + l),~0 \right) \cdot \cinv.
\end{equation}
If a spare component is ordered too late, the component cannot be replaced, and one pays a delay cost of
\begin{equation}
    \label{eq:delay_cost}
        C_{\mathrm{delay}} 
        = \max \left(T_{\mathrm{ord}} + l - \Tf ,~0 \right) \cdot \cunav.
\end{equation}
If a monitored component fails before an order is placed, an ordering is enforced at the next possible order time $T_{\mathrm{ord}} = \ceil*{ \frac{\Tf}{\Delta t}}\cdot \Delta t$. Unavailability costs are accumulated while waiting for the component, and the renewal-reward process continues once the component has arrived after the lead time $l$. 

Next, we introduce a benchmark heuristic ordering rule and then derive the heuristic ordering policy according to our proposed method in \Cref{subsubsec:ord_disc_opt_ass}.

%
%
%%%%%%%%%%%%%%%%%%%%%%%%%%%%%%%%%%%%%%%%%%%%%%%%%%%%%%%%%%%%%%%%%
%%%%%%%%%%%%%%%%%%%%%%%%%%%%%%%%%%%%%%%%%%%%%%%%%%%%%%%%%%%%%%%%%
%
%
\subsubsection{Reference heuristic: simple threshold}
\label{subsubsec:simp_heur_for_comp_ord}
For benchmarking, we formulate a simple heuristic policy, which we denote with $oh1$. It requires ordering if the probability of failure after a time period $v=\ceil*{ \frac{l}{\Delta t}}\cdot \Delta t$ exceeds a threshold:

\begin{equation}
\label{eq:a_ord_k_h1}
    a_k = \Pi^{ oh1}_{p_{\mathrm{thres}}} = 
        \begin{cases}
        \mathrm{PO}, & \text { if } p_{\mathrm{thres}} <\Pro{RUL_{k}\leq v}
        \\ \mathrm{DN} & \text { else},
        \end{cases}
\end{equation}
in which the threshold $p_{\mathrm{thres}}$ is the heuristic parameter, which can be optimized or set to the default value $\cinv/\cunav$.

%
%
%%%%%%%%%%%%%%%%%%%%%%%%%%%%%%%%%%%%%%%%%%%%%%%%%%%%%%%%%%%%%%%%%
%%%%%%%%%%%%%%%%%%%%%%%%%%%%%%%%%%%%%%%%%%%%%%%%%%%%%%%%%%%%%%%%%
%%%%%%%%%%%%%%%%%%%%%%%%%%%%%%%%%%%%%%%%%%%%%%%%%%%%%%%%%%%%%%%%%
%
%
\subsubsection{Proposed policy}
\label{subsubsec:ord_disc_opt_ass}
Following our proposed method, we establish a one-step decision tree in which we consider two options: either order a spare component now at $t_k$ or at the next possible decision time $t_{k+1}$. The decision tree with associated branch probabilities and expected costs (steps i. to iii.) is shown \Cref{fig:order_decision_tree}.  
The branches are again labelled (a) to (d). Example realizations of these four branches are shown in \Cref{fig:order_graphs}, with coinciding subfigure labels. Again, we omit the subscript $k$ in the derivation of the policy for ease of notation, and then reintroduce it in the final policy expression. \vspace{1em}

\orderdecisiontrees

The long-running maintenance cost per unit time, $c_\infty$, is again used to account for the costs associated with the continuation of the renewal-reward process after a component has been replaced. 
Based on this, the expected costs associated with the branches (a) - (d), corresponding to \Cref{subfig:ord_at_tk_fail_bef_del,subfig:ord_at_tk_fail_aft_del,subfig:ord_at_tk1_fail_bef_del,subfig:ord_at_tk1_fail_aft_del} are given in \Cref{eq:h_ord_tk_fail,eq:h_ord_tk_surv,eq:h_ord_tk+1_fail,eq:h_ord_tk+1_surv}:
\begin{equation}
    \label{eq:h_ord_tk_fail}
        c_{O@t_k, \Tf\leq t_k+l}(r) 
        = \left( 
        l - \EE{RUL \mid RUL \leq l}
        \right) \cdot \cunav 
        + c_{\infty} \cdot (r - l)
\end{equation}
\begin{equation}
    \label{eq:h_ord_tk_surv}
    c_{O@t_k, \Tf>t_k+l}(r) = 
    \left(\EE{RUL \mid RUL > l} - l\right) \cdot \cinv +
    c_{\infty} \cdot (r - \EE{RUL \mid RUL > l})
\end{equation}
\begin{equation}
    \label{eq:h_ord_tk+1_fail}
    c_{O@t_{k+1}, \Tf\leq t_{k+1}+l}(r) 
    = \left( 
    l + \Delta t - \EE{RUL \mid RUL \leq l + \Delta t}
    \right) \cdot \cunav
    + c_{\infty} \cdot (r - l - \Delta t)
\end{equation}
\begin{equation}
\label{eq:h_ord_tk+1_surv}
    c_{O@t_{k+1}, \Tf>t_{k+1}+l}(r) 
    = \left( 
    \EE{RUL \mid RUL > l + \Delta t} - l - \Delta t
    \right) \cdot \cinv + 
    c_{\infty} \cdot (r - \EE{RUL \mid RUL > l + \Delta t}).
\end{equation}
\ordergraphs
The resulting total expected costs of the decisions ``ordering at $t_k$'' and ``ordering at $t_{k+1}$'' are
\begin{align}
\label{eq:h_ord_tk}
    c_{O@t_k}(r) &= 
    \Pro{RUL \leq l} \cdot c_{O@t_k, \Tf\leq t_k+l}(r) + 
    \Pro{RUL > l} \cdot c_{O@t_k, \Tf>t_k+l}(r) 
    \\
    c_{O@t_{k+1}}(r) & = 
    \Pro{RUL \leq l + \Delta t} \cdot c_{O@t_{k+1}, \Tf\leq t_{k+1}+l}(r) + 
    \Pro{RUL > l + \Delta t} \cdot c_{O@t_{k+1}, \Tf>t_{k+1}+l}(r).
\end{align}

The heuristic ordering policy is to choose the action with minimal expected cost (step iv):
\begin{equation}
\label{eq:a_ord_k_h4_gen}
    a_k = \Pi^{\mathrm{ord}, doa}_{c_{\infty},r} = 
        \begin{cases}
        \mathrm{PO}, & \text { if } c_{O@t_{k}}(r)<c_{O@t_{k+1}}(r)
        \\ \mathrm{DN} & \text { else}.
        \end{cases}
\end{equation}
Investigating the individual terms in detail, this expression can be ``simplified'' to:
\begin{equation}
\label{eq:a_ord_k_h4}
    a_k = \Pi^{\mathrm{ord}, doa}_{c_{\infty}} =
        \begin{cases}
        \mathrm{PO}, & \text { if } 
        \begin{alignedat}[t]{2}
            0 < ~& \Pro{RUL_k \leq l}\cdot \Delta t \cdot \left(\cunav + \cinv - c_{\infty} \right) - \cinv \cdot \Delta t + 
            \\
            &  \Pro{l<RUL_k\leq l+\Delta t} \cdot \left( 
            \cunav + \cinv  - c_{\infty} \right) \cdot \\
            & \left(l + \Delta t - \EE{RUL_k \mid l< RUL_k \leq l + \Delta t}\right)
        \end{alignedat}
        \\ \mathrm{DN} & \text { else}.
        \end{cases}
\end{equation}
The full derivation of \Cref{eq:a_ord_k_h4} is given in \Cref{sec:app_derivation_oh4_final_costs}.

As in \cite{kamariotis2024metric}, we consider 3 approaches for choosing an appropriate initial guess $c_{\infty}\approx c_{\infty,0}$:

\begin{enumerate}
    \item[]\emph{1. Without monitoring}: 

    Considering a situation without monitoring (NM), renewal theory can be used to find the optimal time for ordering. The objective function for a deterministic ordering time $\tord$ is formulated as:

    \begin{align}
    \label{eq:obj_func_ord_nm}
        g^{NM}(\tord) = \frac{1}{\mutf} [ & 
        \Pro{\Tf \leq \tord + l} \cdot 
        \left( \tord + l - \EE{\Tf \mid \Tf \leq \tord + l} \right) \cdot \cunav + 
        \notag \\
        & \Pro{\Tf > \tord + l} \cdot 
        \left(\EE{\Tf \mid \Tf > \tord + l} - \tord - l \right) \cdot \cinv ],
    \end{align}
    and the corresponding initial estimate for the cost rate is found by solving the optimization problem:
    \begin{equation}
        \label{eq:opt_prob_ord_nm}
        c_{\infty,0} = \min_{\tord} g^{NM}(\tord), \quad \tord \in (0, \infty).
    \end{equation}
    With monitoring and prognostics, the policy's cost rate will be lower than without. Hence, $c_{\infty,0}$ according to \Cref{eq:opt_prob_ord_nm} is an upper bound on the true long-running maintenance cost per unit time. As a consequence, if this value of $c_{\infty,0}$ is used in the policy of \Cref{eq:obj_func_ord_nm}, the penalty for early ordering is overestimated, which results in a PdM policy that prefers late ordering. 

    \item[]\emph{2. Perfect ordering}:
    
    An alternative approach to determining an initial estimate for the cost rate is to assume a perfect ordering policy, which could be obtained when the time of failure is known exactly, i.e., with perfect prognostics. The perfect ordering policy would be to order at $\Tf-l$ if ordering were always possible. However, in our setting, one can only order at discrete decision times, so it must be determined if it is optimal to order at the nearest decision time before or after $\Tf-l$. That depends on when exactly the component fails in this interval. Paying the inventory holding cost for almost the entire period $\Delta t$, can be worse than to pay the unavailability cost for only a fraction of $\Delta t$. 
    
    We derive $c_{\infty,0}$ with such a perfect policy by making a suitable set of assumptions detailed in \Cref{sec:app_derivation_oh4_ectr_option2}. The resulting analytical estimate for the long-running perfect ordering costs per unit time is
    \begin{equation}
    \label{eq:ectorder_opt_2}
        c_{\infty,0} = \frac{\cinv \cdot \cunav}{2(\cinv + \cunav)} \frac{\Delta t}{\mutf}.
    \end{equation}

    The $c_{\infty,0}$ value obtained via \Cref{eq:ectorder_opt_2} represents a lower bound to the cost rate achieved with the actual ordering policy; resulting in a more conservative ordering policy. 

    \item[]\emph{3. Intermediate value:}

    Any $c_{\infty,0}$ value between the upper bound (option 1) and lower bound (option 2) can be chosen. 
\end{enumerate}
The presented initialization options 1. to 3. can either be computed analytically with an analytical time-to-failure distribution or in a sample-based manner based on a set of run-to-failure trajectories. 

Following step v, one can also regard $c_{\infty}$ as a mere optimization parameter, and determine its value via a data-based optimization with respect to the final policy performance $\hat{C}_{\infty}$, i.e., by solving \Cref{eq:optimal_heuristic_strategy}.

%
%
%%%%%%%%%%%%%%%%%%%%%%%%%%%%%%%%%%%%%%%%%%%%%%%%%%%%%%%%%%%%%%%%%
%%%%%%%%%%%%%%%%%%%%%%%%%%%%%%%%%%%%%%%%%%%%%%%%%%%%%%%%%%%%%%%%%
%%%%%%%%%%%%%%%%%%%%%%%%%%%%%%%%%%%%%%%%%%%%%%%%%%%%%%%%%%%%%%%%%
%
%
\section{Numerical investigations on a virtual RUL simulator}
\label{sec:Numerical_investigations_RUL_simulator}
We investigate the performance of the proposed heuristic policies in a controlled setting that enables accurate testing and verification. To this end, we make use of the RUL simulator introduced in \cite{kamariotis2024metric}. The simulator provides \emph{calibrated}\footnote{A well-calibrated probabilistic model outputs predictions whose credible intervals correspond to the observed frequencies in the long run \cite{gneiting2007probabilistic}, i.e., the true RUL lies in the $(1-\alpha)\%$ credible interval in $(1-\alpha)\%$ of cases.} RUL predictions. This simulator enables the generation of an unlimited number of test data, which facilitates the exact evaluation of the different policies' performances.

We assume a Gaussian time-to-failure distribution of components: $\Tf \sim \mathcal{N}(\mu=225, \sigma=40)$. The RUL simulator models prediction errors $\epsilon_k$ at time step $t_k$ by a lognormal process $RUL_{k}=\epsilon_k\cdot RUL_{\mathrm{true},k}$, where the logarithms of the prediction errors $\ln (\epsilon_k)$ follow a zero-mean multivariate normal distribution. This leads to marginal RUL prediction distributions equal to
\begin{equation}
    \label{eq:lnrul_dist}
    \ln\left(RUL_{k}\right) \sim \mathcal{N}\left(\mu=\ln(RUL_{\mathrm{true},k}) + \ln (\epsilon_k), ~\sigma=0.4\right),
\end{equation}
i.e., given an error sample $\epsilon_k^{(i)}$, $f_{RUL_{k}}$ is the PDF of a lognormal distribution with mean $\mu_{RUL_{k}}\approx1.083\epsilon_k^{(i)}\cdot RUL_{\mathrm{true},k}$ and standard deviation $\sigma_{RUL_{k}}\approx0.451\epsilon_k^{(i)}\cdot RUL_{\mathrm{true},k}$. Thus, the simulator yields slightly biased RUL predictions, which is, however, not critical for performance. The RUL predictions at different $t_k$ are correlated, with correlation length equal to $l_{\mathrm{corr}}=50$. 
A detailed description of the RUL simulator can be found in \Cref{sec:virtual_RUL_simulator}. An example realization of RUL predictions is shown in \Cref{fig:decision_setting}.

%
%
%%%%%%%%%%%%%%%%%%%%%%%%%%%%%%%%%%%%%%%%%%%%%%%%%%%%%%%%%%%%%%%%%
%%%%%%%%%%%%%%%%%%%%%%%%%%%%%%%%%%%%%%%%%%%%%%%%%%%%%%%%%%%%%%%%%
%
%

\subsection{Replacement policy performance}
\label{subsec:replacement_performance}
\begin{figure}[H]
    \centering
    \includegraphics[width=0.93\textwidth]{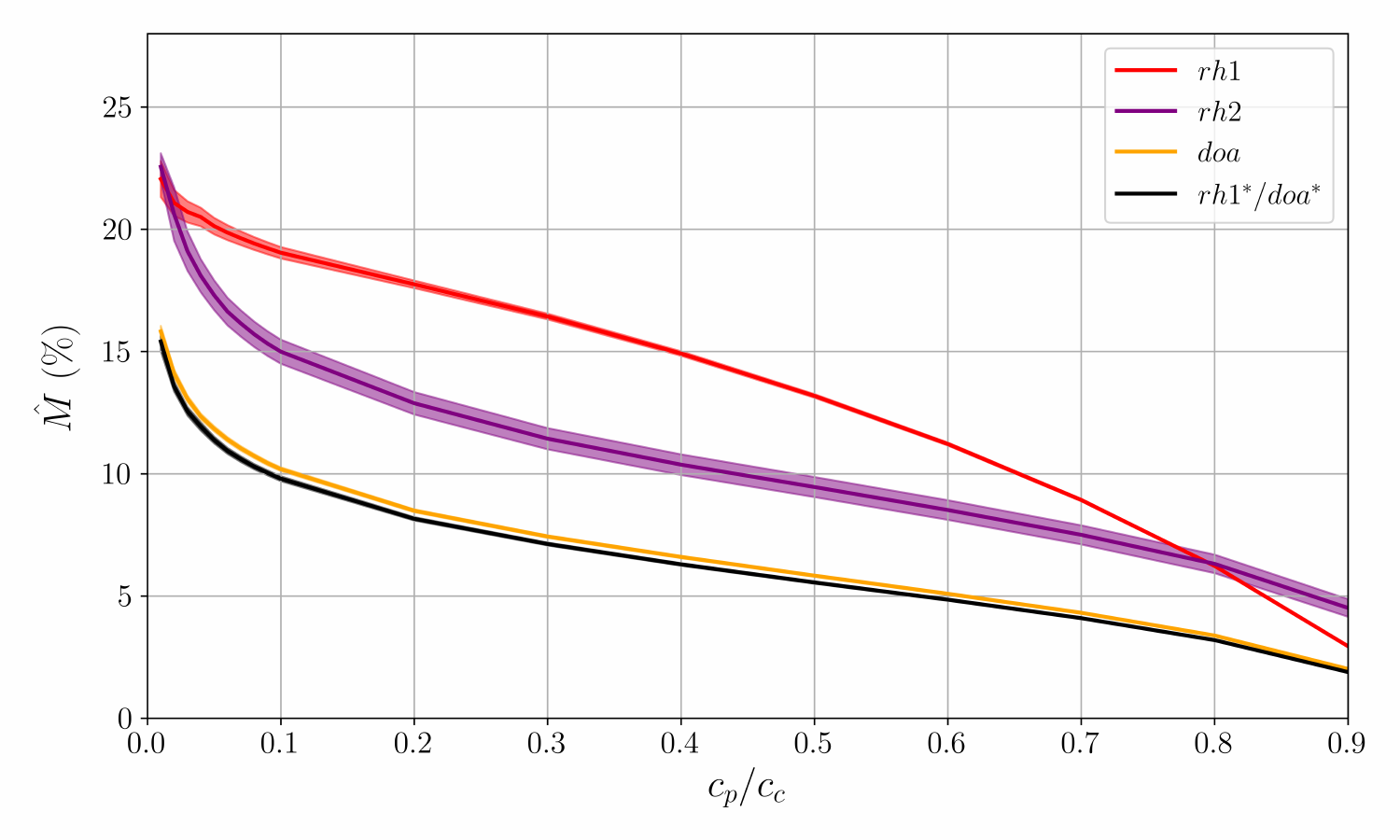}
    \caption{Comparison of our proposed $doa$ replacement policy with the benchmarks $rh1$ and $rh2$. Since the optimized versions $rh1^*$ and $doa^*$ yield the same performance, their performance is plotted as a single black curve. The following number of MC samples were used for the evaluation of each policy: $rh1: 10^6$, $rh2: 1.2\cdot10^4$, $doa: 10^6$. For the optimized policies, $10^6$ MC samples were used for training. The figure shows the mean of each policy as a solid line, as well as the corresponding 95\% credible intervals.}
    \label{fig:rh_perf_overview}
\end{figure}

We compare the performance of our proposed $doa$ replacement policy with the benchmark heuristics $rh1$ and $rh2$ introduced in \Cref{subsec:heuristic_replacement_policies}. For $rh1$ and $doa$, we include the default versions without parameter optimization as well as their optimized counterparts, denoted as $rh1^*$ and $doa^*$, respectively. The comparison is performed for $\Delta t=10$ over a range of $c_p/c_c$ ratios. The results are summarized in \Cref{fig:rh_perf_overview}. The performances of $rh1^*$ and $doa^*$ are indistinguishable when graphed together with the other policies, which is why they are shown as one color.

One can see in \Cref{fig:rh_perf_overview} that the default (unoptimized) version of our proposed $doa$ heuristic achieves much lower cost ratios compared to the default versions of the benchmark policies. The reason why $rh1$ performs worse compared to $doa$ is that it only considers the expected preventive/corrective costs directly attributable to the monitored component, and does not include the costs of the underlying renewal-reward process. $rh2$ improves upon this fallacy by including $c_{\infty}$ in the objective function and solving for the optimal replacement time $\rrep^{*,rh2}$ at each decision time step. The final maintenance decision is taken based on whether $\rrep^{*,rh2}$ lies within the next decision time interval or not. However, this optimal replacement time is not available to the operator, but only the discrete decision times $t_k, ~t_{k+1}, ...$. Hence, having $\rrep^{*,rh2}$ does not necessarily help with taking the optimal maintenance decision regarding the discrete decision times, which ultimately leads to worse performance. 
By contrast, our proposed $doa$ replacement heuristic specifically assesses the associated costs of the discrete replacement options at hand (or a one-step proxy thereof), and chooses the maintenance action accordingly. 
Importantly, the performance of the default $doa$ is only marginally worse compared to the optimized policy. This is an indication that the proposed policy works even with a crude assumption on $c_{\infty}$ and might be used even with no training data to optimize the policy (training data is of course still needed to train the RUL prediction). 

If one optimizes over the free parameter of each replacement heuristic, i.e., $p_{\mathrm{thres}}$ for $rh1$, or $c_{\infty}$ for $doa$, the performances of the optimized policies converge to the same best performance. 
This indicates that it is possible to compensate an initially ``bad'' heuristic with optimization.
However, these results are based on optimizing with a very large number of training data, specifically $10^6$ run-to-failure samples. 
In practice, such large amounts of data are not available, especially for safety-critical systems. Thus, we investigate the effect of training data size on the performances of the respective heuristics. 

\begin{figure}[H]
    \centering
    \begin{subfigure}{0.49\textwidth}
        \includegraphics[width=\textwidth]{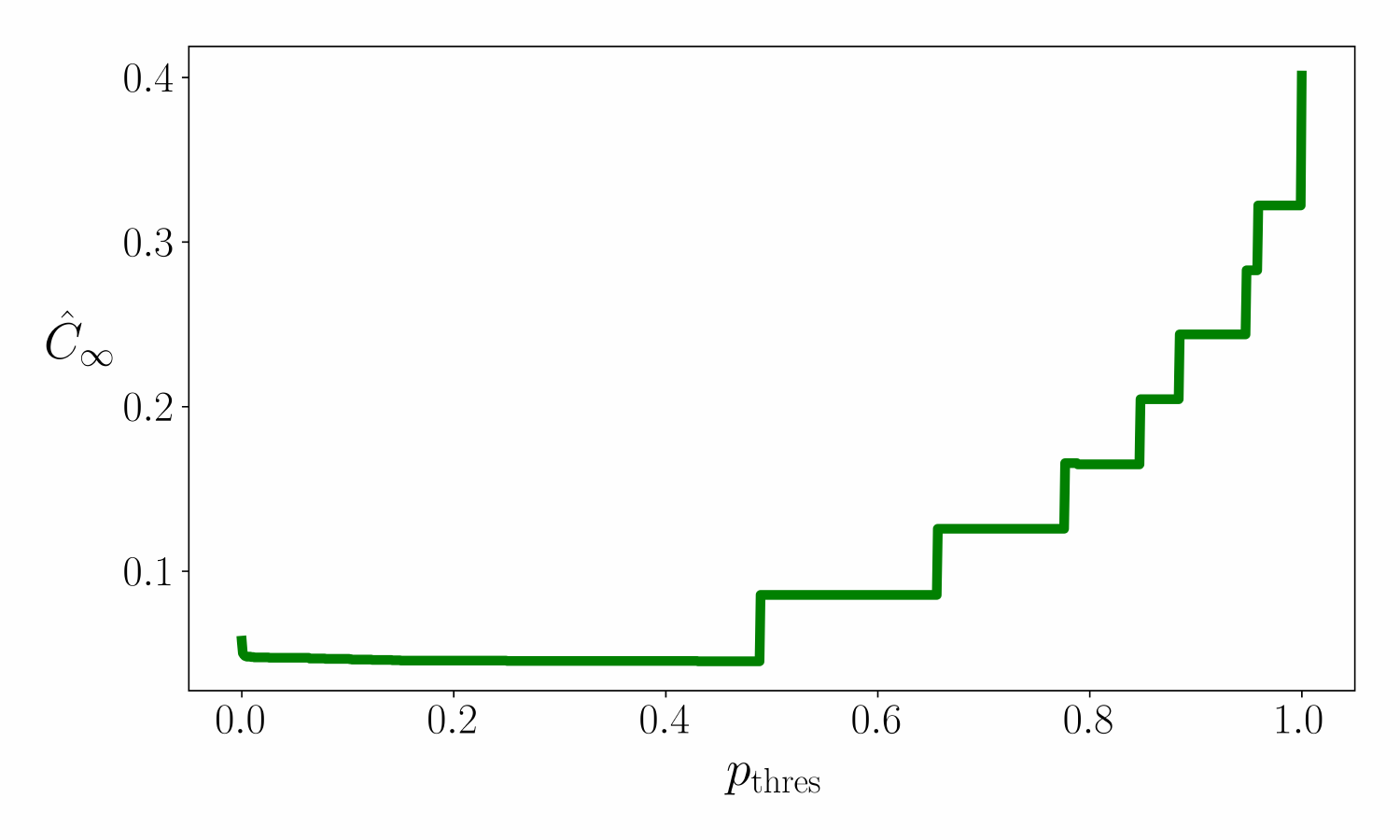}
        \caption{$rh1^*$}
        \label{subfig:rh1_perf_vs_pthres}
    \end{subfigure}
    \\
    \begin{subfigure}{0.49\textwidth}
        \includegraphics[width=\textwidth]{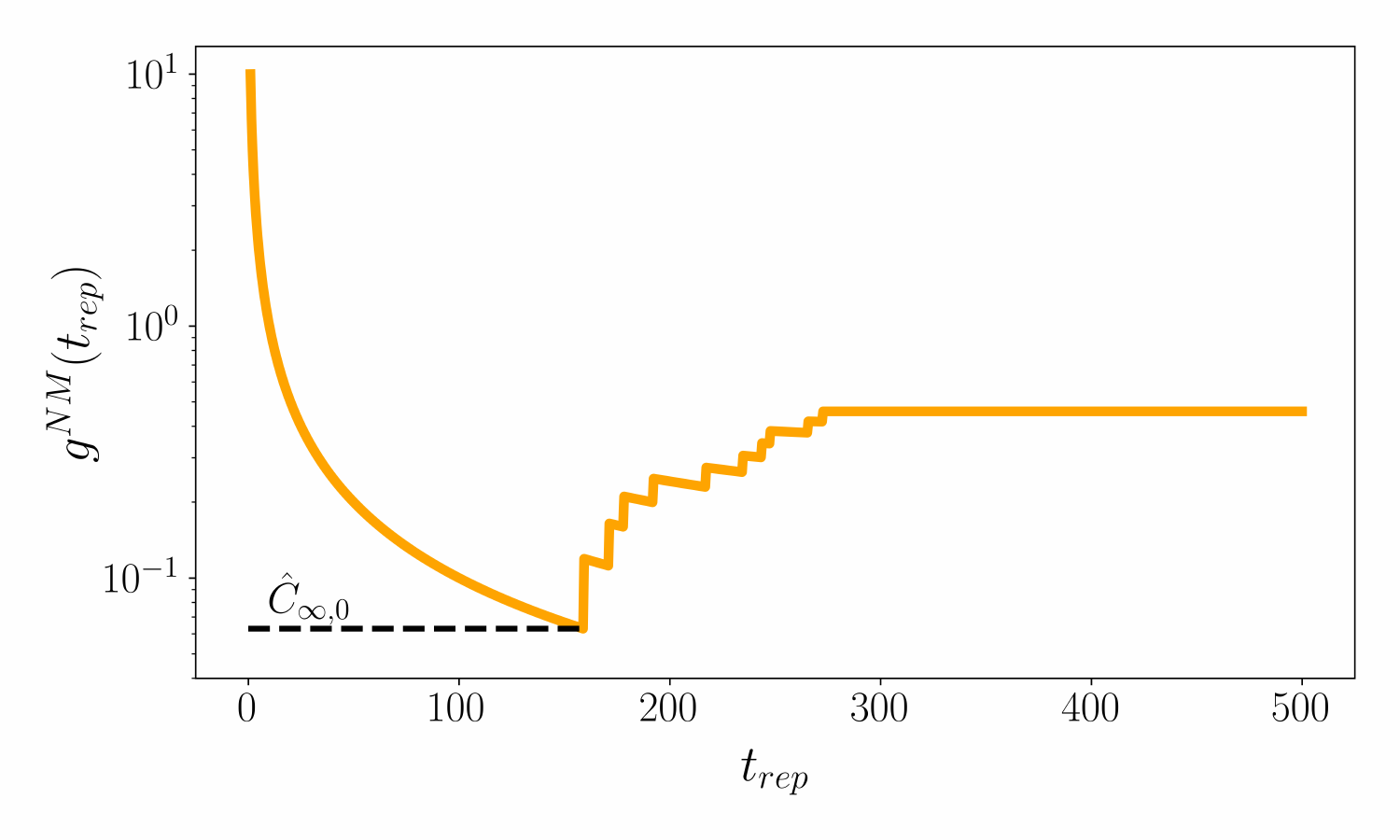}
        \caption{$doa$}
        \label{subfig:doa_ectr_vs_trep}
    \end{subfigure}
    \hfill
    \begin{subfigure}{0.49\textwidth}
        \includegraphics[width=\textwidth]{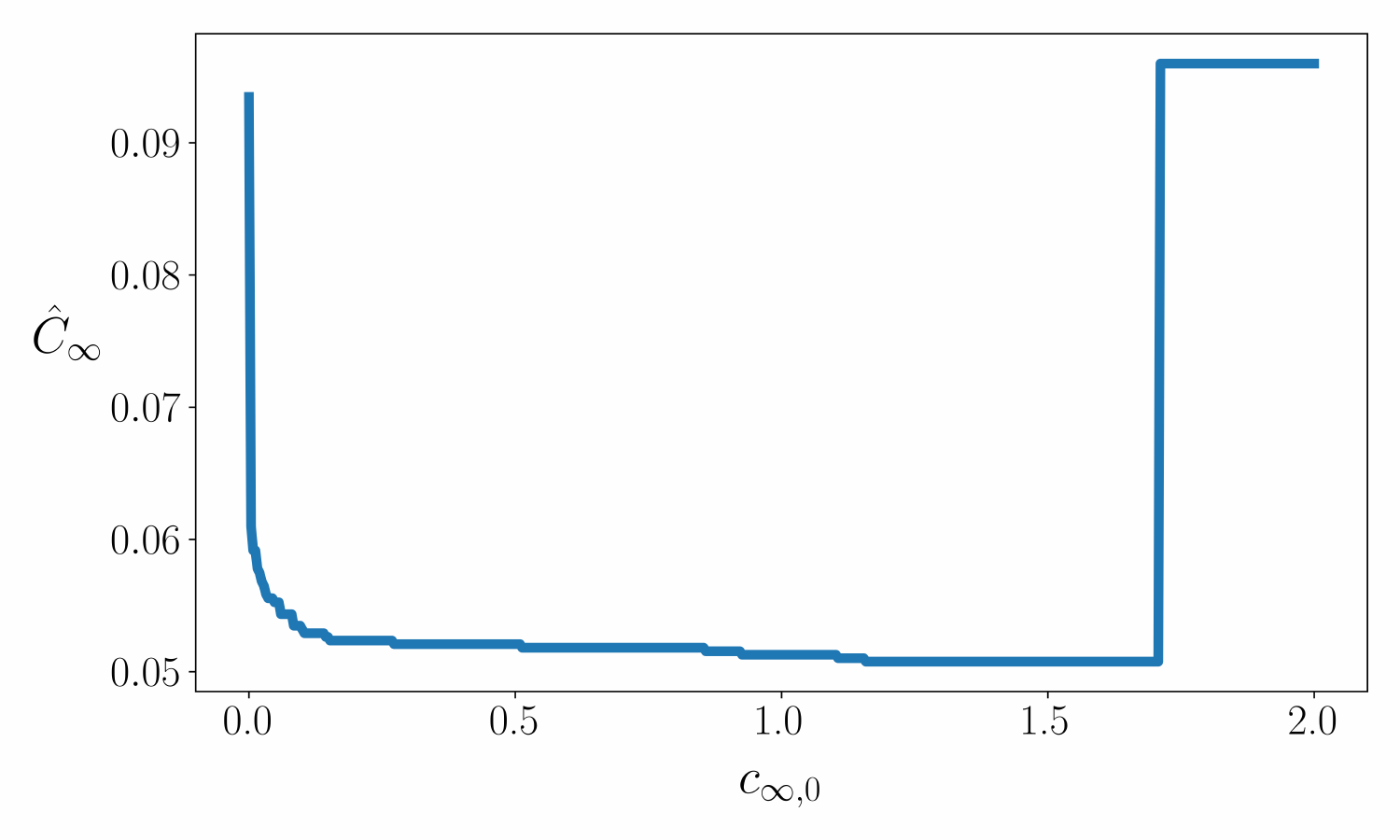}
        \caption{$doa^*$}
        \label{subfig:doa_ectr_vs_ectr}
    \end{subfigure}
    \caption{One example of the objective functions in the optimization of the heuristic parameters when only 10 failure samples are available. The cost ratio here is $\frac{c_p}{c_c}=0.1$. The objective function for $rh1^*$ (top) and $doa^*$ (bottom right) is the final policy performance ($\hat{C}_{\infty}$), whereas for $doa$ it is the initial cost ratio estimate ($c_{\infty,0}$).}
    \label{fig:obj_function_rh1_vs_doa}
\end{figure}

In \Cref{fig:obj_function_rh1_vs_doa}, we examine the optimization for replacement heuristics $rh1^*$, $doa$ and, $doa^*$ in the case where only $10$ sample histories are available. $rh1^*$ and $doa^*$ perform the optimization directly over the policy performance ($\hat{C}_{\infty}$), whereas for $doa$ we perform an optimization over $g(\trep)$ described in \Cref{eq:obj_func_rep_nm} for finding the cost rate estimate $c_{\infty,0}$, which is then subsequently used in the final policy. The differences between these two approaches have substantial consequences. For low amounts of failure samples, the objective function for $rh1^*$ and $doa^*$ is constant over a wide range of $p_\mathrm{thres}$ values; see the plateaus in \Cref{subfig:rh1_perf_vs_pthres,subfig:doa_ectr_vs_ectr}. The locations of the vertical jumps, as well as the final $\hat{C}_{\infty}$ values, depend on the specific set of samples. 
Because of these plateaus, there is no unique solution to the optimization problem. We dealt with this by generating uniform samples over the lowest plateau. But clearly, this is an issue when optimizing a heuristic with respect to the policy performance with limited data and impacts the accuracy of the solution.

In addition, finding the lowest plateau might be nontrivial, as gradient-based optimizers cannot be used and instead $\hat{C}_{\infty}$ has to be evaluated over the whole range of possible $p_{\mathrm{thres}}$ values. For the one parameter in our simulation setting, this is not a problem, but this could pose some issues for more parameters and more expensive policy evaluations. 
By contrast, with $doa$, the objective function for the parameter 
optimization for an initial cost rate estimate is better behaved. Although there are local minima where the optimizer can get stuck, the function is monotonically decreasing between jumps, which helps the optimizer. In our experience, even with these local minima, the optimizer usually finds the global minimum over the whole range of $\frac{c_p}{c_c}$ values. One can even avoid the optimization problem by choosing another initialization option for $c_{\infty,0}$, as presented in \cite{kamariotis2024metric}; here, one simply uses $\frac{c_p}{\hat{\mu}_{\Tf}}$, where $\hat{\mu}_{\Tf}$ is the mean failure time computed over the set of available training samples. We have found that with this second initialization option, the heuristic's performance is slightly worse compared to the first option, but the results are almost indistinguishable. 

\begin{figure}[H]
    \centering
    \includegraphics[width=0.93\textwidth]{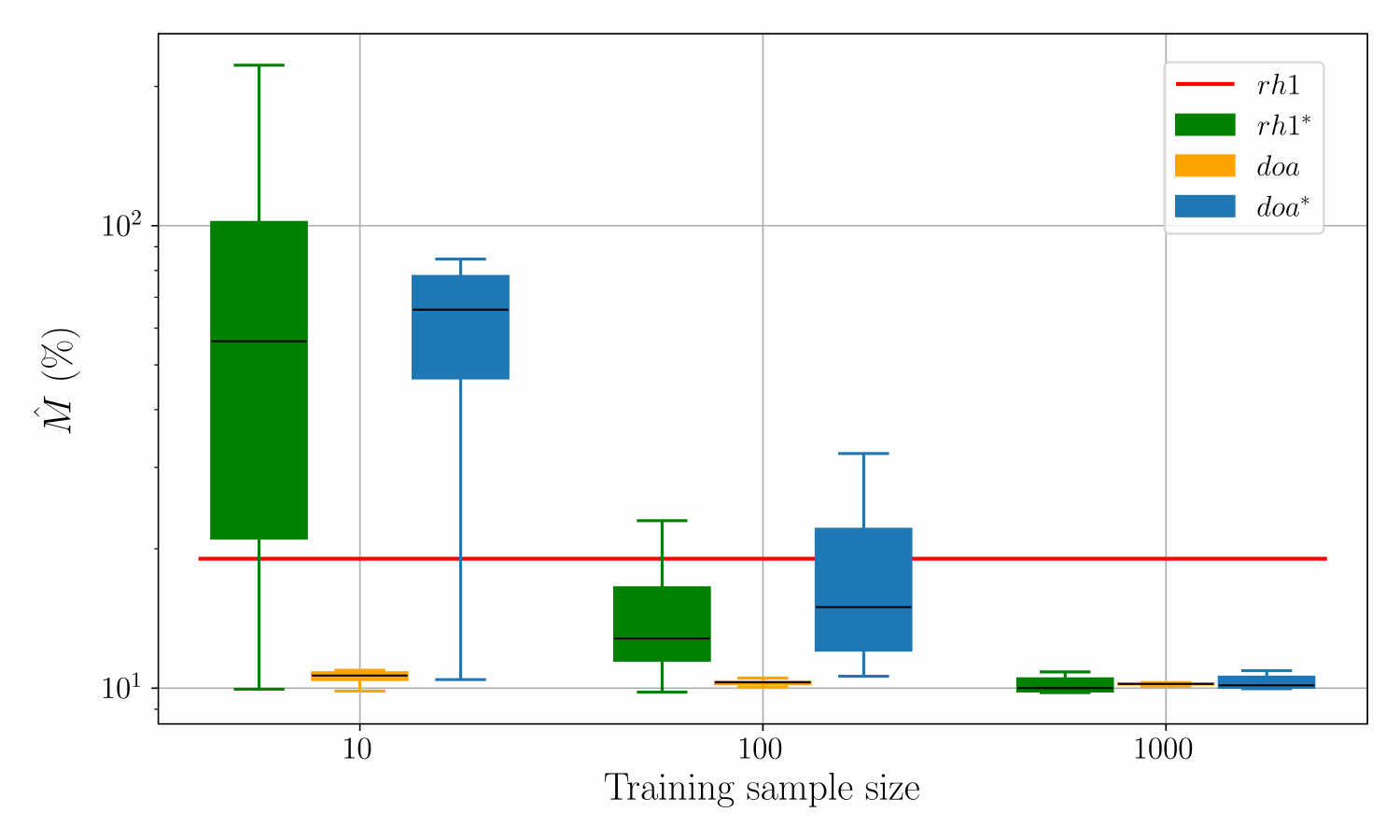}
    \caption{Boxplot comparison of our proposed $doa$ policy with benchmark heuristics $rh1$ and $rh1^*$ for $c_p/c_c=0.1$. For the optimization problem, training sample sizes of [10, 100, 1000] were used; evaluation of the final policies was performed with $10^6$ samples, and the investigation was repeated 20 times, with the respective medians highlighted in black.}
    \label{fig:r_comparison_low_data}
\end{figure}

Ultimately, the key question is how the optimized policies perform with limited training samples.
Therefore, we test the performance of $rh1$ and $doa$ for various training sample sizes for a fixed cost ratio $\frac{c_p}{c_c}=0.1$. The results for 20 independent runs are shown as boxplots in \Cref{fig:r_comparison_low_data}. They demonstrate that if one optimizes $p_\mathrm{thres}$ or $c_{\infty}$ with only small amounts of training data, the resulting policy can have substantially worse performance than with the default parameter values $p_\mathrm{thres}=\frac{c_p}{c_c}$ and $c_{\infty,0}$. This is classical overfitting. One can also see that our proposed $doa$ heuristic is very robust against overfitting, and yields excellent results even in low data regimes. 

%

%
%
%%%%%%%%%%%%%%%%%%%%%%%%%%%%%%%%%%%%%%%%%%%%%%%%%%%%%%%%%%%%%%%%%
%%%%%%%%%%%%%%%%%%%%%%%%%%%%%%%%%%%%%%%%%%%%%%%%%%%%%%%%%%%%%%%%%
%
%
%
\subsection{Ordering policy performance}
\label{subsec:ordering_performance}
\begin{figure}[H]
    \centering
    \includegraphics[width=0.93\textwidth]{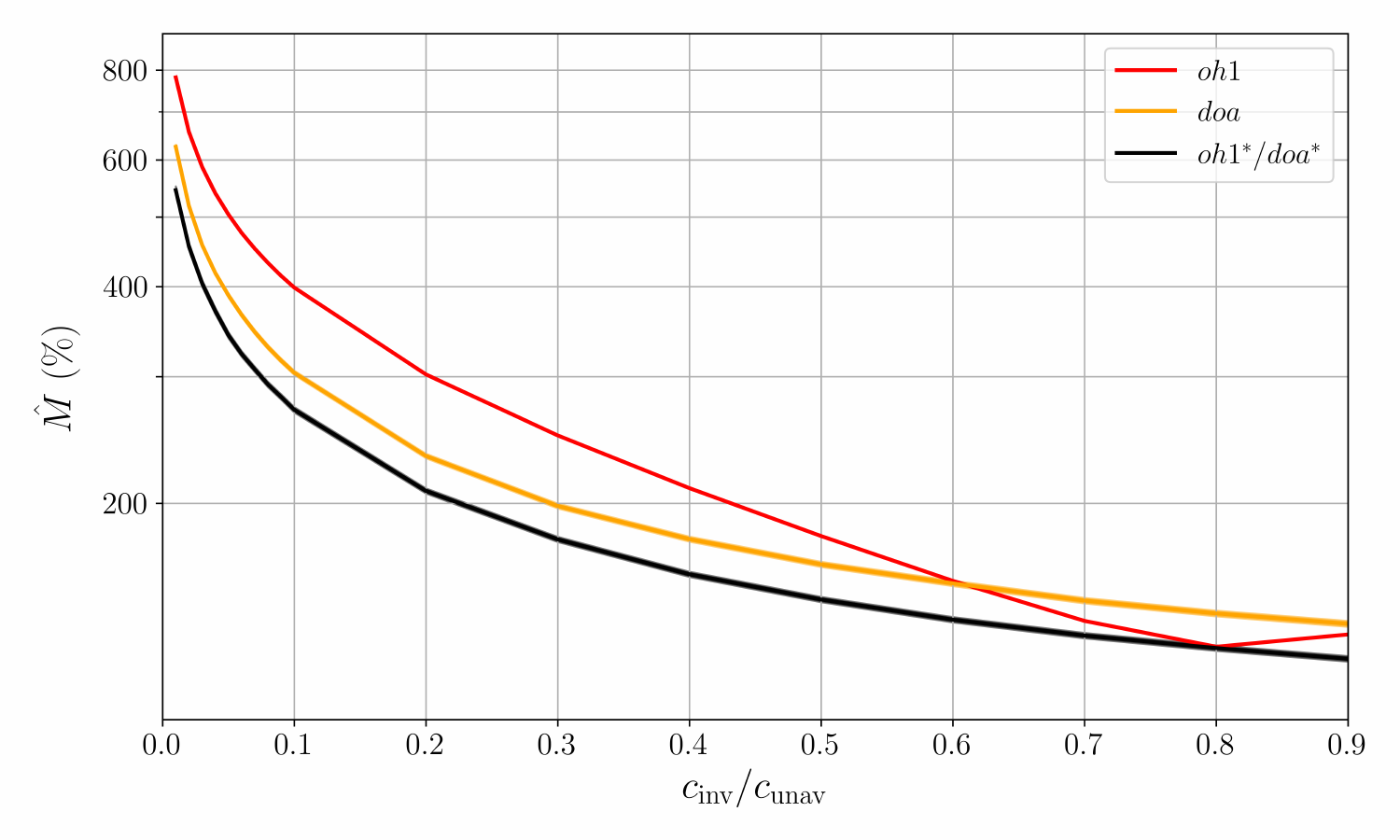}
    \caption{Comparison of our proposed $doa$ ordering policy (\Cref{eq:a_ord_k_h4}) with the simple heuristic ordering policy (\Cref{eq:a_ord_k_h1}). Since the optimized version of both policies yields the same performance, they are represented as a single black line. The $oh1$ and $doa$ policies were evaluated with $10^6$ and $10^5$ MC samples, respectively. The figure shows the mean of each policy as a solid line, as well as the corresponding 95\% credible intervals.}
    \label{fig:oh_perf_overview}
\end{figure}

We perform a similar investigation as in \Cref{subsec:replacement_performance} to compare our proposed $doa$ preventive ordering heuristic with the simple threshold policy $oh1$. We again include for both heuristics the default versions (without parameter optimization) as well as their optimized counterparts, denoted as $oh1^*$ and $doa^*$, respectively. The default threshold for $oh1$ is chosen as $p_{\mathrm{thres}}=\cinv / \cunav$. The results for time step size $\Delta t=10$ and ordering lead time $l=11$ in function of the cost ratio $c_\mathrm{inv}/c_\mathrm{unav}$ are depicted in \Cref{fig:oh_perf_overview}. They are quite similar to the results of the preventive replacement setting in \Cref{fig:rh_perf_overview}, except that the differences between the policies are now much larger (in the order of hundreds of percent), which indicates a more difficult problem. 

Again, the default version of our proposed $doa$ ordering heuristic achieves significantly lower cost ratios compared to the default version of $oh1$ over most of the investigated range. In the range $c_\mathrm{inv}/c_\mathrm{unav}\in[0.62,92]$, $oh1$ with $p_\mathrm{thres}=c_\mathrm{inv}/c_\mathrm{unav}$ achieves better performance than our proposed method. 
This better performance is a coincidence rather than a systematic benefit and happens because for a cost ratio of $c_\mathrm{inv}/c_\mathrm{unav}=0.8$, the optimal threshold is also close to $0.8$. Furthermore, typical decision settings are located on the left side of the graph, i.e., for lower cost ratios, where the $doa$ heuristic is superior. 

Optimizing over the free parameter in each decision heuristic ($p_\mathrm{thres}$ for $oh1$ and $c_{\infty}$ for $doa$), i.e., data-driven compensation of the process cost assumptions, results in the best performances. The curves of both optimized policies are visually indistinguishable.

\begin{figure}[H]
    \centering
    \includegraphics[width=0.93\textwidth]{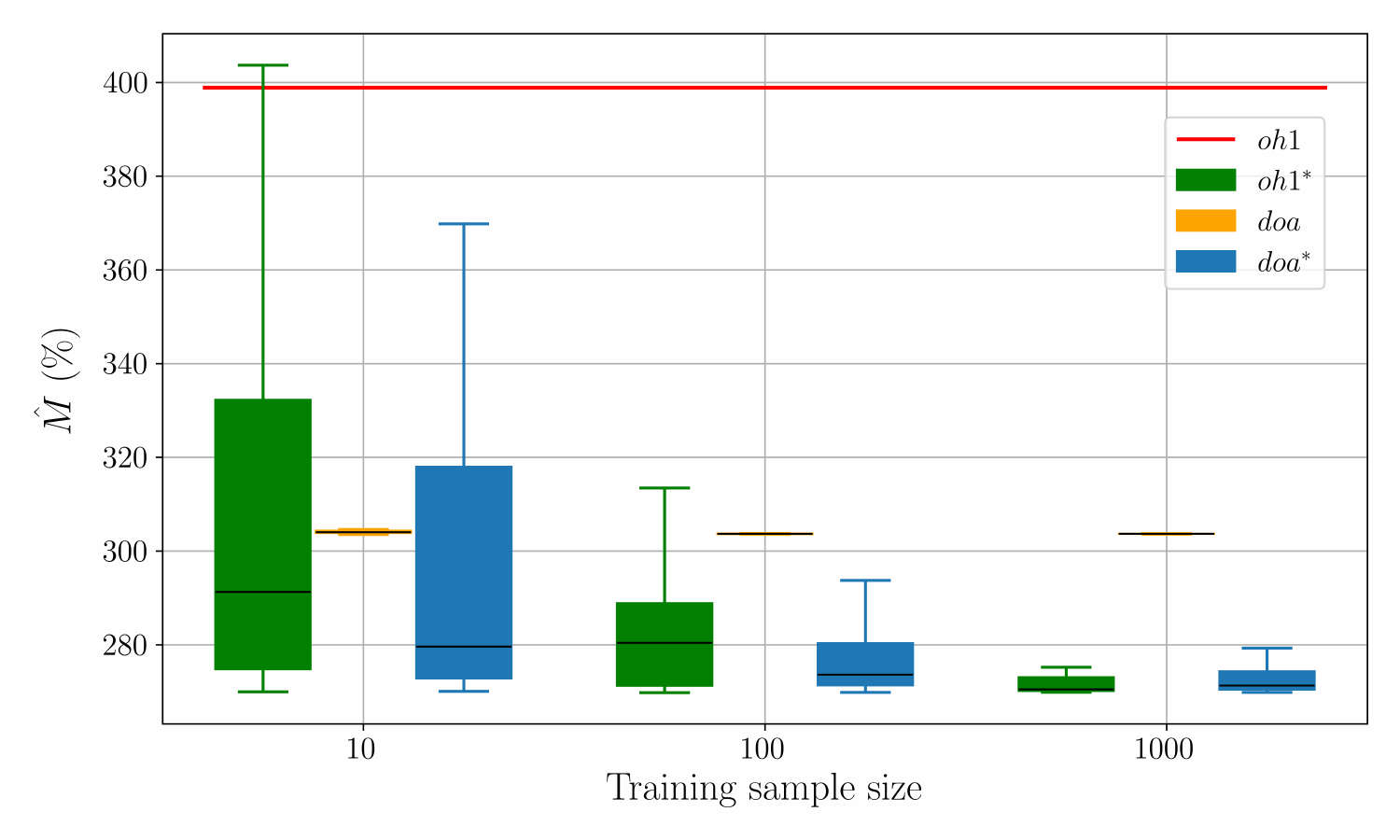}
    \caption{Boxplot comparison of our proposed $doa$ policy with $oh1$ with optimized $p_\mathrm{thres}$ for $c_\mathrm{inv}/c_{\mathrm{unav}}=0.1$. For the optimization problem, training sample sizes of [10, 100, 1000] were used; evaluation of all final policies was performed with $10^5$ samples, and the investigation was repeated 40 times, with the respective medians highlighted in black.}
    \label{fig:comparison_o_low_data}
\end{figure}

Similarly to \Cref{subsec:replacement_performance}, we also investigate the robustness of these heuristics in the low data regime; the results of which are depicted in \Cref{fig:comparison_o_low_data}. They show that $doa$ without optimization is very robust against overfitting, even in low data regimes. However, due to the larger performance gap, also noticeable in \Cref{fig:oh_perf_overview}, the optimized $oh1^*$ policy, achieves a better performance in most cases already for 100 samples and a better median performance already for 10 samples. Overall, the best-performing policy is $doa^*$. While it is less robust than $doa$, it does have significantly lower $\hat{M}$ and is more robust than the $oh1^*$ policy. Nevertheless, in the 10-sample regime, one might also opt for $doa$.

%

%
%
%%%%%%%%%%%%%%%%%%%%%%%%%%%%%%%%%%%%%%%%%%%%%%%%%%%%%%%%%%%%%%%%%
%%%%%%%%%%%%%%%%%%%%%%%%%%%%%%%%%%%%%%%%%%%%%%%%%%%%%%%%%%%%%%%%%
%%%%%%%%%%%%%%%%%%%%%%%%%%%%%%%%%%%%%%%%%%%%%%%%%%%%%%%%%%%%%%%%%
%
%
\section{Case study: preventive ordering of turbofan engines}
\label{sec:validation_on_cmapss_dataset}

In this section, we investigate our proposed method in a more realistic PdM setting. We utilize the C-MAPSS dataset\footnote{C-MAPSS and other prominently used prognostics datasets are available at the NASA Prognostics Center of Excellence Data Set Repository: \url{https://www.nasa.gov/intelligent-systems-division/discovery-and-systems-health/pcoe/pcoe-data-set-repository/}}, a well-known PHM benchmark problem concerning the degradation simulation of turbofan engines under different flight conditions \cite{saxena2008cmapss}.

The data of interest consists of 24 time series, comprising 21 sensor readings and 3 engine operation mode variables. We focus on the training set of FD001, one out of four data subsets, which contains run-to-failure trajectories of different engines under varying initial conditions, degradation parameters, and noise levels. The FD001 training set contains a total of 100 run-to-failure trajectories, which we use for our investigations. 
For a more detailed description of the dataset, the reader is referred to \cite{saxena2008cmapss}.

For our investigation, we split the original training set into a training set and a test set. The training set is used to train a prognostic model as well as to perform potential initialization or optimization of a decision heuristic's free parameter. We suspect that this double use of the training set exacerbates the overfitting issue observed in \Cref{fig:r_comparison_low_data,fig:comparison_o_low_data}.\footnote{Note that this “double use” of the training set relates to a practical finite-sample overfitting effect observed when both the prognostic model training and the parameter optimization of the heuristic policy rely on the same limited training data, rather than a methodological shortcoming of our evaluation protocol.} The final performance of the total PdM bundle, consisting of the prognostic model and the decision heuristic, is then evaluated on the test set. To show the variability of the whole training pipeline, we perform repeated random sub-sampling, which is also known as Monte Carlo cross-validation (MCCV). The full evaluation process is depicted in \Cref{fig:cmapss_evaluation_process}.

\evaldiagram

As a prognostic model, we choose a Long Short-Term Memory (LSTM)-based neural network (NN), which is a standard choice for time series prediction tasks. The specific architecture follows the RUL classifier introduced in \cite{nguyen2019new}, i.e., the NN outputs a class label of the estimated RUL with a corresponding probability. For our ordering setting, the classes are chosen as i) $RUL\leq l$, ii) $RUL\leq l + \Delta t$, and iii) $RUL > l+\Delta t$. From the individual class probabilities, the full RUL-PDF required for the $doa$ heuristic can be obtained by fitting a certain distribution type: for this work, following \cite{kamariotis2024metric}, the lognormal distribution is chosen. A more detailed description of the employed architecture, the configuration, as well as the hyperparameter setup is given in \Cref{sec:app_nn_specs}. 

\begin{figure}[H]
    \centering
    \includegraphics[width=0.93\textwidth]{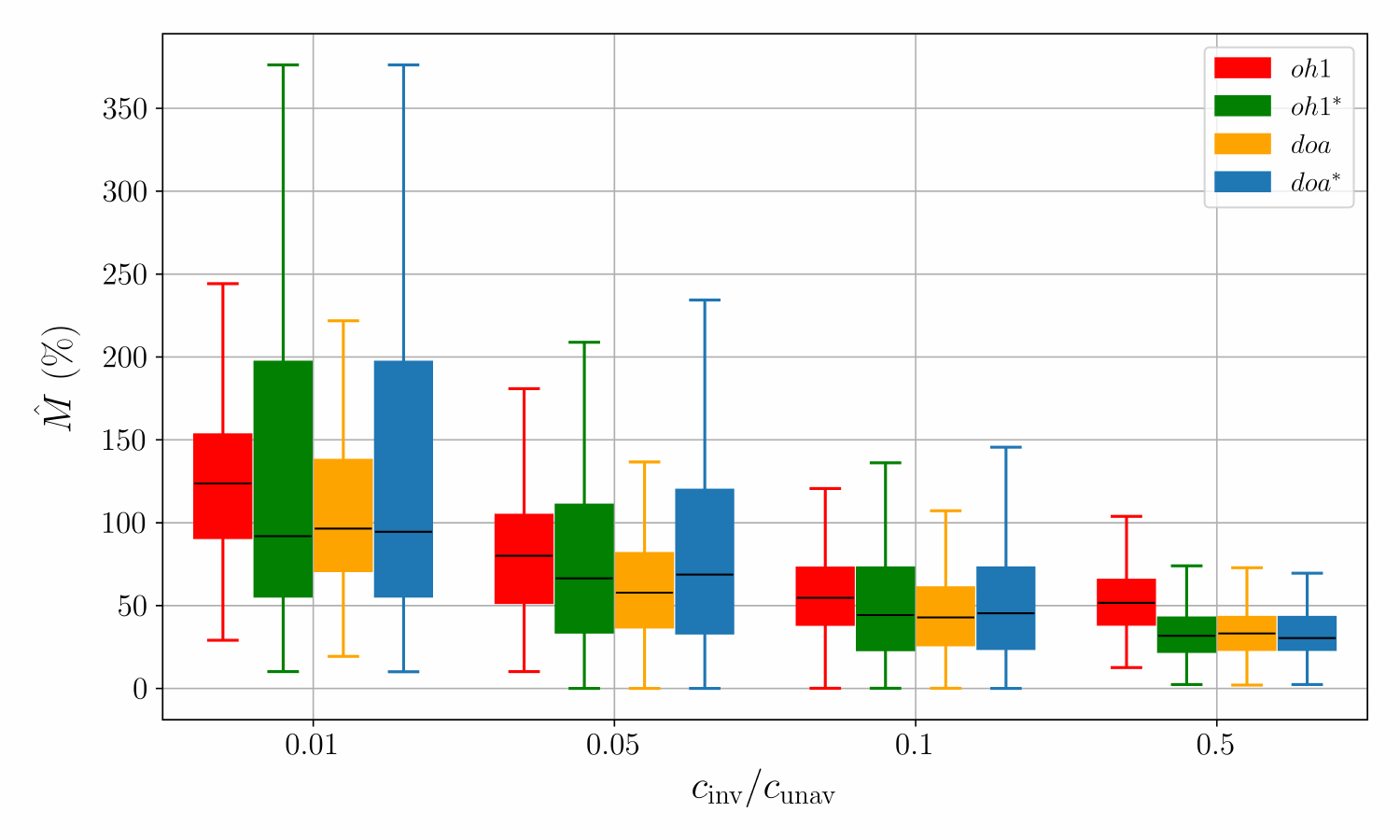}
    \caption{Boxplot comparison of our proposed $doa/doa^*$ policies with benchmark heuristics $rh1/rh1^*$ for $c_p/c_c=[0.01, 0.05, 0.1, 0.5]$. 80 of the available 100 trajectories were used for training of the prognostic model and for initialization/optimization of the heuristics' free parameter; each policy's performance is evaluated on the remaining 20 trajectories. The investigation is repeated 100 times, with the respective medians highlighted in black.}
    \label{fig:comp_cmapss_olc}
\end{figure}

We focus the investigation on the ordering decision setting, because the prognostics algorithm achieves such a good performance for this data set that the replacement decision setting results in $\hat{M}$ values close to $0$, even with a non-optimized heuristic, see \cite{kamariotis2024metric}.
With the presented evaluation pipeline, we compare the proposed $doa/doa^*$ ordering heuristics with the benchmark policies $oh1/oh1^*$. The results are shown for multiple inventory to unavailability cost ratios in \Cref{fig:comp_cmapss_olc}. The box plots reflect the variability among the 100 MCCV repetitions. In each of these, we employ the random 80/20 train+validation/test split of the training FD001 dataset.

The resulting performances show larger variability for all investigated policies and cost ratios than with the RUL simulator of \Cref{sec:virtual_RUL_simulator}. We suspect that this results from the combined effect of multiple factors, in particular the small test data size and uncalibrated RUL predictions from the prognostic model. Note that the shown boxplots in \Cref{fig:comp_cmapss_olc} are based on correlated samples obtained with MCCV; thus, the true performance variances are likely even larger.

We observe that our proposed $doa$ heuristic yields the best or close to the best performance over all investigated $\cinv/\cunav$ values. $rh1$ consistently shows the worst median performance over all cost ratios. As expected, optimizing the decision heuristic on the limited training dataset of 80 trajectories leads to overfitting, which results in the increased performance variations of $rh1^*$ and $doa^*$. This phenomenon is amplified in regions with low cost ratios. Furthermore, $rh1^*$ and $doa^*$ show similar behaviour with comparable median performances and performance variations across the investigated cost ratios.

%

%
%
%%%%%%%%%%%%%%%%%%%%%%%%%%%%%%%%%%%%%%%%%%%%%%%%%%%%%%%%%%%%%%%%%
%%%%%%%%%%%%%%%%%%%%%%%%%%%%%%%%%%%%%%%%%%%%%%%%%%%%%%%%%%%%%%%%%
%%%%%%%%%%%%%%%%%%%%%%%%%%%%%%%%%%%%%%%%%%%%%%%%%%%%%%%%%%%%%%%%%
%
%
\section{Discussion}
\label{sec:discussion}
This paper proposes a novel method to derive heuristic maintenance policies for PHM. We call these maintenance policies \textit{heuristics} because they do not solve the full sequential decision optimization problem that is the maintenance process. We note that most currently existing maintenance policies proposed for PHM in the literature are of this kind. This is due to the limited data availability, which hinders the training of more sophisticated maintenance policies. An added advantage of heuristic policies is their explainability, which is important in practice. 

In contrast to existing PHM maintenance policies, which were mostly derived without much theoretical justification, our proposed method does consider the full decision sequence explicitly. In order to end up with a heuristic policy, it simplifies the full sequential decision process by only considering and comparing the options of performing maintenance at the current and next opportunity. This online forward search with depth $d=1$ is enhanced with a renewal theory-based offline estimation of the future costs (via the use of the long-running maintenance cost per unit time). Thus, combining online and offline approaches, our proposed method for deriving predictive maintenance policies falls in the category of \emph{hybrid planning} \cite{kochenderfer2022algorithms}.

The use of the one-step decision tree assumes a certain monotonicity of the underlying problem. To show this, let again denote $c_{M@t_k}$ as the expected costs of performing maintenance at time step $t_k$.\footnote{For example, for preventive replacement, $c_{M@t_k}$ and $c_{M@t_{k+1}}$ are given in \Cref{eq:h4_c_PR_tk,eq:h4_c_PR_tk+1}, respectively.} If one were to consider maintenance decision options also at all future time steps $t_{k+\tau},~\tau > 1$, the optimal maintenance decision would correspond to the one derived from the one-step tree in all cases, except when $c_{M@t_{k+\tau}}<c_{M@t_k}<c_{M@t_{k+1}}$. In this case, the one-step tree yields the decision to maintain immediately, whereas the multi-step tree would yield the decision to do nothing. 
The one-step tree ignores that this case can occur, and thus assumes a monotonic cost structure of the maintenance actions over time. 
If decision-making is solely based on the latest RUL-PDF (i.e., no further information is extracted from the environment dynamics or from past RUL-PDFs), we believe that the one-step tree yields the optimal decision if the non-monotonic cost case does not occur. 

Therefore, the one-step simplification of the decision process leads to policies that tend to perform a maintenance action too early. The representation of the costs of \emph{maintain later} by $c_{M@t_{k+1}}$ undersells the advantages of waiting with a certain maintenance action. This results from the fact that one does not necessarily perform maintenance at $t_{k+1}$, but instead waits for more information and decides anew. The value of the option to perpetually wait and extend the life of the component for $2,3,...,n$ cycles is not included in the decision rule. 
Nevertheless, the achieved performances of the derived policies in \Cref{sec:Numerical_investigations_RUL_simulator,sec:validation_on_cmapss_dataset} show that the one-step assumption results in good policies for various cost and decision settings. Still, there is potential in further reducing the gap to the optimized policies, which correct for the conservatism of the proposed $doa$ policies using data. 

In future works, different approaches could be pursued to account for the value of the option to wait beyond the next decision time step. One approach could be to extend the one-step decision trees to include more future steps $t_{k+\tau},~\tau > 1$ (i.e., forward search with a depth $d>1$).\footnote{A challenge here is that the assignment of the probabilities and associated expected costs of the subbranches would have to be made based on future RUL-PDFs $f_{RUL_{k+\tau}}$, i.e., on information that only becomes available in the future. Correctly quantifying these probabilities and costs is nontrivial, but using the latest RUL-PDF $f_{RUL_k}$ as an approximation of the true future distributions could be a first avenue to explore.} Another approach would be to keep the one-step formulation of
\emph{maintain now vs. maintain later}, but to use cost expressions that quantify the value of waiting and extending the life of the component for multiple time steps. Examples for such formulations could be the switch to age-based replacement (with the RUL-PDF as the time-to-failure PDF), or the assumption of replacement at the expected value (with $c_{\infty} \cdot \mathbb{E}\left[{RUL}\right]$ as expected cost savings).

The proposed policies are computationally more expensive than the simple threshold-based benchmark policies. The proposed policies require computation of consequences for each branch of the underlying tree, with the main cost drivers being the evaluation of the cumulative probabilities, $\Pr(a<RUL\leq b)$, and conditional expectations, $\mathbb{E}[RUL_k \mid a < RUL\leq b]$. They can either be computed analytically (as in this work) or, in the general case, via 1-D numerical integration. For the single-component decision settings of this paper, the respective decision trees have at most four branches. Comparing the policy expressions for replacement (\Cref{eq:rep_threshold,eq:a_rep_k_h4}) and ordering (\Cref{eq:a_ord_k_h1,eq:a_ord_k_h4}), one can see that these cost-driving terms occur only once in the threshold-based policies, but $2-3$ times in the $doa$ policies. Hence, evaluating the $doa$ policies requires approximately 2–3 times the computation time, which is not critical. Accounting for the computational cost of the parameter optimization for $rh1^*/oh1^*$, the threshold-based policies are approximately 10–15 times more computationally expensive than the $doa$ policies.

Optimized policies still yield the best performance when large amounts of training data are available. Interestingly, it does not seem to matter how the policy was originally formulated; optimization leads to the same final performance for these two relatively simple decision settings of preventive replacement and preventive ordering. Here, future work could explore the true performance limit by comparison with a partially observable Markov Decision Process based solution, assuming full knowledge of the environment \cite[see, e.g.,][]{kochenderfer2015decision}. 

Another advantage of parameter optimization is that it can also compensate for potential inadequacies of the prognostic model. Policies with initialized parameters, among them also the proposed $doa$ policies, do not incorporate optimization and hence require a more reliable prognostic model to achieve a good maintenance performance.

The major drawback of optimized policies is their heavy overfitting in low data regimes, where one can end up with worse performance than with reference values or just simple age-based policies without monitoring. Here, our proposed $doa$ policies with simple parameter initialization provide a remedy. In both simulated settings as well as on the C-MAPSS dataset, our default proposed heuristics are robust against overfitting and thus can be reliably applied even in low data domains. Future work could investigate the application of the presented method to more complicated decision settings, e.g., ordering combined with replacement, or non-deterministic settings for $c_p, c_c, \cinv, \cunav, l$. Of particular interest is the formulation of heuristics for multi-component systems or fleets. An envisioned end result for practitioners could take the form of a look-up table, i.e., a guide on which heuristic to choose for the decision setting at hand. Note that for an M-component system, the number of branches in the $doa$ decision tree grows as $2^{2M}$. Hence, extending the $doa$ policies to multi-component systems poses a greater computational challenge. We leave this for further research.

In this work, we have neglected discounting of future incurred costs, and used $c_{\infty}$ to consider the costs of the underlying renewal-reward process in the $doa$ decision rules (and $rh2$). For applications in which the overall lifetime of components is in the order of up to a few years, this is a reasonable assumption. However, to generalize the applicability to components with longer lifetimes, discounting should be considered.
\citet{van2008renewal} generalize renewal-reward processes to include discounting with respect to any discount function that is non-increasing and monotonic over time. Here, one can formulate an expected equivalent average cost (EEAC) per unit time as a discounted version of $c_{\infty}$. For instance, for exponential discounting described by $e^{-\rho t}$ with discount rate $\rho$ over time $t$, EEAC (in the limit $t\rightarrow\infty$) can be obtained as \cite{van2008renewal}:
\begin{equation}
    \label{eq:eeac_exp}
    EEAC = \frac{
    \rho\mathbb{E}\left[e^{-\rho T_{\mathrm{lc}}}C\right]}{
    1-\mathbb{E}\left[e^{-\rho T_{\mathrm{lc}}}\right]}
\end{equation}
Hence, to include discounting, one could potentially use EEAC instead of $c_{\infty}$ to account for the costs of the underlying renewal-reward process, though it remains to be verified whether the heuristic policies derived in this paper would require modification.

%
%
%%%%%%%%%%%%%%%%%%%%%%%%%%%%%%%%%%%%%%%%%%%%%%%%%%%%%%%%%%%%%%%%%
%%%%%%%%%%%%%%%%%%%%%%%%%%%%%%%%%%%%%%%%%%%%%%%%%%%%%%%%%%%%%%%%%
%%%%%%%%%%%%%%%%%%%%%%%%%%%%%%%%%%%%%%%%%%%%%%%%%%%%%%%%%%%%%%%%%
%
%
\section{Concluding remarks}
\label{sec:conclusion}
We present a novel and general method of deriving heuristic policies for predictive maintenance. The policies are based on a reduced decision tree to model the available decision options and their consequences. The probability of reaching a certain branch is computed with the RUL-PDF obtained from a prognostic model. The expected costs of each action are then computed as the sum of the maintenance action and failure costs, which are directly attributable to the monitored component, as well as the costs of the underlying renewal-reward process. In this way, the approach also accounts for the gain from a potential life extension of a component.

We applied the proposed method to two predictive maintenance settings, namely preventive replacement and preventive ordering. We investigated the performance of our proposed $doa$ policies, together with commonly used benchmarks, on a simulated dataset as well as the C-MAPSS dataset. Our results show that the $doa$ policies achieve lower cost rates than the presented benchmarks. In the case of preventive replacement, our proposed heuristic is almost as good as an optimized heuristic. 
When a large corpus of training data is available, one can compensate a ``bad'' heuristic by parametrization of the decision policy, with subsequent optimization over the parameter(s) with respect to the policy's final performance. However, in the context of low training data availability, this leads to overfitting. Here, the proposed $doa$ policies have shown to be robust against overfitting in all investigated decision settings and datasets. 
If one uses a parametrized version of our proposed heuristic, optimization still leads to overfitting, with more or less variance compared to benchmark policies, depending on the specific decision setting.

On the one hand, this work highlights the importance of data, and how large datasets can lead to great performance through optimization. On the other hand, this work also accentuates that a lot of performance can be gained by taking a step back and performing a more theory-driven investigation. This is essential for applications with a limited amount of run-to-failure data, which is the case for many application domains. The better and more robust policies resulting from a theory-driven approach constitute another step towards optimal use of the data at hand and, ultimately, the improvement of efficiency and safety of engineering systems.

%
%
%%%%%%%%%%%%%%%%%%%%%%%%%%%%%%%%%%%%%%%%%%%%%%%%%%%%%
%%%%%%%%%%%%%%%%%%%%%%%%%%%%%%%%%%%%%%%%%%%%%%%%%%%%%
%%%%%%%%%%%%%%%%%%%%%%%%%%%%%%%%%%%%%%%%%%%%%%%%%%%%%
%
%
\section*{Acknowledgement}
This work was supported by the German Federal Ministry for Economic Affairs and Climate Action (BMWK) through the aviation research program LUFO VI-3 in the project BIG-ROHU.

\bibliographystyle{unsrtnat}
\bibliography{references}  %%% Uncomment this line and comment out the ``thebibliography'' section below to use the external .bib file (using bibtex) .

\newpage
\appendix
\counterwithin{figure}{section}
\counterwithin{table}{section}

\label{appendix}

%
%
%
%%%%%%%%%%%%%%%%%%%%%%%%%%%%%%%%%%%%%%%%%%%%%%%%%%%%%%%
%%%%%%%%%%%%%%%%%%%%%%%%%%%%%%%%%%%%%%%%%%%%%%%%%%%%%%%
%%%%%%%%%%%%%%%%%%%%%%%%%%%%%%%%%%%%%%%%%%%%%%%%%%%%%%%
%
%
\section{Derivations}
For the derivation of the final $doa$ policies, we omit the subscript $k$ for ease of notation.

\subsection{Derivation of doa replacement heuristic}
\label{sec:app_derivation_doa_final_costs}
Starting with the decision criterion in \Cref{eq:a_rep_k_h4_gen}, we can expand all individual terms:
\begin{equation}
\begin{aligned}
    \cancel{c_p} + \cancel{c_{\infty} \cdot r} <~&
    \Pf \cdot c_c + \cancel{\Pf \cdot c_{\infty} \cdot r} - \Pf \cdot c_{\infty} \cdot \EE{RUL \mid RUL < \Delta t} + \\
    & \cancel{c_p} + \cancel{c_{\infty} \cdot r} - c_{\infty} \cdot \Delta t  - %\\
    %&
    \Pf \cdot c_p - \cancel{\Pf \cdot c_{\infty} \cdot r} + \Pf \cdot c_{\infty} \cdot \Delta t,
\end{aligned}
\end{equation}
which yields through reordering
\begin{equation}
\label{eq:doa_crit_final}
    c_{\infty} \cdot \Delta t < \Pf \left[c_c - c_p + c_{\infty}\left(\Delta t -  \EE{RUL \mid RUL < \Delta t}\right) \right].
\end{equation}

%

%
%
%
%%%%%%%%%%%%%%%%%%%%%%%%%%%%%%%%%%%%%%%%%%%%%%%%%%%%%%%
%%%%%%%%%%%%%%%%%%%%%%%%%%%%%%%%%%%%%%%%%%%%%%%%%%%%%%%
%%%%%%%%%%%%%%%%%%%%%%%%%%%%%%%%%%%%%%%%%%%%%%%%%%%%%%%
%
%
\subsection{Derivation of doa ordering heuristic}
\label{sec:app_derivation_oh4_final_costs}
\newcommand{\ctone}{red}
\newcommand{\cttwo}{blue}
\newcommand{\ctthree}{green!70!black}
\newcommand{\ctfour}{orange}
\newcommand{\ctfive}{violet}
\newcommand{\ctsix}{brown}
\newcommand{\ctseven}{gray}
\newcommand{\cteight}{magenta}
\newcommand{\ctnine}{black}
\newcommand{\os}[3]{\overset{\textcolor{#1}{#2}}{\vphantom{I}#3}}
Starting with the decision criterion in \Cref{eq:a_ord_k_h4_gen}, we investigate all individual terms separately, to see which ones cancel out:
\begin{align}
\label{eq:app_cord_at_tk_total}
\begin{split}
    c_{Ord@t_k} =~& 
    \Pro{RUL \leq l} \cdot \bigg[\bigg( 
    \os{\ctone}{1}{l} - 
    \os{\cttwo}{2}{\EE{RUL \mid RUL \leq l}}
    \bigg) \cdot \cunav 
    + \bigg(
    \os{\ctthree}{3}{r} - 
    \os{\ctfour}{4}{l}
    \bigg) \cdot c_{\infty}
    \bigg] + \\[0.35em]
    & \Pro{RUL > l} \cdot \bigg[\bigg( 
    \os{\ctfive}{5}{\EE{RUL \mid RUL > l}} - 
    \os{\ctsix}{6}{l}
    \bigg) \cdot \cinv
    + \bigg(
    \os{\ctthree}{3}{r} - \os{\ctseven}{7}{\EE{RUL \mid RUL > l}}
    \bigg) \cdot c_{\infty}
    \bigg]
\end{split}
\end{align}
\begin{align}
\begin{split}
\label{eq:app_cord_at_tk1_total}
    c_{Ord@t_{k+1}} =~ 
    & \Pro{RUL \leq l + \Delta t} \cdot \bigg( 
    \os{\ctone}{1}{l} + 
    \os{\cteight}{8}{\Delta t} - 
    \os{\cttwo}{2}{\EE{RUL \mid RUL \leq l + \Delta t}}
    \bigg) \cdot \cunav
    + \\[0.35em]
    & \Pro{RUL \leq l + \Delta t} \cdot
    \bigg(
    \os{\ctthree}{3}{r} - 
    \os{\ctfour}{4}{l} - 
    \os{\cteight}{8}{\Delta t}
    \bigg) \cdot c_{\infty} + \\[0.35em]
    & \Pro{RUL > l + \Delta t} \cdot \bigg( 
    \os{\ctfive}{5}{\EE{RUL \mid RUL > l + \Delta t}} - 
    \os{\ctsix}{6}{l} - 
    \os{\ctnine}{9}{\Delta t}
    \bigg) \cdot \cinv + \\[0.35em]
    & \Pro{RUL > l + \Delta t} \bigg(
    \os{\ctthree}{3}{r} - 
    \os{\ctfour}{4}{\EE{RUL \mid RUL > l + \Delta t}}
    \bigg) \cdot c_{\infty}.
\end{split}
\end{align}
Expanding all the terms and comparing them against their counterparts leads to \Cref{tab:oh4_cost_comparison_and_simplification}. Collecting all the terms on the left and right side yields the condition:
\begin{equation}
\begin{aligned}
    0 < ~ & \Pro{l < RUL \leq l + \Delta t} \cdot 
    \left( 
    \cunav + \cinv \right) \cdot \left(l - \EE{RUL \mid l< RUL \leq l + \Delta t} \right) + \\
    & \Pro{l < RUL \leq l + \Delta t} \cdot 
    c_{\infty} \cdot \left(
    \EE{RUL \mid l< RUL \leq l + \Delta t} - l 
    \right) +
    \\
    &
    \Pro{RUL \leq l + \Delta t} \cdot \Delta t \cdot \left(\cunav + \cinv - c_{\infty} \right) - \cinv \cdot \Delta t \\[1em]
    \longleftrightarrow \quad 0 < ~ & 
    \Pro{l < RUL \leq l + \Delta t} \cdot \left( 
    \cunav + \cinv \right) \cdot \left(l + \Delta t - \EE{RUL \mid l< RUL \leq l + \Delta t} \right) + \notag \\
    & \Pro{l < RUL \leq l + \Delta t} \cdot c_{\infty} \cdot \left(\EE{RUL \mid l< RUL \leq l + \Delta t} - l - \Delta t \right) + \notag \\ 
    &\Pro{RUL \leq l} \cdot \Delta t \cdot \left(\cunav + \cinv - c_{\infty} \right) - \cinv \cdot \Delta t .
\end{aligned}
\end{equation}
The terms can finally be simplified to the heuristic ordering condition
\begin{equation}
\begin{aligned}
    0 < ~ & 
    \Pro{l < RUL \leq l + \Delta t} \cdot \left( 
    \cunav + \cinv  - c_{\infty} \right) \cdot \left(l + \Delta t - \EE{RUL \mid l< RUL \leq l + \Delta t} \right) + \\
    &\Pro{RUL \leq l} \cdot \Delta t \cdot \left(\cunav + \cinv - c_{\infty} \right) - \cinv \cdot \Delta t.
\end{aligned}
\end{equation}
%
%\newpage
%
\begin{landscape}
{
\renewcommand{\arraystretch}{1}
\newlength{\mytableskip}
\setlength{\mytableskip}{10pt}
    \begin{tabular}{r@{}c|@{}c}
    %\linespread{1}\selectfont
    & Expected costs for ordering at $t_k$ & Expected costs for ordering at $t_{k+1}$ \\
    \toprule
    % l * c_unav
    \textcolor{\ctone}{1} & $\cancel{\Pro{RUL \leq l} \cdot l \cdot \cunav}$ &
    $\Pro{RUL \leq l + \Delta t} \cdot l \cdot \cunav = %$ \\
    %& & 
    %$
    \left[\cancel{\Pro{RUL \leq l}} + \Pro{l < RUL \leq l + \Delta t}\right] \cdot l \cdot \cunav$ \\[\mytableskip]
    & $\longrightarrow \quad 0$ & $\Pro{l < RUL \leq l + \Delta t} \cdot l \cdot \cunav$ \\
    \midrule
    %  l * c_inv
    \textcolor{\ctsix}{6} & $ \cancel{- \Pro{RUL > l} \cdot l \cdot \cinv}$ & $- \Pro{RUL > l + \Delta t} \cdot l \cdot \cinv = 
    %$ \\
    %& & 
    %$
    - \left[\cancel{\Pro{RUL > l}} - \Pro{l < RUL \leq l + \Delta t} \right] \cdot l \cdot \cinv$ \\[\mytableskip]
    & $\longrightarrow \quad 0$ & $\Pro{l < RUL \leq l + \Delta t} \cdot l \cdot \cinv$ \\
    \midrule
    % \EE{RUL | x} * c_unav
    \textcolor{\cttwo}{2} & $\cancel{-\Pro{RUL \leq l} \cdot \EE{RUL \mid RUL \leq l} \cdot \cunav}$ &
    $ - \Pro{RUL \leq l + \Delta t} \cdot \EE{RUL \mid RUL \leq l + \Delta t} \cdot \cunav = 
    %$ \\
    %& &
    %$
    - \int_0^{l + \Delta t} \frulr \cdot r \dr \cdot \cunav$ \\
    & & $= -\left[\cancel{\int_0^{l} \frulr \cdot r \dr} + \int_l^{l + \Delta t} \frulr \cdot r \dr \right] \cdot \cunav$ \\[\mytableskip]
    & $ \longrightarrow \quad 0$ & 
    $- \Pro{l < RUL \leq l + \Delta t} \cdot \EE{RUL \mid l< RUL \leq l + \Delta t} \cdot \cunav$ \\
    \midrule
    % \EE{RUL | x} * c_inv
    \textcolor{\ctfive}{5} & $\cancel{\Pro{RUL > l} \cdot \EE{RUL \mid RUL > l} \cdot \cinv}$ & 
    $ \Pro{RUL > l + \Delta t} \cdot \EE{RUL \mid RUL > l + \Delta t} \cdot \cinv =
    %$ \\
    %& & 
    %$
    \int_{l + \Delta t}^{\infty} \frulr \cdot r \dr \cdot \cinv $ \\
    & & $=\left[\cancel{\int_l^{\infty} \frulr \cdot r \dr} - \int_l^{l + \Delta t} \frulr \cdot r \dr \right] \cdot \cinv$ \\[\mytableskip]
    & $\longrightarrow \quad 0$ & 
    $- \Pro{l < RUL \leq l + \Delta t} \cdot \EE{RUL \mid l < RUL \leq l + \Delta t} \cdot \cinv$ \\
    \midrule
    % c_{\infty} * r
    \textcolor{\ctthree}{3} & $\left[\Pro{RUL \leq l} + \Pro{RUL > l}\right] \cdot c_{\infty} \cdot r $ & 
    $\left[\Pro{RUL \leq l + \Delta t} + \Pro{RUL > l + \Delta t}\right] \cdot c_{\infty} \cdot r $ \\
    & $=\cancel{c_{\infty} \cdot r}$ & $=\cancel{c_{\infty} \cdot r}$ \\[\mytableskip]
    & $\longrightarrow \quad 0$ & 0 \\
    \midrule
    % c_{\infty} * l
    \textcolor{\ctfour}{4} & $- \Pro{RUL \leq l} \cdot c_{\infty} \cdot l$ & 
    $- \Pro{RUL \leq l + \Delta t} \cdot c_{\infty} \cdot l = %$ \\
    %& & $
    - \left[\cancel{\Pro{RUL \leq l}} + \Pro{l < RUL \leq l + \Delta t} \right] \cdot c_{\infty} \cdot l$ \\[\mytableskip]
    & $\longrightarrow \quad 0$ & $-\Pro{l < RUL \leq l + \Delta t} \cdot c_{\infty} \cdot l$ \\
    \midrule
    % \EE{RUL | x} * c_{\infty}
    \textcolor{\ctseven}{7} & $-\Pro{RUL > l} \cdot \EE{RUL \mid RUL > l} \cdot c_{\infty}$ & 
    $-\Pro{RUL > l + \Delta t} \cdot \EE{RUL \mid RUL > l + \Delta t} \cdot c_{\infty} = $ \\
    & $= \cancel{- \int_l^{\infty} \frulr \cdot r \dr} \cdot c_{\infty} $ & 
    $- \left[\cancel{\int_l^{\infty} \frulr \cdot r \dr} - 
    \int_l^{l+\Delta t} \frulr \cdot r \dr \right] \cdot c_{\infty}$     
    \\[\mytableskip]
    & $\longrightarrow \quad 0$ & 
    $\Pro{l < RUL \leq l + \Delta t} \cdot \EE{RUL \mid l< RUL \leq l + \Delta t} \cdot c_{\infty}$ \\
    \midrule 
    % DT * c_unav - c_{\infty}
    \textcolor{\cteight}{8} & & $\Pro{RUL \leq l + \Delta t} \cdot \Delta t \cdot \left(\cunav - c_{\infty} \right)$ \\ 
    \midrule
    % DT * c_inv
    \textcolor{\ctnine}{9} & & $-\Pro{RUL > l + \Delta t} \cdot \Delta t \cdot \cinv$ \\
    \multicolumn{3}{c}{} \\
    \end{tabular}
    \captionof{table}{Cost comparisons for ordering at $t_k$ vs. $t_{k+1}$, where the respective terms from \Cref{eq:app_cord_at_tk_total,eq:app_cord_at_tk1_total} are sorted by row. The table shows the cancellation of terms and the resulting simplification of the equations.}
    \label{tab:oh4_cost_comparison_and_simplification}
}
\end{landscape}
%

%
%
%
%%%%%%%%%%%%%%%%%%%%%%%%%%%%%%%%%%%%%%%%%%%%%%%%%%%%%%%
%%%%%%%%%%%%%%%%%%%%%%%%%%%%%%%%%%%%%%%%%%%%%%%%%%%%%%%
%%%%%%%%%%%%%%%%%%%%%%%%%%%%%%%%%%%%%%%%%%%%%%%%%%%%%%%
%
%

\subsection{Derivation of initial cost rate estimate option 2 for doa ordering heuristic}
\label{sec:app_derivation_oh4_ectr_option2}
Suppose $t_k<T_F<t_{k+1}$ for a certain component. The decision that the operator now faces is whether they should order earlier and pay the inventory cost, or order later and pay the unavailability cost. Since the goal is to minimize the cost, the operator should choose to pay for the inventory if
\begin{equation}
\label{eq:app_inv_unav_ini_decision}
    \cinv \cdot \left(T_F - t_k\right) < \cunav \cdot \left(t_{k+1} - T_F \right).
\end{equation}
By inserting $t_{k+1}=t_k+\Delta t$ and reordering, the operator should pay the inventory cost if
\begin{equation}
    \label{eq:app_inv_unav_tf_decision}
    T_F < t_k + \frac{\cunav}{\cinv+\cunav} \cdot \Delta t.
\end{equation}
Hence, the decision on which cost to pay boils down to how close the component fails to $t_k$ relative to $t_{k+1}$. We are therefore only interested in the remaining failure time\footnote{This definition is identical to the RUL in the prognostic setting. However, to avoid confusion, we use a different name for the case of no monitoring.}
\begin{equation}
    \label{eq:app_rfi_def}
    R_F=T_F - t_k,
\end{equation}
which, when inserted into \Cref{eq:app_inv_unav_tf_decision}, yields 
\begin{equation}
    \label{eq:app_inv_unav_rf_decision}   
    R_F < \frac{\cunav}{\cinv+\cunav} \cdot \Delta t.
\end{equation}
Therefore, the expected costs can be calculated as an expectation over $R_F$:
\begin{equation}
    \label{eq:app_min_ord_cost_int}   
    \EE{C_{\mathrm{min}}} = \int_0^{\Delta t} C(t)f_{R_F}(t) \dt,
\end{equation}
where $C(t)$ defines the minimal cost function (i.e., whether we pay inventory or unavailability costs). For performing the computation in \Cref{eq:app_min_ord_cost_int}, one either needs the exact distribution for the remaining failure time, or one can resort to a sample-based computation.

We can approximate $R_F \sim \mathcal{U}[0,\Delta t]$. We checked this approximation for the RUL simulator, and it is quite accurate for $\Delta t$ sufficiently smaller than $\mutf$ \cite{kuipers2012uniform}. With this, the expected costs can be calculated as:
\begin{equation}
\begin{aligned}
    \EE{C_{\mathrm{min}}} &= \int_0^{\Delta t} C(t) \frac{1}{\Delta t} \dt
    \\ 
    & = 
    \int_0^{\frac{\cunav}{\cinv + \cunav}\Delta t} t\cdot \cinv \frac{1}{\Delta t} \dt + \int_{\frac{\cunav}{\cinv + \cunav}\Delta t}^{\Delta t} (\Delta t - t) \cdot \cunav \frac{1}{\Delta t} \dt 
    \\
    & = 
    \int_0^{\frac{\cunav}{\cinv + \cunav}\Delta t} t\cdot \cinv \frac{1}{\Delta t} \dt + \int_{0}^{\frac{\cinv}{\cinv + \cunav}\Delta t} t \cdot \cunav \frac{1}{\Delta t} \dt
    \\
    & = 
    \frac{\cinv}{2\Delta t} \cdot \frac{\cunav^2}{(\cinv+\cunav)^2} \Delta t^2 + \frac{\cunav}{2\Delta t} \cdot \frac{\cinv^2}{(\cinv+\cunav)^2} \Delta t^2 \\
    & = 
    \frac{\cinv \cdot \cunav \Delta t}{2(\cinv+\cunav)^2}(\cinv + \cunav) 
    \\
    &= 
    \frac{\cinv \cdot \cunav}{2(\cinv+\cunav)} \Delta t.
\end{aligned}
\end{equation}
Finally, we obtain an analytical expression for an initial estimate of the policy cost rate under perfect ordering:
\begin{equation}
    c_{\infty,0} =  \frac{\cinv \cdot \cunav}{2(\cinv+\cunav)} \frac{\Delta t}{\mutf} .
\end{equation}
%

%
%
%
%%%%%%%%%%%%%%%%%%%%%%%%%%%%%%%%%%%%%%%%%%%%%%%%%%%%%%%
%%%%%%%%%%%%%%%%%%%%%%%%%%%%%%%%%%%%%%%%%%%%%%%%%%%%%%%
%%%%%%%%%%%%%%%%%%%%%%%%%%%%%%%%%%%%%%%%%%%%%%%%%%%%%%%
%
%
%
\section{Virtual RUL simulator}
\label{sec:virtual_RUL_simulator}
This section is more or less an iteration of the RUL simulator description in \cite{kamariotis2024metric}.

We assume an underlying normal distribution of the time-to-failure of the components, i.e., $\Tf \sim\mathcal{N}(\mu=225,\sigma=40)$. At each step $t_k$, we make a prediction about the remaining useful life of the component, denoted as $RUL_{k}$. This prediction deviates from the true RUL, denoted $RUL_{\mathrm{true},k}=\Tf-t_k$, by a certain error $\epsilon_k$. We define this prediction in a logarithmic scale as:
\begin{equation}
    \label{eq:rul_pred_error}
    \ln\left(RUL_{k} \right) = \ln\left(RUL_{\mathrm{true},k}\right) + \ln\left( \epsilon_k \right).
\end{equation}
In this work, we assume a multivariate normal (MVN) distribution for the logarithmic prediction errors over time. Specifically, 
\begin{equation}
    \label{eq:mvn_of_errors}
    [\ln(\epsilon_1), \dots, \ln(\epsilon_n)] \sim \mathrm{MVN}(\mathbf{0}, \mathbf{\Sigma}).
\end{equation}
The mean of the prediction errors is 0, hence the prognostic model is designed to be unbiased in log-scale. On the other hand, the covariance matrix $\mathbf{\Sigma}$ (not to be confused with the strategy $\Sigma$ in \Cref{sec:pdm_policies}!) is constructed as
\begin{equation}
    \label{eq:cov_mat}
    \mathbf{\Sigma} = \mathbf{D} \cdot \mathbf{R} \cdot \mathbf{D},
\end{equation}
where $\mathbf{D}$ is a diagonal matrix containing the standard deviation of the prediction errors $\sigma_{\ln(\epsilon_k)}$ on the main diagonal. For this work, we fix $\sigma_{\ln(\epsilon_k)}=0.4~\forall k$. $\mathbf{R}$ is a correlation matrix. If the predictions of the prognostic model are assumed to be independent at $t_k$, then $\mathbf{R}$ reduces simply to a diagonal matrix with $1$ on the main diagonal. In this work, the prognostic model's predictions are correlated as
\begin{equation}
    \label{eq:corr_mat}
    \mathbf{R} = [\rho_{ij}], \qquad \rho_{ij} = \exp \left(-\frac{|t_i-t_j|}{l} \right),
\end{equation}
where $l$ is the hyperparameter defining the correlation length.

The steps to generate the RUL-PDFs at each time steps are then:
\begin{enumerate}
    \item Draw a failure sample $\Tfi$ from the underlying time to failure distribution.
    \item Compute the true RUL values at each time step via $RUL_{\mathrm{true},k}^{(i)}=\Tfi-t_k$.
    \item Sample an error vector $[\ln(\epsilon_1^{(i)}), \dots, \ln(\epsilon_n^{(i)})]$ from the MVN distribution in \Cref{eq:mvn_of_errors}.
    \item Compute the mean value of the prediction via \Cref{eq:rul_pred_error}: $\mu_{\ln\left(RUL_{k} \right)}^{(i)}=\ln(RUL_{\mathrm{true},k}^{(i)}) + \ln(\epsilon_k^{(i)})$.
    \item The RUL-PDF is then a lognormal distribution $\ln\left(RUL_{k}^{(i)} \right) \sim \mathcal{N}\left(\mu_{\ln\left(RUL_{k} \right)}^{(i)}, \sigma_{\ln(\epsilon_k)}\right)$.
\end{enumerate}
%

%
%
%
%%%%%%%%%%%%%%%%%%%%%%%%%%%%%%%%%%%%%%%%%%%%%%%%%%%%%%%
%%%%%%%%%%%%%%%%%%%%%%%%%%%%%%%%%%%%%%%%%%%%%%%%%%%%%%%
%%%%%%%%%%%%%%%%%%%%%%%%%%%%%%%%%%%%%%%%%%%%%%%%%%%%%%%
%
%
%
\section{Neural Network specifications}
\label{sec:app_nn_specs}
The NN architecture used in this work is depicted in \Cref{fig:NN_architecture}. It includes at least one LSTM layer, at least one dropout layer, one fully connected (FC) layer and a Gaussian error linear unit (GELU) activation function \cite{hendrycks2016gaussian}. The 25-dimensional input consists of 21 sensor measurements, 3 operation modes, and 1 variable indicating the age of the engine, where each input is passed as a window of length 50 to the network. 

The fixed setting choices for all runs are given in \Cref{tab_app:fixed_nn_params}, whereas the hyperparameters optimized with a hyperparameter search (HS) for each Monte Carlo cross-validation split are given in \Cref{tab_app:nn_hyperparams}. 
\nnarchitecture
\begin{table}[H]
    \renewcommand{\arraystretch}{0.9}
    \setlength{\abovecaptionskip}{10pt}
    \centering
    \begin{tabular}{ll}
        \toprule
        Setting & Choice \\
        \midrule
        Optimizer & ADAM \cite{kingma2014adam} \\
        AMSGRAD \cite{reddi2019amsgrad} & Included \\
        \# HS (\# MCCV splits) & 100 \\
        \# runs per HS & 20 \\
        HS type & Bayesian optimization \cite{wandbbayes} \\
        LR scheduler type & Reduce on plateau \\
        Input window size & 50 \\
        \bottomrule
    \end{tabular}
    \caption{Fixed training setting choices for all hyperparameter searches.}
    \label{tab_app:fixed_nn_params}
\end{table}
\begin{table}[H]
    \setlength{\abovecaptionskip}{10pt}
    \centering
    \begin{tabular}{ll}
        \toprule
        Parameter & Distribution \\
        \midrule
        Batch size & $\{64, 128, 256\}$ \\
        LSTM hidden size & $\{64, 128\}$ \\
        \# LSTM blocks & $\{1, 2\}$ \\ 
        Dropout \cite{srivastava2014dropout} probability & $\mathcal{U}[0.0, 0.3]$
        \\
        Initial lr & $\mathcal{LU}[10^{-5}, 10^{-1}]$ \\
        LR scheduler factor & $\mathcal{L}[0.5, 0.9]$ \\
        LR scheduler patience & $\{3, 4, ..., 25\}$ \\
        \bottomrule
    \end{tabular}
    \caption{Distributions for each hyperparameter included in the search, where $\{ \cdot \}$, $\mathcal{U}[\cdot]$, and $\mathcal{LU}[\cdot]$ represent the discrete uniform, uniform, and log-uniform distributions, respectively.}
    \label{tab_app:nn_hyperparams}
\end{table}
%

%
%
%
%%%%%%%%%%%%%%%%%%%%%%%%%%%%%%%%%%%%%%%%%%%%%%%%%%%%%%%
%%%%%%%%%%%%%%%%%%%%%%%%%%%%%%%%%%%%%%%%%%%%%%%%%%%%%%%
%%%%%%%%%%%%%%%%%%%%%%%%%%%%%%%%%%%%%%%%%%%%%%%%%%%%%%%
%
%
%
\section{Further investigations of rh2}
\label{sec:app_further_investigations_rh2}
%

%
%
%%%%%%%%%%%%%%%%%%%%%%%%%%%%%%%%%%%%%%%%%%%%%%%%%%%%%%%
%%%%%%%%%%%%%%%%%%%%%%%%%%%%%%%%%%%%%%%%%%%%%%%%%%%%%%%
%
%
\subsection{Derivation}
\label{subsec:app_derivation_rh2}
Here, we provide a derivation for $rh2$. For ease of notation, we again drop the subscript $k$ wherever possible.

The overall goal is to find a maintenance policy with minimal $c_{\infty}$. It is easier to work directly with costs instead of rates; hence, we fix the time frame, and then we can perform optimization over the costs. Consider time step $t_k$, where one decides whether to replace the component at some future time $t_k+\rrep$ or not. Since we base our cost calculations on the latest RUL-PDF, $f_{RUL}(r)$, the natural reference time frame is chosen as the largest residual time in the support of the distribution, i.e., $\rmax = \max~ r \mid  f_{RUL}(r) > 0$. 

%\vspace{0.5em} \noindent
%Case 1: \newline 
%
\begin{enumerate}
    \item[]\emph{Case 1:} The component fails before replacement, i.e., $RUL \in [0, ~\rrep]$. 
    
    This incurs a corrective replacement cost at $RUL$ and the long-running maintenance cost rate for the remaining time until $\rmax$:
\begin{equation}
    \label{eq:rh2_cc}
    C_c = c_c + c_{\infty} \cdot (\rmax - RUL).
\end{equation}
Using the predicted RUL-PDF, we can calculate the expectation over $C_c$ in the timeframe $[0,~ \rrep]$ as
\begin{align}
\label{eq:Exp_C_c_Tmax}
    \EE{C_c} &= \int_{0}^{\rrep} \left[c_c + c_{\infty} \cdot (\rmax - r)\right] \cdot f_{RUL}(r) \dr \notag \\
    & =  c_c\cdot p_F(\rrep) + c_{\infty} \cdot \int_{0}^{\rrep} (\rmax - r) \cdot f_{RUL}(r) \dr,
\end{align}
where $p_F(\rrep)=\int_{0}^{\rrep} f_{RUL}(r) \dr$.

\vspace{0.5em} \noindent
\item[]\emph{Case 2:} The component fails after the planned replacement, i.e., $RUL \in [\rrep, ~\rmax]$. 

This incurs a preventive replacement cost at $\rrep$ and the long-running maintenance cost rate for the time frame $\rmax-\rrep$:
\begin{equation}
    \label{eq:rh2_cp}
    C_p = c_p + c_{\infty} \cdot (\rmax - \rrep).
\end{equation}
Analogously to \Cref{eq:Exp_C_c_Tmax}, we can calculate the expectation over $C_p$ in the timeframe $[\rrep,~ \rmax]$:
\begin{equation}
\begin{aligned}
\label{eq:Exp_C_p_Tmax}
    \mathrm{E}\left[C_p\right] &= \int_{\rrep}^{\rmax} \left[c_p + c_{\infty} \cdot (\rmax - \rrep)\right] \cdot  f_{RUL}(r) \dr \\
    &= \int_{\rrep}^{\rmax} c_{\infty} \cdot (\rmax - r + r - \rrep) \cdot  f_{RUL}(r) \dr 
    + c_p \cdot (1-p_F(\rrep))\\
    & = c_{\infty} \cdot \left[\int_{\rrep}^{\rmax} (\rmax - r) \cdot  f_{RUL}(r) \dr  + 
    \int_{\rrep}^{\rmax} (r - \rrep) \cdot  f_{RUL}(r) \dr \right] \\
    & + c_p \cdot (1-p_F(\rrep)).
\end{aligned}
\end{equation}
\end{enumerate}
The total expected cost of performing a preventive maintenance action at $\rrep$ is then the sum of expected corrected and preventive maintenance costs in  \Cref{eq:Exp_C_c_Tmax,eq:Exp_C_p_Tmax}, respectively:
\begin{align}
\label{eq:Crepfull}
    \EE{C_{PR@t_k+\rrep}} 
    & = c_c\cdot p_F(\rrep) + c_p \cdot (1-p_F(\rrep)) + c_{\infty} \cdot \left[\int_{0}^{\rrep} (\rmax - r) \cdot f_{RUL}(r) \dr \right. \notag \\
    & \left. + \int_{\rrep}^{\rmax} (\rmax - r) \cdot f_{RUL}(r) \dr  + 
    \int_{\rrep}^{\rmax} (r - \rrep) \cdot f_{RUL}(r) \dr \right] \notag \\
    & = c_c\cdot p_F(\rrep) + c_p \cdot (1-p_F(\rrep)) + c_{\infty} \cdot \left[\int_{0}^{\rmax} (\rmax - r) \cdot f_{RUL}(r) \dr  \right. \notag \\
    & \left.  + \int_{\rrep}^{\rmax} (r - \rrep) \cdot f_{RUL}(r) \dr \right].
\end{align}
Since $\int_{0}^{\rmax} (\rmax - r)\cdot f_{RUL}(r) \dr$ is independent of $\rrep$, the objective function can be split into a part depending on $\rrep$ and a constant part:
\begin{equation}
    \EE{C_{PR@t_k+\rrep}} = \EE{C_{PR@t_k+\rrep}(\rrep)} + \EE{C_{PR@t_k+\rrep}^{\mathrm{const}}}.
\end{equation}
For finding the optimal replacement time $\rrep^*$, the constant term is not of importance and thus can be dropped. The resulting formulation is identical to \Cref{eq:rep_optimization_obj_fun} when letting $\rmax \rightarrow \infty$, with the corresponding decision heuristic defined in \Cref{eq:rep_optimization}.

%
%
%%%%%%%%%%%%%%%%%%%%%%%%%%%%%%%%%%%%%%%%%%%%%%%%%%%%%%%%%%%%%%%%%%%%%%%%%%%%%%%%
%
%
\subsection{Policy Extensions}
\label{subsec:Extension_h3}
For the use of $rh2$, one needs $c_{\infty}$ to define the objective function. The accumulation of individual decisions then makes up the final maintenance policies, whose performance is again $c_{\infty}$. Thus, to compute the optimal replacement time, one would already need to know the final performance of this policy. To circumvent this implicit definition, \citet{kamariotis2024metric} propose three initial estimates. These estimates most likely do \emph{not} coincide with the policy's final performance.

Hence, to align the assumed cost rate value with the actual performance, i.e., $\hat{C}_{\infty}=c_{\infty,0}$, we also investigate an iterative procedure as well as a direct optimization, which are outlined in \Cref{subsubsec:Iterative_procedure,subsubsec:Direct_optimization}.

%
%
%%%%%%%%%%%%%%%%%%%%%%%%%%%%%%%%%%%%%%%%%%%%%%%%%%%%%%%%%%%%%%%%%%%%%%%%%%%%%%%%
%
%
\subsubsection{Iterative procedure}
\label{subsubsec:Iterative_procedure}
The idea for the iterative procedure is to start with one of the $c_{\infty,0}$ options, use it for the objective function at each $t_k$, compute the resulting performance of the policy, $\hat{C}_{\infty}$, which is then used as the improved estimate $c_{\infty,1}$, and so on. Theoretically, if $rh2$ really represents the optimal policy, and the function $\hat{C}_{\infty}(c_{\infty,0})$ is convex/concave, then this fixed-point iteration scheme \cite{hoffman2018numerical} should converge to the globally optimal solution, at which $\hat{C}_{\infty}$, and consequently also $\hat{M}$, is minimal.

%
%
%%%%%%%%%%%%%%%%%%%%%%%%%%%%%%%%%%%%%%%%%%%%%%%%%%%%%%%%%%%%%%%%%%%%%%%%%%%%%%%%
%
%
\subsubsection{Direct optimization}
\label{subsubsec:Direct_optimization}
Another data-driven option is to treat $c_{\infty}$ as a mere parameter to be optimized. Similarly to $p_{\mathrm{thres}}$ of $rh1$, one solves an optimization problem to find the optimal value $c_{\infty}^*$, for which $\hat{C}_{\infty}$ is minimal. If this value resulting from the optimization is substantially different from the converged value of the procedure in \Cref{subsubsec:Iterative_procedure}, then $rh2$ does not accurately describe the costs arising from the underlying renewal-reward process, and there is room for improvement. 
%
%
%%%%%%%%%%%%%%%%%%%%%%%%%%%%%%%%%%%%%%%%%%%%%%%%%%%%%%%%%%%%%%%%%%%%%%%%%%%%%%%%
%
%
\subsubsection{Results}
\label{subsubsec:app_results}
\begin{figure}[H]
    \centering
    \includegraphics[width=0.9\textwidth]{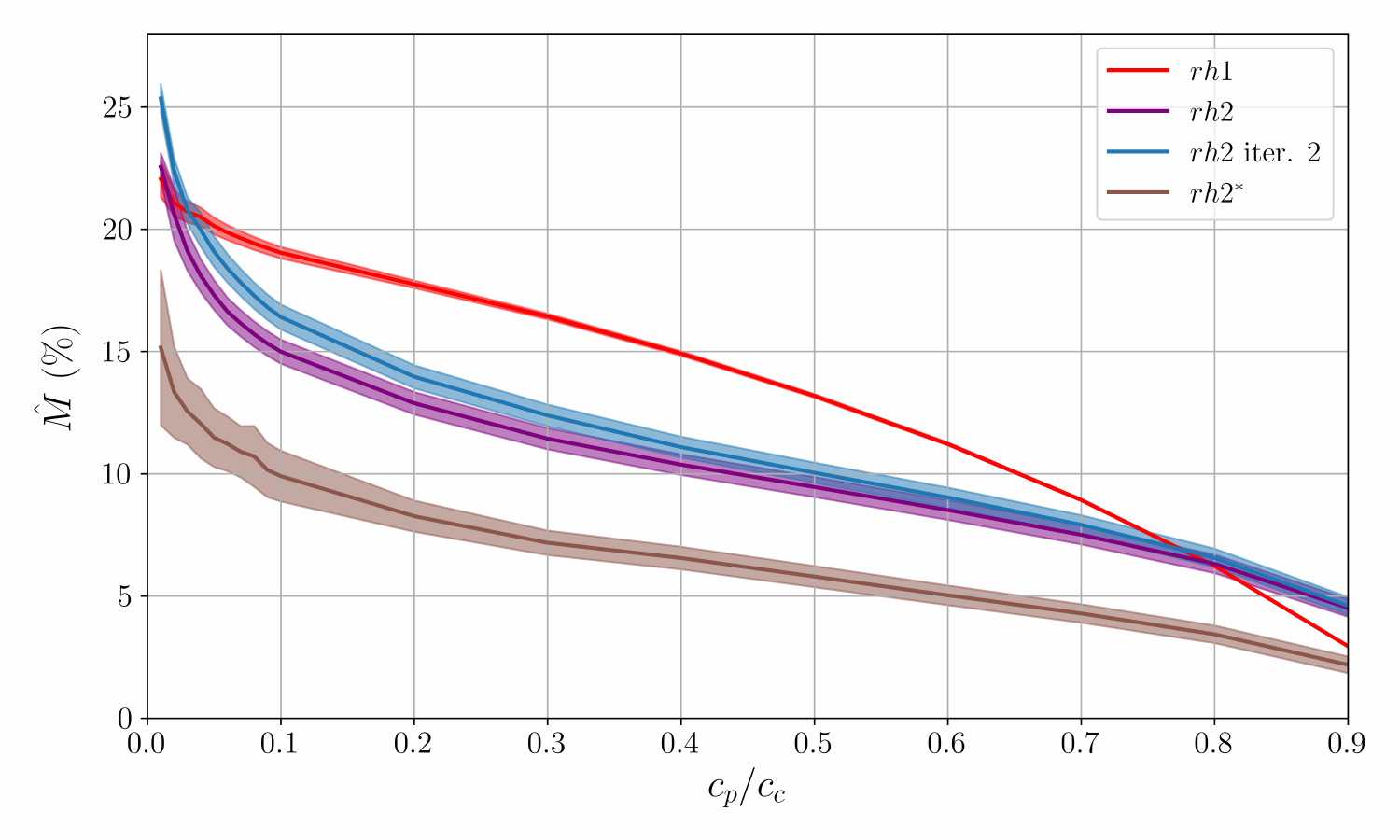}
    \caption{Comparison of the default $rh2$ with its iterative and direct optimization variants. $rh1$ with the default probability threshold is also shown for reference. For all $rh2$ variants, $1.2\cdot10^4$ MC samples were used for training/evaluation. The figure shows the mean of each policy as a solid line, as well as the corresponding 95\% credible intervals.}
    \label{fig:rh_perf_rh2_overview}
\end{figure}

\begin{figure}[H]
    \centering
     \begin{tikzpicture}
        \node[anchor=south west,inner sep=0] (image) at (0,0) {\includegraphics[width=0.83\textwidth]{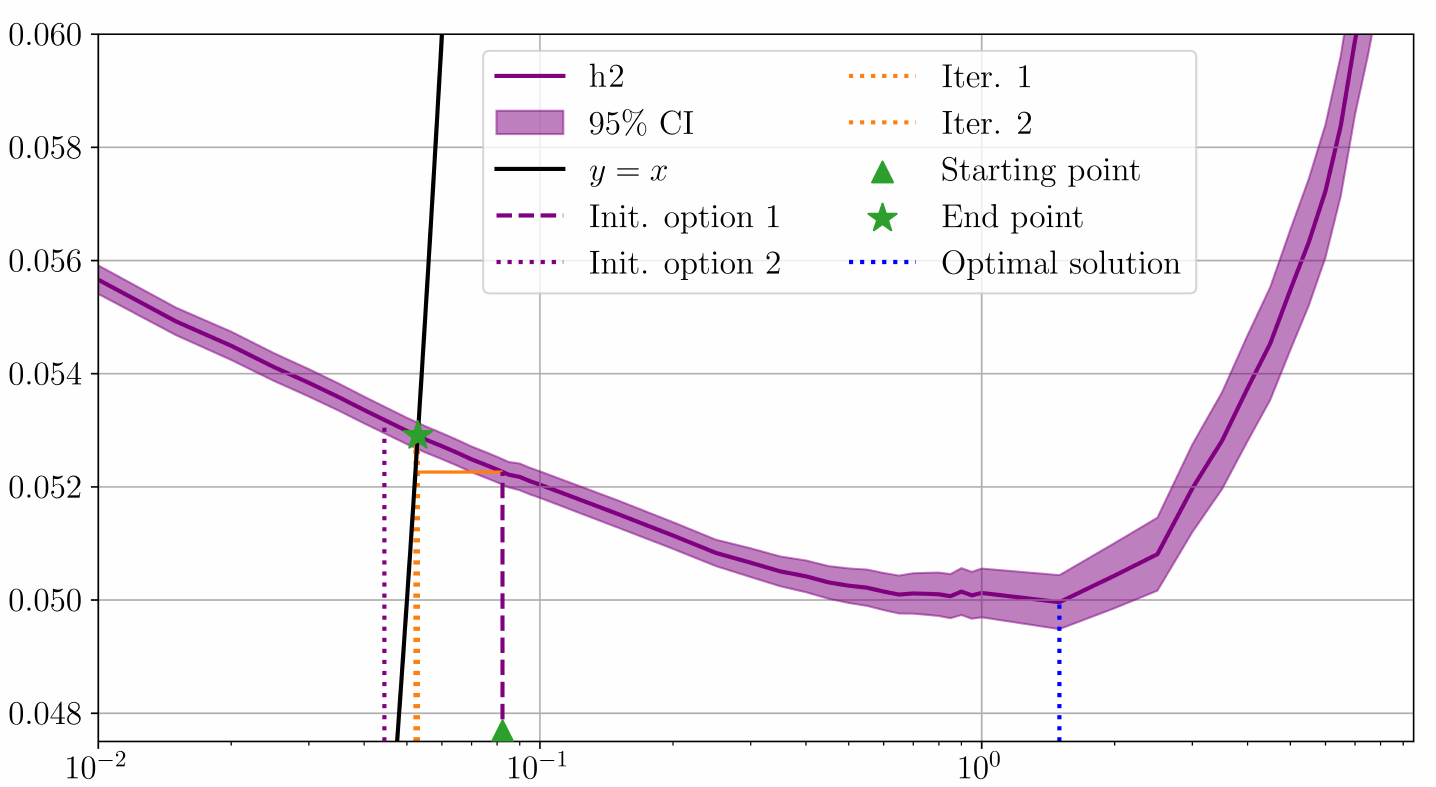}};
        \begin{scope}[x={(image.south east)}, y={(image.north west)}]
            \node[align=center] at (0.5,-0.03) {\Large $c_{\infty,0}$};
            \node[align=center] at (-0.03,0.5) {\Large $\hat{C}_{\infty}$};
        \end{scope}
    \end{tikzpicture}
    \caption{$rh2$ performance with $\frac{c_p}{c_c} = 0.1$ plotted over a range of $c_{\infty,0}$ values, including two iterations as described in \Cref{subsubsec:Iterative_procedure}, initialization options 1 \& 2 from \cite{kamariotis2024metric}, as well as the optimal value. The two iterations lie very close to each other, which is why they are visually overlapping to a single line.}
    \label{fig:rh2_iteration_problem}
\end{figure}

The results for the comparison of the standard $rh2$ approach with the proposed extensions are depicted in \Cref{fig:rh_perf_rh2_overview}. 
Surprisingly, the iteration scheme leads to a worse performance than by using just the default $c_{\infty,0}$ option 1. The reason for this is shown in \Cref{fig:rh2_iteration_problem}. 
Although the iteration procedure converges to a value somewhere in between the bounds specified by initialization options 1 \& 2, the optimal value, $c_{\infty}^*$, for which the final policy performance $\hat{C}_{\infty}$ is minimal, is much larger than the upper bound. Moreover, the performance curve $rh2^*$ overlaps with the optimized policies $rh1^*/doa^*$ in \Cref{fig:rh_perf_overview}, i.e., one can again compensate a ``bad'' heuristic via data-driven optimization.

\end{document}